\documentclass{l4dc2027}

\title[Semi-global Convergence of MPPI]{Semi-global Convergence and Finite-Sample Error Bounds for Model Predictive Path Integral Control}
\usepackage{times}

\newcommand{\R}{\mathbb{R}}

\newcommand{\Rp}{\R_{>0}}

\newcommand{\Nnn}{\mathbb{N}_0}
\newcommand{\Np}{\mathbb{N}_+}

\newcommand{\Ident}{I} 

\newcommand{\ball}{\mathrm{B}} 
\newcommand{\trace}{\mathrm{tr}}

\newcommand{\Pb}{\mathbb{P}}
\newcommand{\E}{\mathbb{E}}
\newcommand{\Var}{\mathrm{Var}}

\newcommand{\Pbs}{\mathcal{P}} 

\newcommand{\distr}{\sim}
\newcommand{\iiddistr}{\overset{\text{i.i.d.}}{\sim}}
\newcommand{\Norm}{\mathcal{N}}

\newtheorem{assumption}{Assumption}

\author{%
 \Name{Christian Fiedler} \Email{christian.fiedler@tum.de} \\
 \Name{Giacomo Borghi} \Email{giacomo.borghi@tum.de}
 \\
 \addr Department of Mathematics, Technical University of Munich (TUM) \\and Munich Center for Machine Learning (MCML)
}

\begin{document}

\maketitle

\begin{abstract}%
  We investigate the global optimisation behaviour of Model Predictive Path Integral Control (MPPI) in the low-temperature setting.
  For objective functions with a unique global minimum and under weak additional assumptions -- allowing non-convexity and non-smoothness --  we establish semi-global convergence of the mean-field or population level MPPI iteration to this global optimum.
  We also establish finite-sample results that show that with high probability the MPPI iteration comes close to the global optimum and provide sample complexity results.
  Under slightly stronger assumptions, we make these results more quantitative and provide concrete bounds on the inverse temperature and number of samples.
  Finally, we illustrate the results on the classic LQR problem.
\end{abstract}


\section{Introduction}\label{sec:intro}
Model Predictive Path-Integral Control (MPPI) is one of the major sampling-based control methods used in robotics, autonomous driving, and model-based reinforcement learning \citep{kazim2024recent}.
While very successful in applications, the theoretical understanding of MPPI is still limited.
Introduced originally by information-theoretic arguments connected to stochastic optimal control \citep{williams2018information},
it has since been recognised that MPPI is essentially a zeroth-order sampling-based optimisation method combined with receding-horizon control.
In particular, a rigorous theoretical analysis must start with the optimisation behaviour of MPPI.
Most existing works have been focused on local perspectives, e.g., interpreting MPPI as a form of gradient descent \citep{fazlyab2026model},
or considering restricted classes of objectives like quadratic functions \citep{yi2024covo}.
This is a severe limitation, since MPPI is particularly interesting for challenging non-convex and even non-smooth problems arising in control tasks.
Motivated by this shortcoming, we investigate MPPI as a sampling-based global optimisation method for non-convex, non-smooth objective functions with a unique global optimum. 
For this, we build on the recent work \cite{fornasier2026consensus} that connects MPPI with \emph{global} non-convex optimisation.
More precisely, this work finds that the plain MPPI iteration on the population level corresponds exactly to \emph{Consensus Hopping (CH)}, a limiting case of \emph{Consensus-Based Optimisation (CBO)}, which is a class of popular global optimisation methods based on interacting particle systems \citep{pinnau2017consensus} with strong theoretical guarantees \citep{fornasier2024consensus,carrillo2018analytical}.
In turn, this connection allows one to prove  strong convergence results for CH and hence MPPI.
We build on this work and establish \emph{semi-global} convergence results for MPPI on the population level.
While already conceptually insightful, MPPI is a sampling-based method, and hence any result with practical relevance needs to take into account the finite-sample error incurred by the sampling process.
To address this important issue, we also prove corresponding high-probability error bounds and establish sample complexity bounds for the actual sampling-based MPPI iteration.
Apart from requiring a unique global optimum, our results rely only on rather general assumptions on the objective function, covering in particular non-convex and non-smooth functions.
Our work provides a significant step towards a deeper theoretical understanding of MPPI and helps to put it on a solid theoretical foundation, on which reliable sampling-based control approaches can be built.

\paragraph{Related Work} MPPI was introduced in the context of autonomous driving \citep{williams2016aggressive} based on information-theoretic considerations \citep{williams2018information} that connect to stochastic optimal control and path-integral control.
However, the connection is mostly heuristic and the rigorous theoretical foundation is lost \citep{homburger2025optimality}.
By now there is a multitude of different MPPI variants, cf.\ \cite{kazim2024recent} for an overview,
and similar to most analyses, we focus on the plain MPPI iteration formally introduced in Section \ref{sec:mppiCH}, and investigate MPPI as a zeroth-order sampling-based optimisation method, cf.\ \cite{jordana2025introduction} for more on this perspective.
In contrast to previous works like \cite{yi2024covo} or \cite{yoon2026finite}, we consider \emph{general non-convex} objective functions.
For this, we build on the recent work \cite{fornasier2026consensus} that identified MPPI on the population level as Consensus Hopping (CH), which is connected to Consensus-Based Optimisation (CBO) \citep{pinnau2017consensus}.
Utilising a quantitative Laplace principle as pioneered by the work \cite{fornasier2024consensus},
one can establish semi-global convergence guarantees for CH and hence for population-level MPPI.
The concept of semi-global results is common in nonlinear control, see \cite[Chapter~6]{gruene2017nonlinear} for more details.
The use of the Laplace principle in connection with MPPI has already been noticed in \cite{homburger2025optimality}; however, a different scaling is used there, and no non-asymptotic results are provided.
Note that it is well-known that MPPI has problems with functions having multiple global optima \citep{liu2026clustering} due to averaging effects, which motivates our focus on objectives with a unique global optimum.
We also provide finite-sample results, which are particularly relevant since in practice MPPI is always sampling-based.
Prior works have investigated this in restricted settings like LQR problems, cf.\ \citep{yoon2022sampling,yoon2026finite}, whereas we work with very general nonconvex objectives.
Instead of relying on propagation of chaos techniques as suggested in \cite{fornasier2026consensus}, we directly analyse the estimation error in a single step from the perspective of self-normalised importance weighted sampling, cf.\ \cite{agapiou2017importance} for an overview and similar results.
While this has been used in prior works on MPPI like \cite{yoon2022sampling,yoon2026finite}, our simple but effective estimation error bounds become particularly powerful in connection with our semi-global convergence results.
In particular, to the best of our knowledge we achieve the first finite-sample complexity bounds for general nonconvex objective functions with only modest regularity assumptions.

\paragraph{Contributions} \emph{First}, we isolate, generalise, and sharpen the mean-field analysis of the MPPI iteration from \cite{fornasier2026consensus}, establishing semi-global convergence to the global optimum.
\emph{Second}, we provide the first semi-global finite-sample accuracy results for MPPI.
\emph{Third}, under additional assumptions, we provide sharpened finite-sample results that include the LQR problem.

\paragraph{Outline} In Section \ref{sec:mppiCH}, we introduce the mean-field and empirical versions of MPPI that we will investigate.
Section \ref{sec:meanField} contains results for the mean-field iteration, and Section \ref{sec:finiteSample} the corresponding finite-sample results.
In the main text, we present simplified results that highlight the main features; the supplementary material contains more detailed and quantitative results with full proofs.
In Section \ref{sec:oc} we illustrate some results with the standard LQR problem.
Finally, Section \ref{sec:conclusion} contains a discussion and outlook.

\section{MPPI and Consensus Hopping}\label{sec:mppiCH}
Consider an \emph{objective function} $f: \R^d \rightarrow \R$ that is measurable and bounded from below.
Let $\alpha\in\Rp$ and let $\Sigma\in\R^{d\times d}$ be a symmetric positive definite matrix.
Throughout, $\sigma_{\min}(\Sigma)$ and $\sigma_{\max}(\Sigma)$ denote the square roots of the smallest and largest eigenvalue of $\Sigma$, respectively.
The \emph{MPPI iteration} with \emph{inverse temperature} $\alpha$ and \emph{sampling covariance matrix} $\Sigma$ is defined as follows.
Starting from an initial guess $\hat x_0 \in \R^d$, in order to compute the next iterate $\hat x_{k+1}$ from the current iterate $\hat x_k$,
we sample $X_1,\ldots,X_N \iiddistr \Norm(\hat x_k,\Sigma)$ and then set
\begin{equation} \label{eq:mppiIteration}
    \hat x_{k+1} = \frac{\sum_{i=1}^N \exp(-\alpha f(X_i))X_i}{\sum_{i=1}^N \exp(-\alpha f(X_i))}.
\end{equation}
Note that in each iteration we draw a new sample.

We start  the analysis of the MPPI iteration on the \emph{mean-field} or \emph{population-level}.
For this, define the transition map
\begin{equation} \label{eq:meanFieldTransitionMap}
\begin{aligned}
    \mathcal{T}_\alpha: \R^d\rightarrow \R^d, 
    \qquad 
    \mathcal{T}_\alpha(x) &:= \frac{\E_{X \distr \Norm(x,\Sigma)}[\exp(-\alpha f(X))X]}{\E_{X\distr\Norm(x,\Sigma)}[\exp(-\alpha f(X))]}\\
    &= \frac{1}{Z_\alpha(x)}\E_{X \distr \Norm(x,\Sigma)}[\exp(-\alpha f(X))X],
\end{aligned}
\end{equation}
where $Z_\alpha(x) = \E_{X\distr\Norm(x,\Sigma)}[\exp(-\alpha f(X))]$ is a normalisation constant,
and for later use we define the \emph{free energy} by
\begin{equation}
  \mathcal{E}_\alpha(x) = - \frac{1}{\alpha}\log Z_\alpha(x) = - \frac{1}{\alpha}\log \E_{X\distr\Norm(x,\Sigma)}[\exp(-\alpha f(X))].
\end{equation}
Given some initial state $x_0\in\R^d$, we define the \emph{mean-field MPPI iteration} by
\begin{equation}
  x_{k+1} = \mathcal{T}_\alpha(x_k), \quad k\in\Nnn.
\end{equation}
As is well-known, each step of the MPPI iteration \eqref{eq:mppiIteration} is a self-normalised importance-weighted estimate of the corresponding mean-field step.
As pointed out in \cite{fornasier2026consensus}, the mean-field MPPI iteration is exactly the Consensus Hopping (CH) iteration from this reference.

Assume now that $f$ has a global minimiser $x_\ast$.
The core question of the present work is in which sense the MPPI iterates $(\hat x_k)_k$ and the mean-field MPPI iterates $(x_k)_k$ approximate the global minimiser $x_\ast$ of $f$.
It is clear that without any additional assumptions on $f$, little progress can be made,
so our analysis will be based mostly on the following properties of $f$.
%
%
\begin{assumption}[Objective function]\label{ass:objective}
The function $f:\R^d\to\R$ satisfies the following.
\begin{enumerate}
    \item[(F1)] $f$ is measurable and attains its global minimum $f_\ast\in\R$ at $x_\ast\in\R^d$.
    \item[(F2)]  $f$ is upper semicontinuous at $x_\ast$.
    \item[(F3)]  For every $\epsilon>0$ there exists $S_\epsilon>0$ such that
    \begin{equation}\label{eq:objective:separation}
        \inf_{\|x-x_\ast\|_2\ge\epsilon} f(x)\;\ge\; f_\ast+S_\epsilon .
    \end{equation}
\end{enumerate}
\end{assumption}
(F2) is a very weak continuity assumption.
(F3) makes the notion precise that the global minimiser $x_\ast$ is sufficiently separated from any other local minimisers and that the global minimum is sufficiently pronounced.
In particular, an elementary computation shows that (F3) implies that $x_\ast$ is indeed the unique global minimiser of $f$.

To obtain explicit thresholds, we also consider the following simplified assumption, which only requires strong convexity close to the global minimiser and a suboptimality gap away from it.
%
%
\begin{assumption}[Local strong convexity and energy gap]\label{ass:simplified}
The function $f:\R^d\to\R$ is continuous and attains its global minimum $f_\ast$ at $x_\ast\in\R^d$, and there are $R,m,\Delta>0$ such that the following hold.
\begin{enumerate}
    \item[(SF1)] The function $y\mapsto f(y)-\frac m2\|y-x_\ast\|_2^2$ is convex on $\ball_R(x_\ast)$.
    \item[(SF2)] $f(y)-f_\ast\ge\Delta$ for all $y\in\R^d\setminus\ball_R(x_\ast)$.
    \item[(SF3)] There is $M\ge m$ such that $f(y)-f_\ast\le\frac M2\|y-x_\ast\|_2^2$ for all $y\in\ball_R(x_\ast)$.
\end{enumerate}
\end{assumption}
(SF1) and (SF2) imply Assumption~\ref{ass:objective}, and (F3) holds with $S_\epsilon=\min\{\frac m2\epsilon^2,\Delta\}$, cf.\ Lemma~\ref{lem:simplifiedImpliesObjective} in the supplementary material.
We stress that outside the well $\ball_R(x_\ast)$, the objective can still be non-convex and non-smooth.
Furthermore, condition (SF3) is only needed for the finite-sample results, and it holds, for instance, if $f$ has an $M$-Lipschitz gradient on $\ball_R(x_\ast)$.

\section{Semi-Global Convergence in Mean-Field}\label{sec:meanField}
\begin{figure}[t]
\centering
\includegraphics[width=\textwidth]{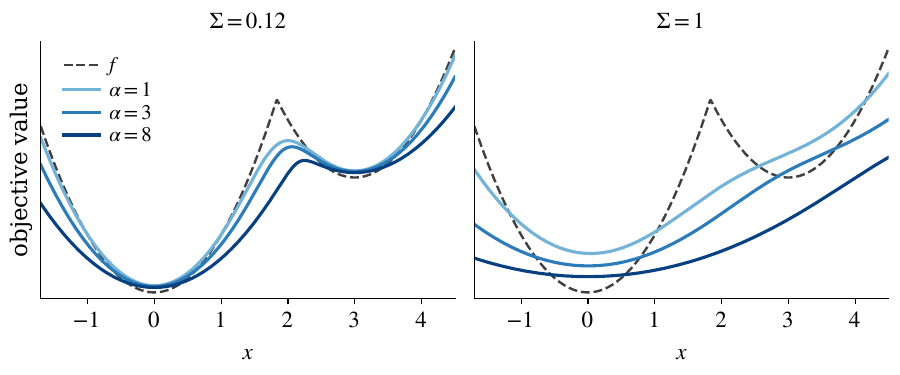}
\caption{The free energy $\mathcal{E}_\alpha$ for the objective \eqref{eq:twoWells} of Example~\ref{ex:twoWells} with $\Sigma=0.35^2\approx0.12$ (left) and $\Sigma=1$ (right), for $\alpha\in\{1,3,8\}$ (solid, darker for larger $\alpha$). 
The dashed curve is the objective $f$. 
}
\label{fig:twoWells}
\end{figure}

Elementary computations show that
\begin{equation}\label{eq:meanFieldGradientStep}
  \mathcal{T}_\alpha(x) = x - \alpha \Sigma\nabla \mathcal{E}_\alpha(x),
\end{equation}
suggesting that mean-field MPPI is a form of preconditioned gradient descent on a smoothed objective function, cf.\ \cite{fazlyab2026model} for an in-depth discussion of this connection.
The following example illustrates this point of view.
%
%
\begin{example}\label{ex:twoWells}
Let $d=1$, and consider the objective
\begin{equation}\label{eq:twoWells}
    f(y)=\min\Big\{\tfrac12y^2,\ 1+\tfrac12(y-3)^2\Big\},
\end{equation}
which has the global minimiser $x_\ast=0$ with $f_\ast=0$, a strict local minimiser at $y=3$ with $f(3)=1$, and a kink at $y=11/6$, where the two branches meet.
Note that $f$ fulfills Assumption \ref{ass:objective}, though it is neither convex nor differentiable. 
Figure~\ref{fig:twoWells} shows $\mathcal{E}_\alpha$ for $\Sigma\in\{0.35^2,1\}$ and $\alpha\in\{1,3,8\}$. 
For $\Sigma=0.35^2$, all three smoothed objectives keep a second well near the local minimiser $3$.
For $\Sigma=1$, the proposal reaches the deeper well from every centre in the displayed interval, and $\mathcal{E}_\alpha$ has a single minimum near $x_\ast$ already for $\alpha=1$.
\end{example}
As expected, Figure \ref{fig:twoWells} shows that $\mathcal{E}_\alpha$ smooths the objective function.
However, we observe that more than that happens: for large sampling covariance and inverse temperature, the smoothed objective looks increasingly convex.
Interpreting mean-field MPPI as preconditioned gradient descent then suggests that we should be able to get some form of global convergence to the global minimiser $x_\ast$, cf.\ Figure~\ref{fig:twoWellsIterates} in the supplementary material for an illustration.

We will now make this precise by building directly on the developments in \cite{fornasier2026consensus}.
The cornerstone of this approach is a \emph{quantitative Laplace principle}, which makes precise the intuition that for large $\alpha>0$ we should have
\begin{equation}
  \frac{\E_{X \distr \Norm(x,\Sigma)}[\exp(-\alpha f(X))X]}{\E_{X\distr\Norm(x,\Sigma)}[\exp(-\alpha f(X))]}
  \approx
  x_\ast 
\end{equation}
even if $x$ is far away from $x_\ast$, cf.\ Section \ref{sec:laplacePrinciple} in the supplementary material for details.
Such a result can be used to show that $\mathcal{T}_\alpha$ is a contraction on a region around $x_\ast$, which in turn leads to our first result: \emph{mean-field MPPI converges semi-globally}.
We would like to stress that this result is a direct generalisation of developments in \cite{fornasier2026consensus}, and we do not claim any novelty in terms of proof techniques.
%
%
\begin{theorem}\label{thm:meanFieldConvergence}
Let $f$ satisfy Assumption~\ref{ass:objective}, and let $\Sigma\in\R^{d\times d}$ be symmetric positive definite.
For all $r>0$ and $\rho\in(0,1)$ there exists $\alpha_0=\alpha_0(f,\Sigma,r,\rho)\ge0$ such that for all $\alpha>0$ with $\alpha\ge\alpha_0$ the following hold.
\begin{enumerate}
    \item[(i)] We have $\mathcal{T}_\alpha\big(\ball_r(x_\ast)\big)\subseteq\ball_r(x_\ast)$, and
    \begin{equation}\label{eq:meanFieldConvergence:contraction}
        \big\|\mathcal{T}_\alpha(x)-\mathcal{T}_\alpha(x')\big\|_2\;\le\;\rho\,\|x-x'\|_2
        \qquad\text{for all }x,x'\in\ball_r(x_\ast).
    \end{equation}
    \item[(ii)] $\mathcal{T}_\alpha$ has exactly one fixed point $m_\alpha$ in $\ball_r(x_\ast)$. For every $x_0\in\ball_r(x_\ast)$, the mean-field iterates $x_{k+1}=\mathcal{T}_\alpha(x_k)$ satisfy $x_k\in\ball_r(x_\ast)$ and
    \begin{equation}\label{eq:meanFieldConvergence:rate}
        \|x_k-m_\alpha\|_2\;\le\;\rho^{k}\,\|x_0-m_\alpha\|_2
        \qquad\text{for all }k\in\Nnn .
    \end{equation}
    Moreover, $\|m_\alpha-x_\ast\|_2\to0$ as $\alpha\to\infty$.
\end{enumerate}
\end{theorem}
This result makes the notion of semi-global convergence precise:
We cannot guarantee global convergence in the sense that for $\alpha\in\Rp$ large enough mean-field MPPI converges for all initial $x_0\in\R^d$.
Rather, for all $r\in\Rp$ we can find a lower bound $\alpha_0$ on the inverse temperature -- depending only on $r$, the contraction factor $\rho$, properties of $f$ as described in Assumption \ref{ass:objective}, and spectral properties of $\Sigma$ -- such that for all $\alpha \geq \alpha_0$ and all $x_0 \in \ball_r(x_\ast)$ we have geometric convergence, uniformly in $\ball_r(x_\ast)$, to a fixed point $m_\alpha$.
Furthermore, $m_\alpha$ does \emph{not} depend on $r$ in the sense that if $0<r\le r'$ and $\alpha$ exceeds both corresponding thresholds, then the fixed point in $\ball_r(x_\ast)$ coincides with the unique fixed point in $\ball_{r'}(x_\ast)$.

Theorem~\ref{thm:meanFieldConvergence} does not yet quantify how close the iterates get to the global minimiser $x_\ast$.
However, for the fixed point from above we have $m_\alpha \rightarrow x_\ast$ for $\alpha \rightarrow \infty$, and this convergence can be made quantitative.
This leads to our second main result: \emph{mean-field MPPI converges semi-globally to the global minimiser $x_\ast$}.
Similar results are contained in \cite{fornasier2026consensus} under stronger assumptions, cf.\ Section \ref{sec:suppl:mainresults} for a discussion.
%
%
\begin{theorem}\label{thm:meanFieldOneStep}
Let $f$ satisfy Assumption~\ref{ass:objective}, and let $\Sigma\in\R^{d\times d}$ be symmetric positive definite.
For every $\beta>0$ there exists $\alpha_{\mathrm b}=\alpha_{\mathrm b}(f,\Sigma,\beta)\ge0$ such that for all $r>0$, $\rho\in(0,1)$, and all $\alpha>0$ with $\alpha\ge\max\{\alpha_0(f,\Sigma,r,\rho),\alpha_{\mathrm b}(f,\Sigma,\beta)\}$, where $\alpha_0$ is from Theorem~\ref{thm:meanFieldConvergence}, the following hold.
\begin{enumerate}
    \item[(i)] For all $x\in\ball_r(x_\ast)$,
    \begin{equation}\label{eq:meanFieldOneStep:affine}
        \big\|\mathcal{T}_\alpha(x)-x_\ast\big\|_2\;\le\;\beta+\rho\,\|x-x_\ast\|_2 .
    \end{equation}
    \item[(ii)] For every $x_0\in\ball_r(x_\ast)$, the mean-field iterates satisfy
    \begin{equation}\label{eq:meanFieldOneStep:iterates}
        \|x_k-x_\ast\|_2\;\le\;\rho^k\,\|x_0-x_\ast\|_2+\frac{\beta}{1-\rho}
        \qquad\text{for all }k\in\Nnn .
    \end{equation}
\end{enumerate}
\end{theorem}
Observe that the threshold $\alpha_0(f,\Sigma,r,\rho)$ depends on the radius $r$, but not on the accuracy parameter $\beta$, 
while $\alpha_{\mathrm b}(f,\Sigma,\beta)$ depends on $\beta$, but not on $r$.
Interestingly, \eqref{eq:meanFieldOneStep:affine} shows that for large $\alpha$ already \emph{one mean-field MPPI} step can lead to good accuracy, which provides a partial theoretical justification, on the mean-field level, for the single-step variant originally proposed in \cite{williams2018information}.
However, \eqref{eq:meanFieldOneStep:iterates} shows that iteration on the mean-field can indeed lead to improvements of the accuracy.
This phenomenon has been pointed out already in \cite[Remark~4.5]{fornasier2026consensus}, albeit under stronger assumptions.
The preceding results imply in particular that for an arbitrary initial guess $x_0$ an arbitrary accuracy can be achieved after enough iterations for a large enough $\alpha>0$, as made precise in the next result.
%
%
\begin{corollary}\label{cor:meanFieldAccuracy}
Let $f$ satisfy Assumption~\ref{ass:objective}, and let $\Sigma\in\R^{d\times d}$ be symmetric positive definite.
Let $r>0$ and $\epsilon>0$, and set $K_0(r,\epsilon):=\max\{0,\lceil\log_2(2r/\epsilon)\rceil\}$.
Then for all $\alpha>0$ with $\alpha\ge\max\{\alpha_0(f,\Sigma,r,\frac12),\alpha_{\mathrm b}(f,\Sigma,\frac\epsilon4)\}$ and all $x_0\in\ball_r(x_\ast)$, the mean-field iterates satisfy
\begin{equation}\label{eq:meanFieldAccuracy:claim}
    \|x_k-x_\ast\|_2\;\le\;\epsilon \qquad \text{for all } k \geq K_0(r,\epsilon).
\end{equation}
\end{corollary}
Note that the number of iterations $K_0(r,\epsilon)$ depends neither on $f$ nor on $\Sigma$ nor on $\alpha$.
We would like to stress that for these results it is enough that $f$ fulfils Assumption \ref{ass:objective}.
In particular, $f$ might be severely non-convex and non-smooth.

Finally, under Assumption~\ref{ass:simplified}, the thresholds become explicit. 
We state them up to logarithmic factors and refer to Corollary~\ref{cor:simplifiedAccuracy} and Lemma~\ref{lem:simplifiedThresholdOrders} in the supplementary material for detailed statements.
%
%
\begin{corollary}\label{cor:meanFieldSimplified}
Let $f$ satisfy (SF1) and (SF2) of Assumption~\ref{ass:simplified}, let $\Sigma\in\R^{d\times d}$ be symmetric positive definite, and let $r>0$ and $\epsilon\in(0,2r]$.
There is a threshold $\alpha^{\mathrm S}_{\mathrm{mf}}(r,\epsilon)$ with
\begin{equation}\label{eq:meanFieldSimplified:order}
    \alpha^{\mathrm S}_{\mathrm{mf}}(r,\epsilon)=\tilde{\mathcal O}\Big(\frac{d}{m\,\epsilon^2}+\frac{(r+R)^2}{\Delta\,\sigma_{\min}(\Sigma)^2}+\frac{1}{m\,\sigma_{\min}(\Sigma)^2}+\frac d\Delta\Big)
\end{equation}
such that for all $\alpha\ge\alpha^{\mathrm S}_{\mathrm{mf}}(r,\epsilon)$ and all $x_0\in\ball_r(x_\ast)$, the mean-field iterates satisfy $\|x_k-x_\ast\|_2\le\epsilon$ for all $k\ge K_0(r,\epsilon)$.
\end{corollary}
Note that the accuracy and the radius enter the inverse temperature additively,
where the dependence on the accuracy is of order $d/(m\epsilon^2)$, and the radius only enters through the term $(r+R)^2/(\Delta\sigma_{\min}(\Sigma)^2)$, which does not depend on $\epsilon$. 

\section{Finite-Sample Results}\label{sec:finiteSample}
While the mean-field results from the preceding section are strong, it is a priori unclear how applicable they are to the actual finite-sample method, i.e., the MPPI iteration \eqref{eq:mppiIteration}.
To investigate this issue, we now turn to finite-sample results.
We focus on finite horizons, i.e., on a fixed number of iterations.
To simplify the notation, for $x\in\R^d$ define the
\emph{empirical transition map}
\begin{equation}\label{eq:empiricalTransitionMap}
    \widehat{\mathcal{T}}_{\alpha,N}(x):=\frac{\sum_{i=1}^N e^{-\alpha f(X_i)}X_i}{\sum_{i=1}^N e^{-\alpha f(X_i)}}
\end{equation}
with $X_1,\dots,X_N$ independent with $X_i\distr\Norm(x,\Sigma)$ as before.

\citet{fornasier2026consensus} suggest deriving finite-sample results via \emph{propagation of chaos} results which are already available for CH \citep{haug2026quantitative}.
However, as pointed out above, it is well-known that $ \widehat{\mathcal{T}}_{\alpha,N}(x)$ is a self-normalised importance-weighted sampling estimate of $\mathcal{T}_\alpha(x)$.
From an abstract standpoint, this means that the MPPI or empirical CH iteration is not a full interacting particle system, but rather a time-homogeneous Markov chain with continuous state space.
We can use this connection to directly derive finite-sample results without first establishing propagation of chaos, similar to earlier works like \cite{yoon2022sampling,yoon2026finite}.
We start with a non-asymptotic error bound for the self-normalised importance sampling estimate $\widehat{\mathcal{T}}_{\alpha,N}(x)$ of $\mathcal{T}_\alpha(x)$.
While we are not aware of this exact result in the literature, such bounds are well-known in the context of self-normalised estimators, see for example \cite[Theorem~2.1]{agapiou2017importance}.
%
%
\begin{proposition}\label{prop:snisErrorBound}
Let $f$ satisfy Assumption~\ref{ass:objective}, let $\Sigma\in\R^{d\times d}$ be symmetric positive definite, and let $\alpha>0$.
Then for all $x\in\R^d$, $N\in\Np$ and $t>0$,
\begin{equation}\label{eq:snisErrorBound:bound}
    \Pb\Big(\big\|\widehat{\mathcal{T}}_{\alpha,N}(x)-\mathcal{T}_\alpha(x)\big\|_2>t\Big)
    \;\le\;\frac{4}{N}\,\frac{1}{\zeta_\alpha(x)}\left(1+\frac{V_\alpha(x)}{t^2}\right),
\end{equation}
where
\begin{equation}\label{eq:snis:quantities}
\begin{gathered}
    \zeta_\alpha(x):=\E_{X\distr\Norm(x,\Sigma)}\big[e^{-\alpha(f(X)-f(x_\ast))}\big]\in(0,1],
    \\
    V_\alpha(x):=\frac{\E_{X\distr\Norm(x,\Sigma)}\big[e^{-\alpha f(X)}\,\|X-\mathcal{T}_\alpha(x)\|_2^2\big]}{Z_\alpha(x)}<\infty .
\end{gathered}
\end{equation}
\end{proposition}
The factor $1/\zeta_\alpha(x)$ plays the role of an inverse effective sample size. It is nondecreasing in the inverse temperature $\alpha$ and typically grows with the distance of $x$ to $x_\ast$, since then fewer samples carry a substantial weight.
We now combine Proposition~\ref{prop:snisErrorBound} with the semi-global mean-field convergence result Theorem~\ref{thm:meanFieldOneStep}.
The key observation for establishing the following result is that the sampling error only needs to be small compared with the current distance to $x_\ast$.
Hence steps far away from $x_\ast$ only require a coarse accuracy, and the fine accuracy only has to be achieved close to $x_\ast$.
%
%
\begin{theorem}\label{thm:finiteSample}
Let $f$ satisfy Assumption~\ref{ass:objective}, let $\Sigma\in\R^{d\times d}$ be symmetric positive definite, and let $(\hat x_k)_{k\in\Nnn}$ be the MPPI iterates \eqref{eq:mppiIteration}.
Let $r>0$, $\rho\in(0,1)$ and $\beta>0$ with $4\beta\le(1-\rho)r$, let $\gamma>0$, and set $\rho':=(1+\rho)/2$.
There exist $\alpha_0'=\alpha_0'(f,\Sigma,r,\rho)$, $\alpha_{\mathrm b}'=\alpha_{\mathrm b}'(f,\Sigma,\beta)$ and $C=C(f,\Sigma,r,\gamma)$ such that for all $\alpha>0$ with $\alpha\ge\max\{\alpha_0',\alpha_{\mathrm b}'\}$, all $K\in\Np$, $\delta\in(0,1)$, $N\in\Np$ with
\begin{equation}\label{eq:finiteSample:N}
    N\;\ge\;\frac{K}{\delta}\,C\,e^{\gamma\alpha},
\end{equation}
and all $\hat x_0\in\ball_r(x_\ast)$, we have
\begin{equation}\label{eq:finiteSample:claim}
    \Pb\Big(\hat x_k\in\ball_r(x_\ast)\ \text{and}\ \|\hat x_k-x_\ast\|_2\le\rho'^k\,\|\hat x_0-x_\ast\|_2+\frac{4\beta}{1-\rho}\ \text{ for all }k=0,\dots,K\Big)\;\ge\;1-\delta .
\end{equation}
In particular, let $\epsilon\in(0,2r]$ and $K_\epsilon:=\max\{0,\lceil\log(2r/\epsilon)/\log(4/3)\rceil\}$. Then for all $\alpha>0$ with $\alpha\ge\max\{\alpha_0'(f,\Sigma,r,\frac12),\alpha_{\mathrm b}'(f,\Sigma,\frac\epsilon{16})\}$, all $K\in\Np$, $\delta\in(0,1)$, $N$ as in \eqref{eq:finiteSample:N}, and $\hat x_0\in\ball_r(x_\ast)$,
\begin{equation}\label{eq:finiteSample:accuracy}
    \Pb\Big(\|\hat x_k-x_\ast\|_2\le\epsilon\ \text{ for all }k\text{ with }K_\epsilon\le k\le K\Big)\;\ge\;1-\delta .
\end{equation}
\end{theorem}
As in the mean-field case, the inverse temperature thresholds are decoupled: $\alpha_0'$ depends on the radius $r$, but not on the accuracy, and $\alpha_{\mathrm b}'$ depends on the accuracy, but not on $r$.
The sample size \eqref{eq:finiteSample:N} grows linearly in $K/\delta$ and like $e^{\gamma\alpha}$ in $\alpha$, where $\gamma>0$ can be chosen arbitrarily small at the price of a larger constant $C$.
Note that since the constant $C$ depends on $r$, the sample size still couples the radius with the inverse temperature, and hence with the final accuracy.
In particular, also this result is semi-global like Theorem \ref{thm:meanFieldOneStep}.

Under Assumption~\ref{ass:simplified}, the sample size becomes explicit as well. The precise statements are Corollary~\ref{cor:simplifiedSampleComplexity} and Lemma~\ref{lem:simplifiedThresholdOrders} in the supplementary material.
%
%
\begin{corollary}\label{cor:finiteSampleSimplified}
Let $f$ satisfy Assumption~\ref{ass:simplified}, let $\Sigma\in\R^{d\times d}$ be symmetric positive definite, and let $r>0$ and $\epsilon\in(0,2r]$.
There is a threshold $\alpha^{\mathrm S}_{\mathrm{fs}}(r,\epsilon)$ of the same order \eqref{eq:meanFieldSimplified:order} as in Corollary~\ref{cor:meanFieldSimplified} such that for all $\alpha\ge\alpha^{\mathrm S}_{\mathrm{fs}}(r,\epsilon)$, $K\in\Np$, $\delta\in(0,1)$, $N\in\Np$ with
\begin{equation}\label{eq:finiteSampleSimplified:N}
    N\;\ge\;\frac{12K}{\delta}\,\big(1+\alpha M\sigma_{\max}(\Sigma)^2\big)^{d/2}\exp\Big(\frac{r^2}{2\sigma_{\min}(\Sigma)^2}\Big),
\end{equation}
and all $\hat x_0\in\ball_r(x_\ast)$, the MPPI iterates satisfy
\begin{equation*}
    \Pb\Big(\|\hat x_k-x_\ast\|_2\le\epsilon\ \text{ for all }k\text{ with }K_0(r,\epsilon)\le k\le K\Big)\;\ge\;1-\delta .
\end{equation*}
\end{corollary}
For small $\epsilon$, the threshold is dominated by the term $384d/(m\epsilon^2)$, cf.\ Lemma~\ref{lem:simplifiedThresholdOrders} in the supplementary material, and for $\alpha$ of this order the sample size \eqref{eq:finiteSampleSimplified:N} is of order
\begin{equation*}
    \frac K\delta\,e^{r^2/(2\sigma_{\min}(\Sigma)^2)}\Big(\frac{384\,dM\sigma_{\max}(\Sigma)^2}{m\epsilon^2}\Big)^{d/2},
\end{equation*}
i.e., it grows like $\epsilon^{-d}$ in the accuracy.

\section{Applications to Optimal Control} \label{sec:oc}
As a concrete illustration of the developments above, let us consider the classic finite-horizon LQR problem.
We consider a deterministic LTI system
\begin{equation}\label{eq:control:lti}
    s_+ = As + Ba,
\end{equation}
where $A\in\R^{n\times n}$ is the state transition matrix,
$B\in\R^{n\times p}$ the control input matrix,
$s\in\R^n$ the current state, $s_+$ the next state, and $a\in\R^p$ the control input.
For an initial state $s_0\in\R^n$, a horizon $T\in\Np$ and a control sequence $a=(a_0,\dots,a_{T-1})\in(\R^p)^T$,
let $s(t;s_0,a)$ denote the state at time $t=0,\ldots,T$, i.e., $s(0;s_0,a)=s_0$ and $s(t+1;s_0,a)=As(t;s_0,a)+Ba_t$ for $t=0,\dots,T-1$.
If no confusion can arise, we also write $s_t=s(t;s_0,a)$.
Note that we use RL-type notation to avoid a notation clash with our preceding analysis.

For horizon $T\in\Np$, consider the total cost functional
\begin{equation}\label{eq:control:cost}
    J_T(s_0, a) = \sum_{t=0}^{T-1} s_t^\top Q s_t + a_t^\top R a_t + 2s_t^\top S a_t
\end{equation}
for matrices $Q\in\R^{n\times n}$, $R\in\R^{p\times p}$ and $S\in\R^{n\times p}$.
We are interested in the optimisation problem
\begin{equation}\label{eq:control:lqocp}
    \min_{a\in(\R^p)^T} J_T(s_0,a).
\end{equation}
Without loss of generality, we assume that $Q$ and $R$ are symmetric, since only the quadratic forms $s^\top Qs$ and $a^\top Ra$ enter the cost.
As is well-known, the objective function $J_T$ can be rewritten as
\begin{equation}\label{eq:control:quadratic}
    J_T(s_0,a)=a^\top H_Ta+2a^\top F_Ts_0+s_0^\top Y_Ts_0 ,
\end{equation}
where the control sequence $a$ is identified with the stacked vector $(a_0^\top,\dots,a_{T-1}^\top)^\top\in\R^{pT}$, and $H_T\in\R^{pT\times pT}$ (symmetric), $F_T\in\R^{pT\times n}$ and $Y_T\in\R^{n\times n}$ depend only on $A$, $B$, $Q$, $R$, $S$ and $T$, cf.\ Section~\ref{app:proofsOC} in the supplementary material.
For simplicity, assume that $\begin{pmatrix} Q & S \\ S^\top & R\end{pmatrix}$ is positive definite (which implies that $R$ is positive definite), so that $H_T$ is also positive definite, cf.\ Lemma~\ref{lem:control:HTbounds} in the supplementary material, and hence \eqref{eq:control:lqocp} has the unique solution $a_\ast(s_0)=-H_T^{-1}F_Ts_0$.
Finally, let $\Sigma_c\in\R^{p\times p}$ be symmetric positive definite, so that also $\Sigma=\mathrm{diag}(\Sigma_c,\ldots,\Sigma_c)$
is symmetric positive definite.
For fixed $s_0\in\R^n$ we can then set $d=pT$ and apply the MPPI iteration to the objective function $f=J_T(s_0,\cdot)$ with $\Sigma$ from above and inverse temperature $\alpha$.
Let us write $a_\ast:=a_\ast(s_0)$, $\bar a_{k+1}:=\mathcal{T}_\alpha(\bar a_k)$ for the mean-field iterates, and $(\hat a_k)_{k\in\Nnn}$ for the MPPI iterates \eqref{eq:mppiIteration} with this objective $f$.
%
%
\begin{corollary}[Explicit results for the LQR problem]\label{cor:control:explicit}
Let $r>0$ and $\epsilon\in(0,2r]$.
\begin{enumerate}
    \item[(i)] For all $\alpha>\frac{16pT}{\lambda_{\min}(H_T)\epsilon^2}+\frac{6}{\lambda_{\min}(H_T)\sigma_{\min}(\Sigma_c)^2}$ and $\bar a_0\in\ball_r(a_\ast)$, the mean-field iterates satisfy $\|\bar a_k-a_\ast\|_2\le\epsilon$ for all $k\ge K_0(r,\epsilon)$.
    \item[(ii)] For all $\alpha>\frac{192pT}{\lambda_{\min}(H_T)\epsilon^2}+\frac{10}{\lambda_{\min}(H_T)\sigma_{\min}(\Sigma_c)^2}$, $K\in\Np$, $\delta\in(0,1)$, $N\in\Np$ with
    \begin{equation}\label{eq:control:explicit:N}
        N\;\ge\;\frac{12K}{\delta}\big(1+2\alpha\lambda_{\max}(H_T)\sigma_{\max}(\Sigma_c)^2\big)^{pT/2}\exp\Big(\frac{r^2}{2\sigma_{\min}(\Sigma_c)^2}\Big),
    \end{equation}
    and $\hat a_0\in\ball_r(a_\ast)$, the MPPI iterates satisfy
    \begin{equation*}
        \Pb\Big(\|\hat a_k-a_\ast\|_2\le\epsilon\ \text{ for all }k\text{ with }K_0(r,\epsilon)\le k\le K\Big)\;\ge\;1-\delta .
    \end{equation*}
\end{enumerate}
\end{corollary}
For a quadratic objective, the radius $r$ does not enter the inverse temperature, since the well of Assumption~\ref{ass:simplified} can be taken arbitrarily large.
For small $\epsilon$ and $\alpha$ at the threshold in (ii), the sample size \eqref{eq:control:explicit:N} is of order $\frac K\delta\,e^{r^2/(2\sigma_{\min}(\Sigma_c)^2)}\big(384\,pT\kappa\,\sigma_{\max}(\Sigma_c)^2/\epsilon^2\big)^{pT/2}$ with $\kappa:=\lambda_{\max}(H_T)/\lambda_{\min}(H_T)$, i.e., the exponent of $1/\epsilon$ grows linearly in the horizon $T$.
Since $\lambda_{\min}(H_T)\ge\lambda_{\min}\big(\begin{smallmatrix}Q&S\\S^\top&R\end{smallmatrix}\big)$ by Lemma~\ref{lem:control:HTbounds}, the thresholds on $\alpha$ remain valid with this lower bound, and then depend on the dynamics only through the dimension $pT$. By the same lemma, $\kappa$ is bounded uniformly in $T$ if $A$ is Schur stable.
Finally, we would like to stress that for a quadratic objective function as arising in the LQR problem,
one can conduct a more tailored analysis, as has been done in related work like \cite{yi2024covo} or \cite{yoon2026finite}.
Since we use the LQR problem as a mere illustration, we do not pursue this here.

\section{Conclusion}\label{sec:conclusion}
We investigated MPPI as a sampling-based global optimisation method for non-convex and non-smooth objectives.
For objective functions with a unique global minimum, under weak additional assumptions, we established semi-global convergence to the minimiser of the mean-field MPPI iteration, as well as semi-global finite-sample results for the MPPI iteration, both in the low-temperature setting.
To the best of our knowledge, these are the first such results for general non-smooth and non-convex objective functions.
We would like to stress that we considered only the low-temperature (large inverse temperature) setting.
We expect that for general objective functions as considered here, no reliable semi-global convergence can be shown at low inverse temperatures,
but for restricted objective classes such results should be possible.
Furthermore, we point out again that our finite-sample results are based on effective, but simple error bounds for the self-normalised importance weighted estimator, and we expect that with more sophisticated techniques and for more restricted function classes better results can be achieved.
Regarding the practical applicability of our results, we stress that our results are quantitative and require only reasonable assumptions, but we do not claim that the various constants necessary to apply the results can be easily found.
Rather, our work contributes to a better conceptual understanding of MPPI, and similarly to no-regret results in the bandit literature, our results show that there are algorithmic parameters (in particular, inverse temperature thresholds and sample bounds) that are benign and that ensure the algorithmic goal (semi-global convergence) is achieved.
Let us close by pointing out two important directions for future work.
First, we focused on the popular plain MPPI iteration \eqref{eq:mppiIteration} in the low-temperature setting.
However, there is by now a multitude of MPPI variants \citep{kazim2024recent} and also other scaling limits as in \cite{homburger2025optimality}.
Investigating more variants with our techniques is therefore an interesting and relevant problem.
Second, we focused on the optimisation behaviour of MPPI. 
Since its main use is in control applications, it is important to investigate control-theoretic properties that make our results applicable to the optimal control objective function, similar to our illustration in Section \ref{sec:oc} for LQR. 
Furthermore, the suboptimality needs to be accounted for in the control scheme.
We expect that our results will be an important step towards such results.


\clearpage

\subsection*{Generative AI Use}
The present project and the high-level approach have been conceived by the authors without AI use.
The mean-field analysis is based on initial work by the authors without AI use,
then refined using Claude Opus 5.5 (in particular, the reduction to Assumption \ref{ass:objective} was found by this AI model), with final editing and rewriting by the authors.
The finite-sample results as well as the results based on the simplified Assumption \ref{ass:simplified} have been developed interactively with ChatGPT-6, and the write-up of technical intermediate results as well as their proofs was supported by Claude Opus 5.5.
The LQR example was developed by the authors with the writeup supported by Claude Opus 5.5.
All illustrations have been conceptualised by the authors without AI, and then typeset and programmed with Claude Opus 5.5.
The main text has been written by the authors; several parts (in particular Assumption~\ref{ass:simplified} and its discussion, Corollaries~\ref{cor:meanFieldSimplified}, \ref{cor:finiteSampleSimplified} and~\ref{cor:control:explicit}, parts of Section~\ref{sec:oc}, and the figure captions) were drafted with Claude Opus 5.5 and subsequently revised by the authors.

\acks{
We would like to thank Jona Klemenc for very helpful discussions, especially related to \cite{fornasier2026consensus}.
C.F.\ acknowledges funding from DFG Project FO 767/10-2 (eBer-24-32734) ``Implicit Bias in Adversarial Training''.
G.B.\ acknowledges support from the Munich Center for Machine Learning.
\begin{center}
    \raisebox{0pt}{\includegraphics[height=1.5cm]{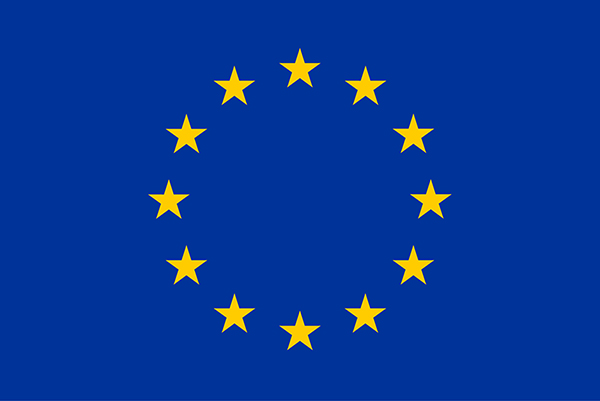}}%
    \hspace{1em}%
    \raisebox{0pt}{\includegraphics[height=1.5cm]{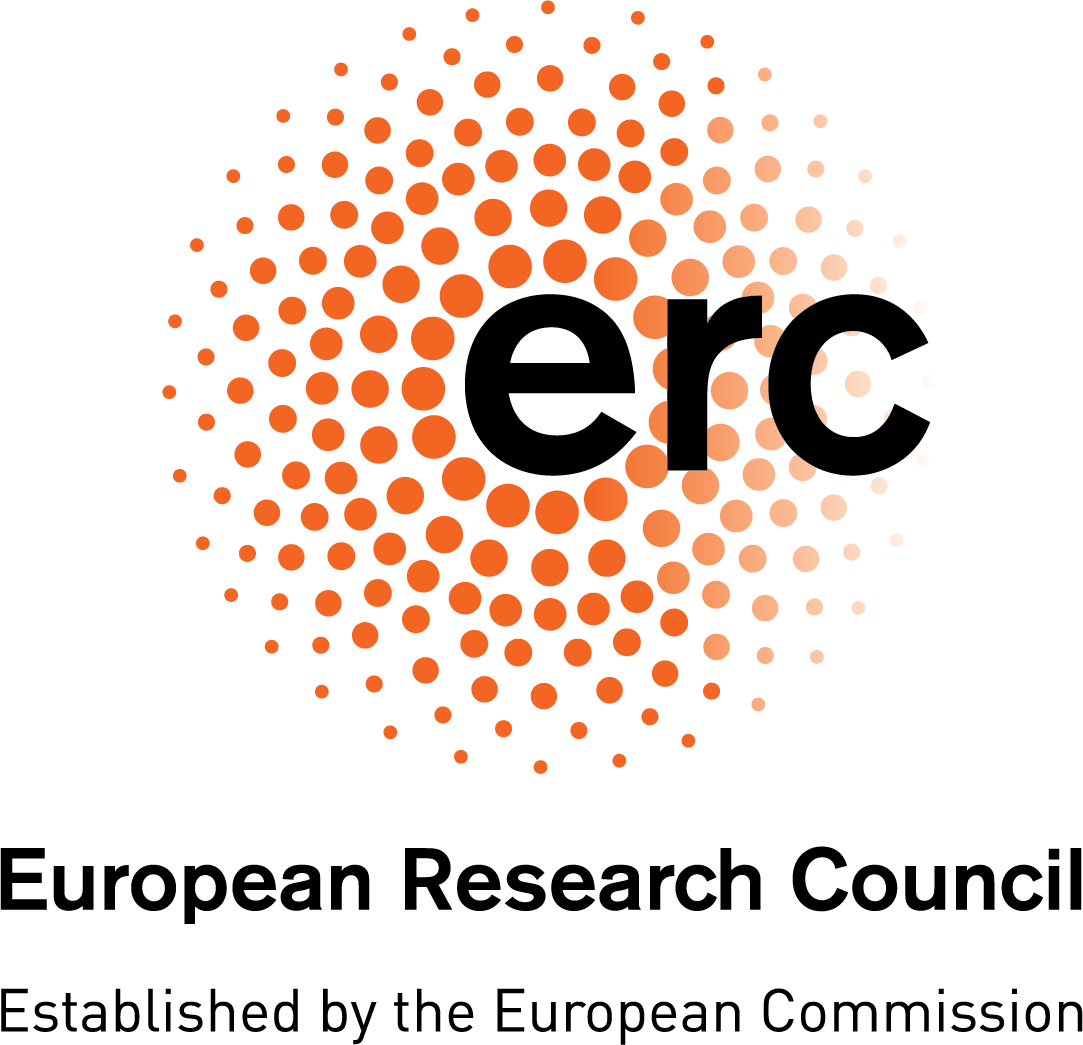}}%
\end{center}
Funded by the European Union. Views and opinions expressed are however those of the author(s) only and do not necessarily reflect those of the European Union or the European Research Council Executive Agency. Neither the European Union nor the granting authority can be held responsible for them. This project has received funding from the European Research Council (ERC) under the European Union’s Horizon Europe research and innovation programme (grant agreement No. 101198055, project acronym NEITALG).
}

\bibliography{refs}

\clearpage
\appendix

\section{Supplementary Material}\label{app:supplementary}

\paragraph{Notation}
We denote the Euclidean norm of vectors by $\|\cdot\|_2$ and the spectral norm of matrices by $\|\cdot\|$.
For $x\in\R^d$ and $r\in\Rp$, we define $\ball_r(x)=\{x' \in\R^d \mid \|x-x'\|_2\le r \}$.
We denote the set of Borel probability measures on $\R^d$ by $\Pbs(\R^d)$,
and with measurable we always mean Borel measurable.
We denote the Gaussian measure on $\R^d$ with mean $x$ and (symmetric positive definite) covariance matrix $\Sigma$ by $\Norm(x,\Sigma)$,
and the corresponding Gaussian density evaluated at $y\in\R^d$ by $\Norm(y \mid x,\Sigma)$.
For a symmetric matrix $A\in\R^{n\times n}$ we write $\lambda_{\min}(A)$ and $\lambda_{\max}(A)$ for the smallest and largest eigenvalue, respectively (note that eigenvalues of a symmetric matrix are necessarily real).
For a symmetric positive definite matrix $\Sigma$, we write $\sigma_{\min}(\Sigma):=\lambda_{\min}(\Sigma)^{1/2}$ and $\sigma_{\max}(\Sigma):=\lambda_{\max}(\Sigma)^{1/2}$ for the square roots of its smallest and largest eigenvalue, and $\kappa_\Sigma:=\sigma_{\max}(\Sigma)^2/\sigma_{\min}(\Sigma)^2$. In particular, $\|\Sigma\|=\sigma_{\max}(\Sigma)^2$ and $\|\Sigma^{-1}\|=\sigma_{\min}(\Sigma)^{-2}$.
For $\nu\in\Pbs(\R^d)$ and $\alpha>0$, we write $Z_\alpha(\nu):=\int_{\R^d}e^{-\alpha f(y)}\,\nu(\mathrm dy)$, so that $Z_\alpha(x)=Z_\alpha(\Norm(x,\Sigma))$ in the notation of the main text; in Lemma~\ref{lem:laplacePrinciples}, where $\nu$ is fixed, we simply write $Z_\alpha$.
Finally, as usual the Landau symbol $\tilde{\mathcal{O}}(\cdot)$ suppresses poly-logarithmic terms.

\subsection{Details for Example \ref{ex:twoWells}}
Let us verify that the function $ f(y)=\min\Big\{\tfrac12y^2,\ 1+\tfrac12(y-3)^2\Big\}$ fulfils Assumption \ref{ass:objective}.
Clearly, $f$ is continuous and has a unique global minimum, so (F1) and (F2) hold. 
Moreover, if $|y|\ge\epsilon$, then either $y\le11/6$ and $f(y)=\frac12y^2\ge\frac12\epsilon^2$, or $y>11/6$ and $f(y)\ge1$. 
Hence (F3) holds with $S_\epsilon=\min\{\epsilon^2/2,1\}$, and $f$ satisfies Assumption~\ref{ass:objective}, although it is neither convex nor differentiable.

\begin{figure}[!htbp]
\centering
\includegraphics[width=\textwidth]{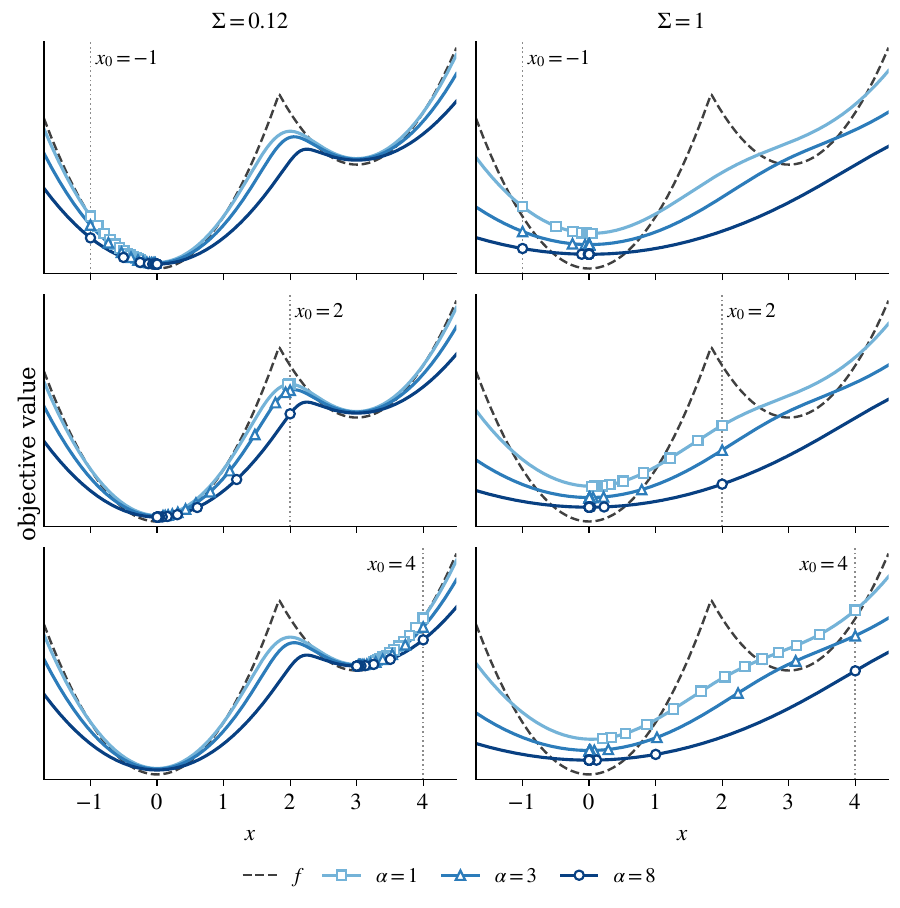}
\caption{Mean-field MPPI iterates for the objective \eqref{eq:twoWells} of Example~\ref{ex:twoWells} with $\Sigma=0.35^2\approx0.12$ (left) and $\Sigma=1$ (right), started at $x_0=-1$, $2$ and $4$ (top to bottom; the dotted vertical line marks $x_0$).
The markers show $(x_k,\mathcal{E}_\alpha(x_k))$, $k=0,\ldots,12$, for $\alpha=1$ (squares), $\alpha=3$ (triangles) and $\alpha=8$ (circles), and lie on the solid curve $\mathcal{E}_\alpha$ of the same colour (darker for larger $\alpha$). The dashed curve is the objective $f$.
}
\label{fig:twoWellsIterates}
\end{figure}

Figure~\ref{fig:twoWellsIterates} shows the mean-field MPPI iterates $x_{k+1}=\mathcal{T}_\alpha(x_k)$, $k=0,\ldots,11$, started at $x_0\in\{-1,2,4\}$.
For $\Sigma=0.35^2$, the iterates started at $x_0=4$ converge to the fixed point near the local minimiser $3$ for all three values of $\alpha$, in line with the second well of $\mathcal{E}_\alpha$ in Figure~\ref{fig:twoWells}.
Started at $x_0=2$, the iterates reach a neighbourhood of $x_\ast$ for $\alpha\in\{3,8\}$, whereas for $\alpha=1$ they have barely moved after twelve steps ($x_{12}\approx1.97$).
For $\alpha=1$, the small step size $\alpha\Sigma\approx0.12$ in \eqref{eq:meanFieldGradientStep} also slows down the iterates started at $x_0=-1$.
For $\Sigma=1$, all iterates approach $x_\ast$. For $\alpha=1$, the iterates started at $-1$ and $2$ settle at about $0.04$, while the one started at $4$ is still at about $0.2$ after twelve steps; for $\alpha=3$, all iterates settle at about $0.004$. This remaining offset reflects that the fixed point $m_\alpha$ differs from $x_\ast$ for finite $\alpha$.

\subsection{General Mean-Field Analysis}
\subsubsection{Laplace Principle} \label{sec:laplacePrinciple}
Our analysis is based on the following version of the quantitative Laplace principle, which is a variation of \cite[Lemma~3.3]{fornasier2026consensus} and \cite[Theorem~4.3]{huang2026faithful}.
Its proof is based directly on arguments from these references, but since we use slightly different assumptions, we provide a complete and self-contained proof below.
%
%
\begin{lemma}[Quantitative Laplace principle under separation]\label{lem:laplacePrinciples}
Let $f:\R^d\to\R$ be a Borel measurable function with global minimum $f_\ast\in\R$, attained at $x_\ast\in\R^d$, i.e.\ $f_\ast=f(x_\ast)\le f(x)$ for all $x\in\R^d$.
Furthermore, assume that $f$ is upper semicontinuous at $x_\ast$,
and that for every $\epsilon>0$ there exists $S_\epsilon>0$ such that
\begin{equation}\label{eq:laplacePrinciples:separation}
  \inf_{\|x-x_\ast\|_2\ge\epsilon} f(x)\;\ge\; f_\ast+S_\epsilon,
\end{equation}
which implies that $x_\ast$ is the unique global minimiser of $f$. 

Let $\nu$ be a Borel probability measure on $\R^d$ such that
\begin{equation}\label{eq:laplacePrinciples:mass}
  \nu\big(\ball_s(x_\ast)\big)\;\ge\; C_s>0\qquad\text{for all } s>0.
\end{equation}
Then for all $\alpha>0$, 
\begin{equation}\label{eq:laplacePrinciples:eta}
  \eta^\alpha_\nu(A)\;:=\;\frac{1}{Z_\alpha}\int_A e^{-\alpha f(x)}\,\nu(dx),
  \qquad
  Z_\alpha:=\int_{\R^d}e^{-\alpha f(x)}\,\nu(dx),
  \qquad A\in\mathcal B(\R^d),
\end{equation}
is a well-defined Borel probability measure,
and for every $\epsilon>0$, there exists $s_f(\epsilon)>0$ with
\begin{equation}\label{eq:laplacePrinciples:epsf}
  \sup_{\|x-x_\ast\|_2\le s_f(\epsilon)} f(x)\;\le\; f_\ast+\tfrac12 S_\epsilon ,
\end{equation}
such that the following hold.
\begin{enumerate}
    \item[(i)] For all $\epsilon>0$ and $\alpha>0$,
        \begin{equation}\label{eq:laplacePrinciples:tailmass}
            \eta^\alpha_\nu\big(\{x\in\R^d:\ \|x-x_\ast\|_2\ge\epsilon\}\big) \;\le\;\frac{e^{-\alpha S_\epsilon/2}}{C_{s_f(\epsilon)}} .
        \end{equation}
        In particular, for every $\delta>0$ and
        \begin{equation}\label{eq:laplacePrinciples:alpha0}
            \alpha_0=\alpha_0\big(\delta,S_\epsilon,C_{s_f(\epsilon)}\big)
                :=\frac{2}{S_\epsilon}\,\max\left\{0,\;\log\frac{1}{\delta\,C_{s_f(\epsilon)}}\right\},
        \end{equation}
        we have
        \begin{equation}\label{eq:laplacePrinciples:alpha0bound}
            \frac{e^{-\alpha S_\epsilon/2}}{C_{s_f(\epsilon)}}\;\le\;\delta
            \qquad\text{for all }\alpha\ge\alpha_0,
        \end{equation}
        so that the left-hand side of \eqref{eq:laplacePrinciples:tailmass} is at most $\delta$ for all $\alpha\ge\alpha_0$. Moreover,
        \begin{equation}
            \eta^\alpha_\nu\big(\{x\in\R^d:\ \|x-x_\ast\|_2\ge\epsilon\}\big) \rightarrow 0 \qquad \text{for} \qquad \alpha \rightarrow \infty.
        \end{equation}
    \item[(ii)] \emph{Moments.} Let $p\in[1,\infty)$ and assume that
        \begin{equation}\label{eq:laplacePrinciples:moment}
            \E_{X\distr\nu}\big[\|X-x_\ast\|_2^p\big]\;\le\; M_p<\infty .
        \end{equation}
        Then $\E_{X\distr\eta^\alpha_\nu}\big[\|X\|_2^p\big]<\infty$ and, for all $\epsilon>0$ and $\alpha>0$,
        \begin{equation}\label{eq:laplacePrinciples:momentbound}
            \E_{X\distr\eta^\alpha_\nu}\big[\|X-x_\ast\|_2^p\big]
                    \;\le\;\epsilon^p+\frac{e^{-\alpha S_\epsilon/2}}{C_{s_f(\epsilon)}}\,M_p .
        \end{equation}
        Consequently, $\E_{X\distr\eta^\alpha_\nu}\big[\|X-x_\ast\|_2^p\big]\to0$ as $\alpha\to\infty$.
        In particular, if \eqref{eq:laplacePrinciples:moment} holds for $p=1$, then the mean $\E_{X\distr\eta^\alpha_\nu}[X]$ exists, and for all $\epsilon>0$ and $\alpha>0$ we have
        \begin{equation}\label{eq:laplacePrinciples:mean}
            \big\|\E_{X\distr\eta^\alpha_\nu}[X]-x_\ast\big\|_2
                    \;\le\;\epsilon+\frac{e^{-\alpha S_\epsilon/2}}{C_{s_f(\epsilon)}}\,M_1 ,
        \end{equation}
        and $\E_{X\distr\eta^\alpha_\nu}[X]\to x_\ast$ as $\alpha\to\infty$.
\end{enumerate}
\end{lemma}

\begin{proof}
\textbf{Uniqueness of global minimiser.} Consider some $x\neq x_\ast$.
Using \eqref{eq:laplacePrinciples:separation} with $\epsilon=\|x-x_\ast\|_2>0$ results in $f(x)\ge f_\ast+S_\epsilon>f_\ast$,
so $x_\ast$ is indeed the unique global minimiser.

\textbf{Existence of $s_f(\epsilon)$.} Let $\epsilon>0$.
By definition of upper semicontinuity at $x_\ast$ there is $r>0$ with $f(x)<f_\ast+\tfrac12S_\epsilon$ for all $x$ with $\|x-x_\ast\|_2<r$, 
and any $s_f(\epsilon)\in(0,r)$ satisfies \eqref{eq:laplacePrinciples:epsf}.

For the remainder of the proof, let $\epsilon>0$ be arbitrary, but fixed, and define
\begin{equation}\label{eq:laplacePrinciples:tailset}
  A_\epsilon:=\{x\in\R^d:\ \|x-x_\ast\|_2\ge\epsilon\}.
\end{equation}
Furthermore, to improve readability, we will write $s_f=s_f(\epsilon)$.

\textbf{$\eta^\alpha_\nu$ is well defined.} We have to show that $Z_\alpha$ is finite and positive.
Since $f\ge f_\ast$ and $f$ is measurable, the integrand $e^{-\alpha f}$ is measurable and bounded by $e^{-\alpha f_\ast}$, 
so $Z_\alpha\le e^{-\alpha f_\ast}<\infty$. 
To show positivity, restrict the integral to $\ball_{s_f}(x_\ast)$ and use \eqref{eq:laplacePrinciples:epsf} and \eqref{eq:laplacePrinciples:mass},
\begin{equation}\label{eq:laplacePrinciples:normaliser}
  Z_\alpha\;\ge\;\int_{\ball_{s_f}(x_\ast)}e^{-\alpha f(x)}\,\nu(dx)
  \;\ge\;e^{-\alpha\left(f_\ast+\frac12S_\epsilon\right)}\,\nu\big(\ball_{s_f}(x_\ast)\big)
  \;\ge\;e^{-\alpha\left(f_\ast+\frac12S_\epsilon\right)}\,C_{s_f}\;>\;0 .
\end{equation}
In particular, $\eta^\alpha_\nu$ in \eqref{eq:laplacePrinciples:eta} is a well-defined Borel probability measure.

\textbf{A density bound away from the minimiser.} 
Let $g:\R^d\to[0,\infty]$ be measurable. 
By \eqref{eq:laplacePrinciples:separation}, $e^{-\alpha f(x)}\le e^{-\alpha(f_\ast+S_\epsilon)}$ for all $x\in A_\epsilon$, so
\begin{equation}\label{eq:laplacePrinciples:numerator}
  \int_{A_\epsilon}g(x)\,e^{-\alpha f(x)}\,\nu(dx)
  \;\le\;e^{-\alpha(f_\ast+S_\epsilon)}\int_{A_\epsilon}g(x)\,\nu(dx).
\end{equation}
Combining this with the lower bound \eqref{eq:laplacePrinciples:normaliser} for $Z_\alpha$ results in
\begin{equation}\label{eq:laplacePrinciples:densitybound}
    \int_{A_\epsilon}g\,d\eta^\alpha_\nu
    = \frac{1}{Z_\alpha}\int_{A_\epsilon}g(x)\,e^{-\alpha f(x)}\,\nu(dx)
    \leq \frac{e^{-\alpha S_\epsilon/2}}{C_{s_f}}\int_{A_\epsilon}g\,d\nu.
\end{equation}

\textbf{Part (i).}
Choosing $g\equiv1$ in \eqref{eq:laplacePrinciples:densitybound} and using $\nu(A_\epsilon)\le1$ yields
\begin{equation*}
    \eta^\alpha_\nu\big(\{x\in\R^d:\ \|x-x_\ast\|_2\ge\epsilon\}\big) = \int_{A_\epsilon} d\eta^\alpha_\nu
    \;\le\;\frac{e^{-\alpha S_\epsilon/2}}{C_{s_f}}\int_{A_\epsilon}1\,d\nu
    \leq \frac{e^{-\alpha S_\epsilon/2}}{C_{s_f}},
\end{equation*}
so we established \eqref{eq:laplacePrinciples:tailmass}. 

Let $\delta>0$. 
Since $C_{s_f}>0$, $\alpha_0$ in \eqref{eq:laplacePrinciples:alpha0} is well defined and nonnegative.
Let now $\alpha\ge\alpha_0$.
If $\delta\,C_{s_f}<1$, then $\alpha_0=\frac{2}{S_\epsilon}\log\frac{1}{\delta\,C_{s_f}}$, and we get
\begin{equation}\label{eq:laplacePrinciples:alpha0check}
  \frac{e^{-\alpha S_\epsilon/2}}{C_{s_f}}
  \;\le\;\frac{e^{-\alpha_0 S_\epsilon/2}}{C_{s_f}}
  \;=\;\frac{\delta\,C_{s_f}}{C_{s_f}}
  \;=\;\delta .
\end{equation}
If $\delta\,C_{s_f}\ge1$, then $e^{-\alpha S_\epsilon/2}/C_{s_f}\le 1/C_{s_f}\le\delta$.
This shows \eqref{eq:laplacePrinciples:alpha0bound}, and together with \eqref{eq:laplacePrinciples:tailmass} that for $\alpha \geq \alpha_0$ the left-hand side of \eqref{eq:laplacePrinciples:tailmass} is bounded by $\delta$.
Furthermore, since $\delta>0$ was arbitrary, this also shows $\eta^\alpha_\nu(A_\epsilon)\to0$ as $\alpha\to\infty$.

\textbf{Part (ii).}
Let $p\in[1,\infty)$ be such that \eqref{eq:laplacePrinciples:moment} holds.
Choosing $g(x)=\|x-x_\ast\|_2^p$ in \eqref{eq:laplacePrinciples:densitybound} and using \eqref{eq:laplacePrinciples:moment} gives
\begin{equation}\label{eq:laplacePrinciples:tailmoment}
  \int_{A_\epsilon}\|x-x_\ast\|_2^p\,\eta^\alpha_\nu(dx)
  \;\le\;\frac{e^{-\alpha S_\epsilon/2}}{C_{s_f}}\int_{A_\epsilon}\|x-x_\ast\|_2^p\,\nu(dx)
  \;\le\;\frac{e^{-\alpha S_\epsilon/2}}{C_{s_f}}\,M_p .
\end{equation}
On the complement $A_\epsilon^c$ we have $\|x-x_\ast\|_2^p<\epsilon^p$, and since $\eta^\alpha_\nu$ is a probability measure we also have $\eta^\alpha_\nu(A_\epsilon^c)\le1$.
This implies that the integral over $A_\epsilon^c$ is at most $\epsilon^p$.
Together with \eqref{eq:laplacePrinciples:tailmoment} we get
\begin{equation}\label{eq:laplacePrinciples:split}
  \E_{X\distr\eta^\alpha_\nu}\big[\|X-x_\ast\|_2^p\big]
  =\int_{A_\epsilon^c}\|x-x_\ast\|_2^p\,\eta^\alpha_\nu(dx)
   +\int_{A_\epsilon}\|x-x_\ast\|_2^p\,\eta^\alpha_\nu(dx)
  \;\le\;\epsilon^p+\frac{e^{-\alpha S_\epsilon/2}}{C_{s_f}}\,M_p<\infty ,
\end{equation}
which is \eqref{eq:laplacePrinciples:momentbound}.
In particular, since $\|x\|_2^p\le2^{p-1}\big(\|x_\ast\|_2^p+\|x-x_\ast\|_2^p\big)$ by convexity of $t\mapsto t^p$ on $[0,\infty)$, we have $\E_{X\distr\eta^\alpha_\nu}[\|X\|_2^p]<\infty$.

It remains to show convergence.
Let $\gamma>0$ and set $\epsilon:=(\gamma/2)^{1/p}$.
If $M_p=0$, then \eqref{eq:laplacePrinciples:momentbound} gives $\E_{X\distr\eta^\alpha_\nu}[\|X-x_\ast\|_2^p]\le\gamma/2$ for all $\alpha>0$.
Otherwise, \eqref{eq:laplacePrinciples:alpha0bound} with $\delta=\gamma/(2M_p)$ shows that
\begin{equation*}
     \frac{e^{-\alpha S_\epsilon/2}}{C_{s_f(\epsilon)}}\,M_p \;\leq\; \frac{\gamma}{2}
     \qquad\text{for all }\alpha\geq\alpha_0\Big(\frac{\gamma}{2M_p},S_\epsilon,C_{s_f(\epsilon)}\Big),
\end{equation*}
and \eqref{eq:laplacePrinciples:momentbound} gives $\E_{X\distr\eta^\alpha_\nu}[\|X-x_\ast\|_2^p]\le\epsilon^p+\gamma/2=\gamma$ for all these $\alpha$.
Since $\gamma>0$ was arbitrary, $\E_{X\distr\eta^\alpha_\nu}[\|X-x_\ast\|_2^p]\to0$ as $\alpha\to\infty$.

Finally, let \eqref{eq:laplacePrinciples:moment} hold for $p=1$.
Then $\E_{X\distr\eta^\alpha_\nu}[\|X\|_2]<\infty$, so the mean $\E_{X\distr\eta^\alpha_\nu}[X]$ exists.
Jensen's inequality and \eqref{eq:laplacePrinciples:momentbound} with $p=1$ give
\begin{equation}\label{eq:laplacePrinciples:jensen}
  \big\|\E_{X\distr\eta^\alpha_\nu}[X]-x_\ast\big\|_2
  =\big\|\E_{X\distr\eta^\alpha_\nu}[X-x_\ast]\big\|_2
  \;\le\;\E_{X\distr\eta^\alpha_\nu}\big[\|X-x_\ast\|_2\big]
  \;\le\; \epsilon+\frac{e^{-\alpha S_\epsilon/2}}{C_{s_f}}\,M_1,
\end{equation}
which is \eqref{eq:laplacePrinciples:mean}. Since the middle term in \eqref{eq:laplacePrinciples:jensen} tends to $0$ as $\alpha\to\infty$, so does the left-hand side.
\end{proof}
Let us briefly discuss the differences to \cite[Lemma~3.3]{fornasier2026consensus} and \cite[Theorem~4.3]{huang2026faithful}.
Assumption \eqref{eq:laplacePrinciples:mass} is also present in these references and used in the same way to get a lower bound on the normalisation constant $Z_\alpha$.
Assumption \eqref{eq:laplacePrinciples:separation} is implied by the setting of these references, cf.\ \cite[Remark~4.2]{huang2026faithful}. 
The splitting technique in our proof is essentially the approach from the proof of \cite[Theorem~4.3]{huang2026faithful}.

We now show that $\Norm(x,\Sigma)$ for any $x\in\R^d$ is a suitable probability measure and that we can apply a quantitative Laplace principle, cf.\ Lemma~\ref{lem:laplacePrinciples}.
This result is essentially the first step of the proof of \cite[Theorem~3.5]{fornasier2026consensus}, generalised from a scaled identity matrix to a general positive definite covariance matrix.
%
%
\begin{lemma} \label{lem:normalSuitableForLaplacePrinciple}
    Let $x_\ast\in\R^d$, $r>0$, and let $\Sigma\in\R^{d\times d}$ be symmetric positive definite.
    For all $x\in\ball_r(x_\ast)$ and all $s\in\Rp$ it holds
    \begin{equation}
        \Norm(x, \Sigma)(\ball_s(x_\ast)) \geq \frac{\pi^{d/2}s^d}{\Gamma\!\left(\frac d2+1\right)}\,(2\pi\sigma_{\max}(\Sigma)^2)^{-d/2}\exp\!\left(-\frac{(r+s)^2}{2\sigma_{\min}(\Sigma)^2}\right)
        =: C_s>0.
    \end{equation}
    Furthermore, for $X\distr\Norm(x,\Sigma)$ we also have
    \begin{align}
        \E\big[\|X-x_\ast\|_2\big] & \leq  r+\sigma_{\max}(\Sigma) \sqrt d =: M_1 \\
        \E\big[\|X-x_\ast\|_2^2\big] & \leq  r^2+d\,\sigma_{\max}(\Sigma)^2 =: M_2.
    \end{align}
\end{lemma}
\begin{proof}
Let $x\in \ball_r(x_\ast)$ be arbitrary.
Since $\Sigma$ is symmetric positive definite, we have $\det(\Sigma) \leq \lambda_{\max}(\Sigma)^{d} = \sigma_{\max}(\Sigma)^{2d}$, and for all $v\in\R^d$ we have $v^\top \Sigma^{-1}v \leq \lambda_{\min}(\Sigma)^{-1}\|v\|_2^2 = \sigma_{\min}(\Sigma)^{-2}\|v\|_2^2$.

Let now $s\in\Rp$ and $y\in\ball_s(x_\ast)$ be arbitrary.
Since $\|y-x\|_2\le\|y-x_\ast\|_2+\|x_\ast-x\|_2\le s+r$, we have
\begin{align*}
    \Norm(y \mid x,\Sigma) & = \frac{1}{\sqrt{(2\pi)^d \det(\Sigma)}}\exp\left(-\frac12 (y-x)^\top \Sigma^{-1}(y-x) \right) \\
    & \geq (2\pi\sigma_{\max}(\Sigma)^2)^{-d/2}\exp\!\left(-\frac{(r+s)^2}{2\sigma_{\min}(\Sigma)^2}\right).
\end{align*}
Integrating over $\ball_s(x_\ast)$ and using $\mathrm{Vol}(\ball_s(x_\ast))=\pi^{d/2}s^d/\Gamma(\frac d2+1)$ leads to
\begin{align*}
    \Norm(x, \Sigma)(\ball_s(x_\ast)) & = \int_{\ball_s(x_\ast)} \frac{1}{\sqrt{(2\pi)^d \det(\Sigma)}}\exp\left(-\frac12 (y-x)^\top \Sigma^{-1}(y-x) \right) \mathrm{d}y \\
    & \geq (2\pi\sigma_{\max}(\Sigma)^2)^{-d/2}\exp\!\left(-\frac{(r+s)^2}{2\sigma_{\min}(\Sigma)^2}\right) \mathrm{Vol}(\ball_s(x_\ast)) \\
    & = \frac{\pi^{d/2}s^d}{\Gamma\!\left(\frac d2+1\right)}\,(2\pi\sigma_{\max}(\Sigma)^2)^{-d/2}\exp\!\left(-\frac{(r+s)^2}{2\sigma_{\min}(\Sigma)^2}\right) > 0.
\end{align*}
For the moment bounds, note that with $X=x+Z$ with $Z\distr\Norm(0,\Sigma)$ we have $X \distr \Norm(x,\Sigma)$.
By the triangle inequality, Jensen's inequality, $\E\big[\|Z\|_2^2\big]=\trace(\Sigma)$ and $\trace(\Sigma)\leq d\,\lambda_{\max}(\Sigma)$, we get
\begin{align*}
    \E\big[\|X-x_\ast\|_2\big] & \;\le\;\|x-x_\ast\|_2+\E\big[\|Z\|_2\big] \\
    & \;\le\; r+\big(\E\big[\|Z\|_2^2\big]\big)^{1/2} \\
    & \;=\; r+(\trace(\Sigma))^{1/2} \\
    & \;\le\; r+\sigma_{\max}(\Sigma)\sqrt d=M_1 .
\end{align*}
For the second moment, expanding the square and using $\E[Z]=0$ gives
\begin{align*}
    \E\big[\|X-x_\ast\|_2^2\big]
    & = \|x-x_\ast\|_2^2 + 2(x-x_\ast)^\top\E[Z] + \E\big[\|Z\|_2^2\big] \\
    & = \|x-x_\ast\|_2^2 + \trace(\Sigma)
    \;\le\; r^2 + d\,\sigma_{\max}(\Sigma)^2 = M_2 .
\end{align*}
\end{proof}

\subsubsection{Mean-Field Transition Map as a Contraction} \label{sec:transitionMapAsContraction}
Our next goal is to show that the transition map $\mathcal{T}_\alpha$ is a contraction for $\alpha$ large enough.
We follow the strategy of the proof of \cite[Theorem~3.5]{fornasier2026consensus} and bound the norm of the gradient of $\mathcal{T}_\alpha$.
For this, we need the following technical result, which is contained, under slightly stronger assumptions, in Step 2 of the aforementioned proof,
and can be proved using exactly the same arguments.
%
%
\begin{lemma}[Differentiability of the transition map]\label{lem:consensusMapDifferentiable}
Let $f:\R^d\to\R$ be Borel measurable and bounded below, $\underline f:=\inf_{y\in\R^d}f(y)>-\infty$.
Let $\Sigma\in\R^{d\times d}$ be symmetric positive definite and $\alpha>0$. 
\begin{enumerate}
    \item[(i)] For every $x\in\R^d$ we have $0<Z_\alpha(\Norm(x,\Sigma))\le e^{-\alpha\underline f}$. Moreover, the Gibbs tilt $\eta^\alpha_{\Norm(x,\Sigma)}$, defined as in \eqref{eq:laplacePrinciples:eta} with $\nu=\Norm(x,\Sigma)$, is a Borel probability measure with finite moments of all orders, and for every Borel measurable $h:\R^d\to[0,\infty]$,
    \begin{equation}\label{eq:consensusMapDifferentiable:changeofmeasure}
        \E_{Y\distr\eta^\alpha_{\Norm(x,\Sigma)}}\big[h(Y)\big]=\frac{1}{Z_\alpha(\Norm(x,\Sigma))}\,\E_{X\distr\Norm(x,\Sigma)}\big[h(X)\,e^{-\alpha f(X)}\big].
    \end{equation}
    In particular, $\mathcal{T}_\alpha(x)$ is well defined and coincides with the mean $\E_{Y\distr\eta^\alpha_{\Norm(x,\Sigma)}}[Y]$ of the Gibbs tilt.
    \item[(ii)] $\mathcal{T}_\alpha\in C^1(\R^d;\R^d)$, and its Jacobian $\nabla \mathcal{T}_\alpha(x)=\big(\partial_{x_j}\mathcal{T}_{\alpha,i}(x)\big)_{i,j}$ is
    \begin{equation}\label{eq:consensusMapDifferentiable:jacobian}
        \nabla \mathcal{T}_\alpha(x)=\mathrm{Cov}\big(\eta^\alpha_{\Norm(x,\Sigma)}\big)\,\Sigma^{-1}
        \qquad\text{for all }x\in\R^d,
    \end{equation}
    where $\mathrm{Cov}(\eta^\alpha_{\Norm(x,\Sigma)}):=\E_{Y\distr\eta^\alpha_{\Norm(x,\Sigma)}}\big[(Y-\mathcal{T}_\alpha(x))(Y-\mathcal{T}_\alpha(x))^\top\big]$.
\end{enumerate}
\end{lemma}

We can now bound the norm of the gradient of $\mathcal{T}_\alpha$.
%
%
\begin{lemma}[Jacobian bound]\label{lem:consensusMapJacobianBound}
Let $f$ satisfy Assumption~\ref{ass:objective}, let $\Sigma\in\R^{d\times d}$ be symmetric positive definite, and let $r>0$.
Let $C_s$ and $M_2$ be the constants of Lemma~\ref{lem:normalSuitableForLaplacePrinciple} for $r$ and $\Sigma$, and let $s_f(\epsilon)$ be as in \eqref{eq:laplacePrinciples:epsf}.
Then, with $\mathcal{T}_\alpha$ as in Lemma~\ref{lem:consensusMapDifferentiable}, for all $\epsilon>0$, $\alpha>0$ and $x\in\ball_r(x_\ast)$,
\begin{equation}\label{eq:consensusMapJacobianBound:bound}
    \big\|\nabla \mathcal{T}_\alpha(x)\big\|
    \;\le\;L_\alpha(\epsilon)
    :=\frac{1}{\sigma_{\min}(\Sigma)^2}\left(\epsilon^2+\frac{e^{-\alpha S_\epsilon/2}}{C_{s_f(\epsilon)}}\,M_2\right).
\end{equation}
\end{lemma}

\begin{proof}
By (F1), $f$ is measurable and bounded below by $f_\ast$, so Lemma~\ref{lem:consensusMapDifferentiable} applies.
Let $x\in\ball_r(x_\ast)$, $\epsilon>0$ and $\alpha>0$ be arbitrary. 
Using \eqref{eq:consensusMapDifferentiable:jacobian}, we can start with
\begin{equation*}
    \big\|\nabla \mathcal{T}_\alpha(x)\big\|
    \le\big\|\mathrm{Cov}(\eta^\alpha_{\Norm(x,\Sigma)})\big\| \,\|\Sigma^{-1}\|
    \le\frac{\trace\,\mathrm{Cov}(\eta^\alpha_{\Norm(x,\Sigma)})}{\sigma_{\min}(\Sigma)^2} .
\end{equation*}
Let now $Y\distr\eta^\alpha_{\Norm(x,\Sigma)}$. Then we have
\begin{align*}
    \trace\,\mathrm{Cov}(\eta^\alpha_{\Norm(x,\Sigma)}) = \E[\|Y - \mathcal{T}_\alpha(x)\|_2^2]
    \leq \E\big[\|Y-x_\ast\|_2^2\big],
\end{align*}
where we used the definition of the covariance matrix and of the trace together with the fact $\E[Y]=\mathcal{T}_\alpha(x)$ in the first step,
and in the second step we used the well-known identity $\E\|Y-\mathcal{T}_\alpha(x)\|_2^2=\E\|Y-z\|_2^2-\|\mathcal{T}_\alpha(x)-z\|_2^2$,
which holds for any $z\in\R^d$.

We know from Lemma~\ref{lem:normalSuitableForLaplacePrinciple} that $\Norm(x,\Sigma)$ satisfies \eqref{eq:laplacePrinciples:mass} with the constants $C_s$ and \eqref{eq:laplacePrinciples:moment} for $p=2$ with $M_2$, so Lemma~\ref{lem:laplacePrinciples}(ii) gives
\begin{equation*}
     \E\big[\|Y-x_\ast\|_2^2\big]\;\le\;\epsilon^2+\frac{e^{-\alpha S_\epsilon/2}}{C_{s_f(\epsilon)}}\,M_2.
\end{equation*}
Altogether, we have
\begin{equation*}
     \big\|\nabla \mathcal{T}_\alpha(x)\big\| \leq \frac{\trace\,\mathrm{Cov}(\eta^\alpha_{\Norm(x,\Sigma)})}{\sigma_{\min}(\Sigma)^2} 
     \leq \frac{1}{\sigma_{\min}(\Sigma)^2}\left(\epsilon^2+\frac{e^{-\alpha S_\epsilon/2}}{C_{s_f(\epsilon)}}\,M_2\right).
\end{equation*}
\end{proof}
As an immediate corollary of the previous bound, we find that $\mathcal{T}_\alpha$ is a contraction for all $\alpha$ large enough.
%
%
\begin{corollary}[Contraction of the transition map]\label{cor:consensusMapContraction}
Let the assumptions of Lemma~\ref{lem:consensusMapJacobianBound} hold, let $L\in(0,1)$, and set $\epsilon:=\sigma_{\min}(\Sigma)\sqrt{L/2}$.
With $\alpha_0$ from \eqref{eq:laplacePrinciples:alpha0}, we have
\begin{equation}\label{eq:consensusMapContraction:choice}
    \alpha_0\Big(\frac{L\,\sigma_{\min}(\Sigma)^2}{2M_2},\,S_\epsilon,\,C_{s_f(\epsilon)}\Big)
    =\frac{2}{S_\epsilon}\max\left\{0,\;\log\frac{2M_2}{L\,\sigma_{\min}(\Sigma)^2\,C_{s_f(\epsilon)}}\right\}.
\end{equation}
Then, for all $\alpha\ge\alpha_0\big(\frac{L\sigma_{\min}(\Sigma)^2}{2M_2},S_\epsilon,C_{s_f(\epsilon)}\big)$,
\begin{equation}\label{eq:consensusMapContraction:supbound}
    \sup_{x\in\ball_r(x_\ast)}\big\|\nabla \mathcal{T}_\alpha(x)\big\|\;\le\;L ,
\end{equation}
and $\mathcal{T}_\alpha$ is Lipschitz continuous on $\ball_r(x_\ast)$ with constant $L$, i.e.
\begin{equation}\label{eq:consensusMapContraction:lipschitz}
    \big\|\mathcal{T}_\alpha(x)-\mathcal{T}_\alpha(x')\big\|_2\;\le\;L\,\|x-x'\|_2
    \qquad\text{for all }x,x'\in\ball_r(x_\ast).
\end{equation}
In particular, $\sup_{x\in\ball_r(x_\ast)}\|\nabla \mathcal{T}_\alpha(x)\|\to0$ as $\alpha\to\infty$.
\end{corollary}

\begin{proof}
With $\epsilon=\sigma_{\min}(\Sigma)\sqrt{L/2}$, the first term of $L_\alpha(\epsilon)$ in \eqref{eq:consensusMapJacobianBound:bound} equals $\epsilon^2/\sigma_{\min}(\Sigma)^2=L/2$. For $\alpha\ge\alpha_0\big(\frac{L\sigma_{\min}(\Sigma)^2}{2M_2},S_\epsilon,C_{s_f(\epsilon)}\big)$, \eqref{eq:laplacePrinciples:alpha0bound} with $\delta=L\sigma_{\min}(\Sigma)^2/(2M_2)$ gives
\begin{equation}\label{eq:consensusMapContraction:secondterm}
    \frac{1}{\sigma_{\min}(\Sigma)^2}\,\frac{e^{-\alpha S_\epsilon/2}}{C_{s_f(\epsilon)}}\,M_2
    \;\le\;\frac{1}{\sigma_{\min}(\Sigma)^2}\,\frac{L\,\sigma_{\min}(\Sigma)^2}{2M_2}\,M_2
    \;=\;\frac L2 .
\end{equation}
Hence $L_\alpha(\epsilon)\le L$, and Lemma~\ref{lem:consensusMapJacobianBound} yields \eqref{eq:consensusMapContraction:supbound}.

Let $x,x'\in\ball_r(x_\ast)$. Since $\ball_r(x_\ast)$ is convex, $x'+t(x-x')\in\ball_r(x_\ast)$ for all $t\in[0,1]$. Since $\mathcal{T}_\alpha\in C^1$ by Lemma~\ref{lem:consensusMapDifferentiable}, the fundamental theorem of calculus and \eqref{eq:consensusMapContraction:supbound} give
\begin{equation}\label{eq:consensusMapContraction:mvt}
    \big\|\mathcal{T}_\alpha(x)-\mathcal{T}_\alpha(x')\big\|_2
    =\Big\|\int_0^1\nabla \mathcal{T}_\alpha\big(x'+t(x-x')\big)(x-x')\,\mathrm dt\Big\|_2
    \le L\,\|x-x'\|_2 .
\end{equation}
Since $L\in(0,1)$ was arbitrary, the final claim follows.
\end{proof}
We also need to show that a region near the global minimiser is invariant under the transition map for $\alpha$ large enough.
\begin{lemma}[Invariance of balls around the minimiser]\label{lem:consensusMapInvariance}
Let $f$ satisfy Assumption~\ref{ass:objective}, let $\Sigma\in\R^{d\times d}$ be symmetric positive definite, and let $r>0$.
Let $C_s$ and $M_1$ be the constants of Lemma~\ref{lem:normalSuitableForLaplacePrinciple} for $r$ and $\Sigma$, and let $s_f(\epsilon)$ be as in \eqref{eq:laplacePrinciples:epsf}. Then the following hold.
\begin{enumerate}
    \item[(i)] For all $\epsilon>0$, $\alpha>0$ and $x\in\ball_r(x_\ast)$,
    \begin{equation}\label{eq:consensusMapInvariance:bound}
        \big\|\mathcal{T}_\alpha(x)-x_\ast\big\|_2\;\le\;\epsilon+\frac{e^{-\alpha S_\epsilon/2}}{C_{s_f(\epsilon)}}\,M_1 .
    \end{equation}
    \item[(ii)] Let $\rho>0$. With $\alpha_0$ from \eqref{eq:laplacePrinciples:alpha0}, we have
    \begin{equation}\label{eq:consensusMapInvariance:alpha0}
        \alpha_0\Big(\frac{\rho}{4M_1},\,S_{\rho/2},\,C_{s_f(\rho/2)}\Big)=\frac{2}{S_{\rho/2}}\max\left\{0,\;\log\frac{4M_1}{\rho\,C_{s_f(\rho/2)}}\right\},
    \end{equation}
    and for all $\alpha\ge\alpha_0\big(\frac{\rho}{4M_1},S_{\rho/2},C_{s_f(\rho/2)}\big)$,
    \begin{equation}\label{eq:consensusMapInvariance:inclusion}
        \sup_{x\in\ball_r(x_\ast)}\big\|\mathcal{T}_\alpha(x)-x_\ast\big\|_2\;\le\;\tfrac34\rho ,
        \qquad\text{in particular}\qquad
        \mathcal{T}_\alpha\big(\ball_r(x_\ast)\big)\subseteq\ball_\rho(x_\ast).
    \end{equation}
    \item[(iii)] In particular, for $\rho\in(0,r]$ and all $\alpha\ge\alpha_0\big(\frac{\rho}{4M_1},S_{\rho/2},C_{s_f(\rho/2)}\big)$, the ball $\ball_r(x_\ast)$ is invariant under $\mathcal{T}_\alpha$:
    \begin{equation}\label{eq:consensusMapInvariance:invariance}
        \mathcal{T}_\alpha\big(\ball_r(x_\ast)\big)\subseteq\ball_\rho(x_\ast)\subseteq\ball_r(x_\ast).
    \end{equation}
    Moreover, $\sup_{x\in\ball_r(x_\ast)}\|\mathcal{T}_\alpha(x)-x_\ast\|_2\to0$ as $\alpha\to\infty$.
\end{enumerate}
\end{lemma}

\begin{proof}
\emph{Part (i).}
Let $x\in\ball_r(x_\ast)$ and write $\nu:=\Norm(x,\Sigma)$ for brevity. 
By Lemma~\ref{lem:normalSuitableForLaplacePrinciple}, $\nu$ satisfies \eqref{eq:laplacePrinciples:mass} with the constants $C_s$, and \eqref{eq:laplacePrinciples:moment} for $p=1$ with $M_1$. 
Both constants depend only on $r$, $\Sigma$ and $d$, not on $x$. 
By Lemma~\ref{lem:laplacePrinciples}(ii) with $p=1$, more precisely \eqref{eq:laplacePrinciples:mean}, for all $\epsilon>0$ and $\alpha>0$,
\begin{equation}\label{eq:consensusMapInvariance:laplace}
    \big\|\mathcal{T}_\alpha(x)-x_\ast\big\|_2
    =\big\|\E_{X\distr\eta^\alpha_\nu}[X]-x_\ast\big\|_2
    \;\le\;\epsilon+\frac{e^{-\alpha S_\epsilon/2}}{C_{s_f(\epsilon)}}\,M_1 ,
\end{equation}
which is \eqref{eq:consensusMapInvariance:bound}.

\emph{Part (ii).}
Let $\rho>0$ and $\epsilon:=\rho/2$. For $\alpha\ge\alpha_0\big(\frac{\rho}{4M_1},S_{\rho/2},C_{s_f(\rho/2)}\big)$, \eqref{eq:laplacePrinciples:alpha0bound} with $\delta=\rho/(4M_1)$ gives
\begin{equation}\label{eq:consensusMapInvariance:secondterm}
    \frac{e^{-\alpha S_{\rho/2}/2}}{C_{s_f(\rho/2)}}\,M_1
    \;\le\;\frac{\rho}{4M_1}\,M_1
    \;=\;\frac{\rho}{4} .
\end{equation}
Inserting \eqref{eq:consensusMapInvariance:secondterm} into \eqref{eq:consensusMapInvariance:bound} with $\epsilon=\rho/2$ yields $\|\mathcal{T}_\alpha(x)-x_\ast\|_2\le\frac\rho2+\frac\rho4=\frac34\rho<\rho$ for all $x\in\ball_r(x_\ast)$. This proves \eqref{eq:consensusMapInvariance:inclusion}.

\emph{Part (iii).}
For $\rho\in(0,r]$, \eqref{eq:consensusMapInvariance:invariance} follows from \eqref{eq:consensusMapInvariance:inclusion} and $\ball_\rho(x_\ast)\subseteq\ball_r(x_\ast)$. Since $\rho>0$ was arbitrary in (ii), $\sup_{x\in\ball_r(x_\ast)}\|\mathcal{T}_\alpha(x)-x_\ast\|_2\le\frac34\rho$ for all $\alpha\ge\alpha_0\big(\frac{\rho}{4M_1},S_{\rho/2},C_{s_f(\rho/2)}\big)$, and hence the supremum tends to $0$ as $\alpha\to\infty$.
\end{proof}

Contraction and invariance alone do not quantify the distance of the iterates to $x_\ast$ independently of the radius.
For this, we evaluate the bias of the transition map at the minimiser itself, where we can use the constants of Lemma~\ref{lem:normalSuitableForLaplacePrinciple} with radius $0$,
\begin{equation}\label{eq:biasAtMinimiser:constants}
    C^\circ_s:=\frac{\pi^{d/2}s^d}{\Gamma\!\left(\frac d2+1\right)}\,\big(2\pi\sigma_{\max}(\Sigma)^2\big)^{-d/2}\exp\!\left(-\frac{s^2}{2\sigma_{\min}(\Sigma)^2}\right),
    \qquad
    M^\circ_1:=\sigma_{\max}(\Sigma)\sqrt d ,
\end{equation}
which depend only on $\Sigma$, $d$ and $s>0$.
%
%
\begin{lemma}[Bias at the minimiser]\label{lem:biasAtMinimiser}
Let $f$ satisfy Assumption~\ref{ass:objective}, and let $\Sigma\in\R^{d\times d}$ be symmetric positive definite. For $\beta>0$, define
\begin{equation}\label{eq:biasAtMinimiser:alphab}
    \alpha_{\mathrm b}(f,\Sigma,\beta):=\alpha_0\Big(\frac{\beta}{2M^\circ_1},\,S_{\beta/2},\,C^\circ_{s_f(\beta/2)}\Big)
\end{equation}
with $\alpha_0$ from \eqref{eq:laplacePrinciples:alpha0}. Then $\big\|\mathcal{T}_\alpha(x_\ast)-x_\ast\big\|_2\le\beta$ for all $\alpha>0$ with $\alpha\ge\alpha_{\mathrm b}(f,\Sigma,\beta)$.
\end{lemma}
\begin{proof}
Let $\nu^\circ:=\Norm(x_\ast,\Sigma)$. For every $r'>0$, Lemma~\ref{lem:normalSuitableForLaplacePrinciple} with radius $r'$ and $x=x_\ast\in\ball_{r'}(x_\ast)$ gives lower bounds on $\nu^\circ(\ball_s(x_\ast))$ and upper bounds on $\E_{X\distr\nu^\circ}[\|X-x_\ast\|_2]$, which depend continuously on $r'$, while the left-hand sides do not depend on $r'$. Letting $r'\to0$ yields $\nu^\circ(\ball_s(x_\ast))\ge C^\circ_s>0$ for all $s>0$ and $\E_{X\distr\nu^\circ}[\|X-x_\ast\|_2]\le M^\circ_1$.
Hence Lemma~\ref{lem:laplacePrinciples}(ii) with $p=1$ applies to $\nu^\circ$ with the constants $C^\circ_s$ and $M^\circ_1$. By Lemma~\ref{lem:consensusMapDifferentiable}(i), $\mathcal{T}_\alpha(x_\ast)$ is the mean of $\eta^\alpha_{\nu^\circ}$, so \eqref{eq:laplacePrinciples:mean} with $\epsilon=\beta/2$ gives
\begin{equation*}
    \big\|\mathcal{T}_\alpha(x_\ast)-x_\ast\big\|_2
    \;\le\;\frac{\beta}{2}+\frac{e^{-\alpha S_{\beta/2}/2}}{C^\circ_{s_f(\beta/2)}}\,M^\circ_1 .
\end{equation*}
For $\alpha\ge\alpha_{\mathrm b}(f,\Sigma,\beta)$, \eqref{eq:laplacePrinciples:alpha0bound} with $\delta=\beta/(2M^\circ_1)$ bounds the second term by $\beta/2$.
\end{proof}

\subsubsection{Main Results: Proofs and Discussion} \label{sec:suppl:mainresults}
We now prove the results of Section~\ref{sec:meanField}.
First, we establish an invariance and contraction property of the mean-field transition map.
This result follows from the results above using essentially the same arguments as in \cite[Section~4]{fornasier2026consensus}.
\begin{proof}[Proof of Theorem~\ref{thm:meanFieldConvergence}]
Let $r>0$ and $\rho\in(0,1)$, let $C_s$ and $M_2$ be the constants of Lemma~\ref{lem:normalSuitableForLaplacePrinciple} for $r$ and $\Sigma$, let $\epsilon_\rho:=\sigma_{\min}(\Sigma)\sqrt{\rho/2}$, and define
\begin{equation}\label{eq:meanFieldConvergence:alpha0}
    \alpha_0(f,\Sigma,r,\rho):=\max\Big\{\alpha_0\Big(\frac{\rho\,\sigma_{\min}(\Sigma)^2}{2M_2},\,S_{\epsilon_\rho},\,C_{s_f(\epsilon_\rho)}\Big),\ \alpha_{\mathrm b}\big(f,\Sigma,(1-\rho)r\big)\Big\},
\end{equation}
with $\alpha_0(\cdot,\cdot,\cdot)$ from \eqref{eq:laplacePrinciples:alpha0} and $\alpha_{\mathrm b}$ from \eqref{eq:biasAtMinimiser:alphab}. It depends only on $f$ (through $S_\epsilon$ and $s_f(\epsilon)$), on $\Sigma$, $r$, $\rho$ and the dimension $d$. Let $\alpha>0$ with $\alpha\ge\alpha_0(f,\Sigma,r,\rho)$.

\emph{Ad (i).}
The first threshold in \eqref{eq:meanFieldConvergence:alpha0} is \eqref{eq:consensusMapContraction:choice} with $L=\rho$, so Corollary~\ref{cor:consensusMapContraction} yields \eqref{eq:meanFieldConvergence:contraction}.
By Lemma~\ref{lem:biasAtMinimiser}, $\|\mathcal{T}_\alpha(x_\ast)-x_\ast\|_2\le(1-\rho)r$. Hence, for $x\in\ball_r(x_\ast)$,
\begin{equation}\label{eq:meanFieldConvergence:invariance}
    \big\|\mathcal{T}_\alpha(x)-x_\ast\big\|_2
    \;\le\;\big\|\mathcal{T}_\alpha(x)-\mathcal{T}_\alpha(x_\ast)\big\|_2+\big\|\mathcal{T}_\alpha(x_\ast)-x_\ast\big\|_2
    \;\le\;\rho r+(1-\rho)r=r ,
\end{equation}
i.e.\ $\mathcal{T}_\alpha(\ball_r(x_\ast))\subseteq\ball_r(x_\ast)$.

\emph{Ad (ii).}
The closed ball $\ball_r(x_\ast)$ is compact, $\mathcal{T}_\alpha$ maps it into itself by (i), and it is a contraction on it by \eqref{eq:meanFieldConvergence:contraction}. Banach's fixed-point theorem yields a unique fixed point $m_\alpha\in\ball_r(x_\ast)$.
Let $x_0\in\ball_r(x_\ast)$. By (i) and induction, $x_k\in\ball_r(x_\ast)$ for all $k\in\Nnn$. Since $m_\alpha=\mathcal{T}_\alpha(m_\alpha)$, \eqref{eq:meanFieldConvergence:contraction} gives $\|x_{k+1}-m_\alpha\|_2\le\rho\|x_k-m_\alpha\|_2$, and induction yields \eqref{eq:meanFieldConvergence:rate}.
Finally, since $m_\alpha\in\ball_r(x_\ast)$ and $m_\alpha=\mathcal{T}_\alpha(m_\alpha)$, we have
\begin{equation*}
    \|m_\alpha-x_\ast\|_2\;\le\;\sup_{x\in\ball_r(x_\ast)}\big\|\mathcal{T}_\alpha(x)-x_\ast\big\|_2 ,
\end{equation*}
and the right-hand side tends to $0$ as $\alpha\to\infty$ by Lemma~\ref{lem:consensusMapInvariance}(iii).
\end{proof}
We can now turn to Theorem \ref{thm:meanFieldOneStep}.
The first part is similar to \cite[Theorem~4.4]{fornasier2026consensus}, though we show it under weaker assumptions with a more direct argument.
The second part has been anticipated in \cite[Remark~4.5]{fornasier2026consensus}, but not formally stated and proven yet.
\begin{proof}[Proof of Theorem~\ref{thm:meanFieldOneStep}]
Let $\beta>0$, let $\alpha_{\mathrm b}(f,\Sigma,\beta)$ be as in \eqref{eq:biasAtMinimiser:alphab}, which only depends on $f$, $\Sigma$ and $\beta$, and let $r>0$, $\rho\in(0,1)$ and $\alpha\ge\max\{\alpha_0(f,\Sigma,r,\rho),\alpha_{\mathrm b}(f,\Sigma,\beta)\}$.
Set $\beta_r:=\min\{\beta,(1-\rho)r\}$. Since $\alpha_{\mathrm b}(f,\Sigma,\beta_r)$ equals either $\alpha_{\mathrm b}(f,\Sigma,\beta)$ or $\alpha_{\mathrm b}(f,\Sigma,(1-\rho)r)\le\alpha_0(f,\Sigma,r,\rho)$, Lemma~\ref{lem:biasAtMinimiser} gives $\|\mathcal{T}_\alpha(x_\ast)-x_\ast\|_2\le\beta_r\le\beta$.

\emph{Ad (i).}
Let $x\in\ball_r(x_\ast)$. Since also $x_\ast\in\ball_r(x_\ast)$, \eqref{eq:meanFieldConvergence:contraction} gives
\begin{equation*}
    \big\|\mathcal{T}_\alpha(x)-x_\ast\big\|_2
    \;\le\;\big\|\mathcal{T}_\alpha(x)-\mathcal{T}_\alpha(x_\ast)\big\|_2+\big\|\mathcal{T}_\alpha(x_\ast)-x_\ast\big\|_2
    \;\le\;\rho\,\|x-x_\ast\|_2+\beta .
\end{equation*}

\emph{Ad (ii).}
By Theorem~\ref{thm:meanFieldConvergence}(ii), $x_k\in\ball_r(x_\ast)$ for all $k\in\Nnn$. Writing $e_k:=\|x_k-x_\ast\|_2$, (i) gives $e_{k+1}\le\beta+\rho\,e_k$, and induction yields
\begin{equation*}
    e_k\;\le\;\rho^ke_0+\beta\sum_{j=0}^{k-1}\rho^j\;\le\;\rho^ke_0+\frac{\beta}{1-\rho} .
\end{equation*}
\end{proof}
As an immediate consequence, we get Corollary \ref{cor:meanFieldAccuracy}.
\begin{proof}[Proof of Corollary~\ref{cor:meanFieldAccuracy}]
Apply Theorem~\ref{thm:meanFieldOneStep} with $\rho=\frac12$ and $\beta=\frac\epsilon4$. For $x_0\in\ball_r(x_\ast)$ and $k\ge K_0(r,\epsilon)$, \eqref{eq:meanFieldOneStep:iterates} gives $\|x_k-x_\ast\|_2\le2^{-k}r+\frac\epsilon2\le\epsilon$, since $2^k\ge2r/\epsilon$ if $K_0(r,\epsilon)>0$, and $r\le\epsilon/2$ otherwise.
\end{proof}
Finally, let us discuss in detail how the preceding results show that iterating in mean-field MPPI can indeed improve over just a single step of MPPI.
This phenomenon has already been described in \cite[Remark~4.5]{fornasier2026consensus}, but using two different thresholds for the inverse temperature allows us to discuss this quantitatively.
\begin{remark}[Iterations help]\label{rem:meanFieldIterationsHelp}
To reach accuracy $\epsilon$ in a single step from every point of $\ball_r(x_\ast)$, Lemma~\ref{lem:consensusMapInvariance}(ii) requires
\begin{equation*}
    \alpha\;\ge\;\alpha_0\Big(\frac{\epsilon}{4M_1},S_{\epsilon/2},C_{s_f(\epsilon/2)}\Big)
    \;\ge\;\frac{2}{S_{\epsilon/2}}\log\frac{4M_1}{\epsilon\,C_{s_f(\epsilon/2)}},
\end{equation*}
and by Lemma~\ref{lem:normalSuitableForLaplacePrinciple},
\begin{equation*}
    \log\frac{1}{C_s}=\frac{(r+s)^2}{2\sigma_{\min}(\Sigma)^2}-\log\frac{\pi^{d/2}s^d}{\Gamma(\frac d2+1)}+\frac d2\log\big(2\pi\sigma_{\max}(\Sigma)^2\big).
\end{equation*}
Hence this threshold contains the product of $1/S_{\epsilon/2}$, which is large for small $\epsilon$, and $r^2/(2\sigma_{\min}(\Sigma)^2)$, which is large for large $r$. In Corollary~\ref{cor:meanFieldAccuracy}, the radius only enters the threshold $\alpha_0(f,\Sigma,r,\frac12)$. There, $r^2/(2\sigma_{\min}(\Sigma)^2)$ is divided by quantities that do not depend on $\epsilon$. The accuracy only enters the threshold $\alpha_{\mathrm b}(f,\Sigma,\frac\epsilon4)$, which does not depend on $r$. If, for instance, $S_\epsilon=(c\,\epsilon)^q$ for small $\epsilon$ with constants $c,q>0$, then up to logarithmic factors the single-step threshold grows like $r^2\epsilon^{-q}$, whereas the threshold of Corollary~\ref{cor:meanFieldAccuracy} grows like $\max\{r^2,\epsilon^{-q}\}$. The price is $K_0(r,\epsilon)=\mathcal O(\log(r/\epsilon))$ iterations.
\end{remark}

\subsection{Finite-Sample Analysis}
We start with a finite-sample bound for the estimation error of the empirical transition map interpreted as a self-normalised importance-weighted estimator of the mean-field transition map.
For clarity, we state and prove a more general result that makes clear which properties of the Gibbs tilt and the self-normalised importance weighted estimator are actually used.
%
%
\begin{lemma}[One-step concentration of the empirical transition map]\label{lem:empiricalTransitionMapConcentration}
Let $f:\R^d\to\R$ be Borel measurable and bounded below, $\underline f:=\inf_{y\in\R^d}f(y)>-\infty$, let $\alpha>0$, let $\Sigma\in\R^{d\times d}$ be symmetric positive definite, let $N\in\Np$ and $x\in\R^d$, and let $X_1,\dots,X_N$ be independent with $X_i\distr\Norm(x,\Sigma)$.
\begin{enumerate}
    \item[(i)] For all $t>0$,
    \begin{equation}\label{eq:empiricalTransitionMapConcentration:general}
        \Pb\Big(\big\|\widehat{\mathcal{T}}_{\alpha,N}(x)-\mathcal{T}_\alpha(x)\big\|_2>t\Big)
        \;\le\;\frac{4}{N}\,\frac{e^{-\alpha\underline f}}{Z_\alpha(\Norm(x,\Sigma))}\left(1+\frac{V_\alpha(x)}{t^2}\right),
    \end{equation}
    where
    \begin{equation}\label{eq:empiricalTransitionMapConcentration:variance}
        V_\alpha(x):=\E_{Y\distr\eta^\alpha_{\Norm(x,\Sigma)}}\big[\|Y-\mathcal{T}_\alpha(x)\|_2^2\big]<\infty
    \end{equation}
    is the total variance of the Gibbs tilt $\eta^\alpha_{\Norm(x,\Sigma)}$.
    \item[(ii)] Let $f$ in addition satisfy Assumption~\ref{ass:objective}, let $r>0$, let $C_s$ and $M_2$ be the constants of Lemma~\ref{lem:normalSuitableForLaplacePrinciple} for $r$ and $\Sigma$, and let $s_f(\epsilon)$ be as in \eqref{eq:laplacePrinciples:epsf}. Then for all $x\in\ball_r(x_\ast)$, $\epsilon,\epsilon'>0$ and $t>0$,
    \begin{equation}\label{eq:empiricalTransitionMapConcentration:uniform}
        \Pb\Big(\big\|\widehat{\mathcal{T}}_{\alpha,N}(x)-\mathcal{T}_\alpha(x)\big\|_2>t\Big)
        \;\le\;\frac{4}{N}\,\frac{e^{\alpha S_\epsilon/2}}{C_{s_f(\epsilon)}}\left(1+\frac{1}{t^2}\Big(\epsilon'^2+\frac{e^{-\alpha S_{\epsilon'}/2}}{C_{s_f(\epsilon')}}\,M_2\Big)\right).
    \end{equation}
\end{enumerate}
\end{lemma}
\begin{proof}
\emph{Ad (i).}
Write $\nu:=\Norm(x,\Sigma)$, $w(y):=e^{-\alpha f(y)}$, $W:=e^{-\alpha\underline f}$, $Z:=Z_\alpha(\nu)$ and $c:=\mathcal{T}_\alpha(x)$, and let $X\distr\nu$.
Then $0<w\le W$, and by Lemma~\ref{lem:consensusMapDifferentiable}(i) we have $0<Z\le W$.
Moreover, by Lemma~\ref{lem:consensusMapDifferentiable}(i), $c$ is well defined and is the mean of $\eta^\alpha_\nu$, the definition \eqref{eq:meanFieldTransitionMap} gives $\E[w(X)X]=Z\,c$, and since $\eta^\alpha_\nu$ has finite second moments, $V_\alpha(x)\le2\E_{Y\distr\eta^\alpha_\nu}\big[\|Y\|_2^2\big]+2\|c\|_2^2<\infty$.

We first collect two moment bounds. Since $0<w\le W$,
\begin{equation}\label{eq:empiricalTransitionMapConcentration:varw}
    \Var\big(w(X)\big)\;\le\;\E\big[w(X)^2\big]\;\le\;W\,\E\big[w(X)\big]\;=\;W Z .
\end{equation}
Moreover, \eqref{eq:consensusMapDifferentiable:changeofmeasure} with $h(y)=\|y-c\|_2^2$ gives
\begin{equation}\label{eq:empiricalTransitionMapConcentration:weightedsecond}
    \E\big[w(X)^2\|X-c\|_2^2\big]
    \;\le\;W\,\E\big[w(X)\|X-c\|_2^2\big]
    \;=\;W Z\,V_\alpha(x) .
\end{equation}

Now define
\begin{equation*}
    Y_i:=w(X_i)(X_i-c),
    \qquad
    A_N:=\frac1N\sum_{i=1}^N Y_i,
    \qquad
    B_N:=\frac1N\sum_{i=1}^N w(X_i).
\end{equation*}
The $Y_i$ are independent and identically distributed, square-integrable by \eqref{eq:empiricalTransitionMapConcentration:weightedsecond}, and centred, since $\E[Y_i]=\E[w(X)X]-c\,\E[w(X)]=Zc-cZ=0$.
Since $B_N>0$, we have
\begin{equation*}
    \widehat{\mathcal{T}}_{\alpha,N}(x)-c=\frac{A_N}{B_N}.
\end{equation*}
Let $E_N:=\{B_N\ge Z/2\}$. Then
\begin{equation*}
    \Pb\Big(\big\|\widehat{\mathcal{T}}_{\alpha,N}(x)-c\big\|_2>t\Big)
    \;\le\;\Pb\big(E_N^c\big)+\Pb\Big(\big\|\widehat{\mathcal{T}}_{\alpha,N}(x)-c\big\|_2>t,\ E_N\Big).
\end{equation*}
Since $\E[B_N]=Z$ and $\Var(B_N)=\Var(w(X))/N$, Chebyshev's inequality and \eqref{eq:empiricalTransitionMapConcentration:varw} yield
\begin{equation*}
    \Pb\big(E_N^c\big)
    \;\le\;\Pb\big(|B_N-Z|\ge Z/2\big)
    \;\le\;\frac{4\Var(w(X))}{NZ^2}
    \;\le\;\frac{4}{N}\,\frac{W}{Z}.
\end{equation*}
On $E_N$, we have $\|\widehat{\mathcal{T}}_{\alpha,N}(x)-c\|_2=\|A_N\|_2/B_N\le 2\|A_N\|_2/Z$.
Since the $Y_i$ are independent and centred, the cross terms vanish, and with \eqref{eq:empiricalTransitionMapConcentration:weightedsecond} we get
\begin{equation*}
    \E\big[\|A_N\|_2^2\big]
    =\frac1{N^2}\sum_{i,j=1}^N\E\big[\langle Y_i,Y_j\rangle\big]
    =\frac1N\E\big[\|Y_1\|_2^2\big]
    \;\le\;\frac{W Z\,V_\alpha(x)}{N}.
\end{equation*}
Markov's inequality applied to $\|A_N\|_2^2$ therefore gives
\begin{equation*}
    \Pb\Big(\big\|\widehat{\mathcal{T}}_{\alpha,N}(x)-c\big\|_2>t,\ E_N\Big)
    \;\le\;\Pb\big(\|A_N\|_2> tZ/2\big)
    \;\le\;\frac{4\,\E\big[\|A_N\|_2^2\big]}{t^2Z^2}
    \;\le\;\frac{4}{N}\,\frac{W}{Z}\,\frac{V_\alpha(x)}{t^2}.
\end{equation*}
Adding the two bounds proves \eqref{eq:empiricalTransitionMapConcentration:general}.

\emph{Ad (ii).}
By (F1), $\underline f=f_\ast$. Let $x\in\ball_r(x_\ast)$, $\nu:=\Norm(x,\Sigma)$, and $\epsilon,\epsilon'>0$.
By Lemma~\ref{lem:normalSuitableForLaplacePrinciple}, $\nu$ satisfies \eqref{eq:laplacePrinciples:mass} with the constants $C_s$ and \eqref{eq:laplacePrinciples:moment} for $p=2$ with $M_2$.
Restricting the integral to $\ball_{s_f(\epsilon)}(x_\ast)$ and using \eqref{eq:laplacePrinciples:epsf} and \eqref{eq:laplacePrinciples:mass}, we obtain, as in \eqref{eq:laplacePrinciples:normaliser},
\begin{equation*}
    Z_\alpha(\nu)
    \;\ge\;\int_{\ball_{s_f(\epsilon)}(x_\ast)}e^{-\alpha f(y)}\,\nu(\mathrm dy)
    \;\ge\;e^{-\alpha(f_\ast+S_\epsilon/2)}\,\nu\big(\ball_{s_f(\epsilon)}(x_\ast)\big)
    \;\ge\;e^{-\alpha(f_\ast+S_\epsilon/2)}\,C_{s_f(\epsilon)},
\end{equation*}
and hence $e^{-\alpha\underline f}/Z_\alpha(\nu)\le e^{\alpha S_\epsilon/2}/C_{s_f(\epsilon)}$.
Furthermore, since the mean minimises the expected squared distance, Lemma~\ref{lem:laplacePrinciples}(ii) with $p=2$ and $\epsilon'$ in place of $\epsilon$ yields
\begin{equation*}
    V_\alpha(x)
    \;\le\;\E_{Y\distr\eta^\alpha_\nu}\big[\|Y-x_\ast\|_2^2\big]
    \;\le\;\epsilon'^2+\frac{e^{-\alpha S_{\epsilon'}/2}}{C_{s_f(\epsilon')}}\,M_2 .
\end{equation*}
Since the right-hand side of \eqref{eq:empiricalTransitionMapConcentration:general} is nondecreasing in $e^{-\alpha\underline f}/Z_\alpha(\nu)$ and in $V_\alpha(x)$, inserting both bounds proves \eqref{eq:empiricalTransitionMapConcentration:uniform}.
\end{proof}

We can now establish Proposition~\ref{prop:snisErrorBound}.
Note that this result only uses (F1) of Assumption~\ref{ass:objective}, namely that $f$ is measurable and bounded from below with $\inf_{y\in\R^d}f(y)=f(x_\ast)$. 
More generally, by Lemma~\ref{lem:empiricalTransitionMapConcentration}(i), the bound \eqref{eq:snisErrorBound:bound} holds for every measurable $f$ that is bounded from below, if $f(x_\ast)$ in the definition \eqref{eq:snis:quantities} of $\zeta_\alpha(x)$ is replaced by $\underline f=\inf_{y\in\R^d}f(y)$.
\begin{proof}[Proof of Proposition~\ref{prop:snisErrorBound}]
By (F1), $f$ is measurable and bounded from below with $\underline f:=\inf_{y\in\R^d}f(y)=f(x_\ast)$, so Lemmas~\ref{lem:consensusMapDifferentiable} and~\ref{lem:empiricalTransitionMapConcentration} apply.
Let $x\in\R^d$ and $\nu:=\Norm(x,\Sigma)$. Since $0<e^{-\alpha(f(y)-f(x_\ast))}\le1$ for all $y\in\R^d$, we have $\zeta_\alpha(x)\in(0,1]$. Moreover, $Z_\alpha(x)=Z_\alpha(\nu)$, and hence $\zeta_\alpha(x)=e^{\alpha f(x_\ast)}Z_\alpha(\nu)$, i.e.\ $1/\zeta_\alpha(x)=e^{-\alpha\underline f}/Z_\alpha(\nu)$.
Furthermore, \eqref{eq:consensusMapDifferentiable:changeofmeasure} with $h(y)=\|y-\mathcal{T}_\alpha(x)\|_2^2$ shows that $V_\alpha(x)$ in \eqref{eq:snis:quantities} coincides with the total variance \eqref{eq:empiricalTransitionMapConcentration:variance} of the Gibbs tilt $\eta^\alpha_\nu$, which is finite by Lemma~\ref{lem:empiricalTransitionMapConcentration}(i).
Hence \eqref{eq:snisErrorBound:bound} is \eqref{eq:empiricalTransitionMapConcentration:general}.
\end{proof}

We can now turn to sample complexity results.
To illustrate the approach, we first consider only a single step in the MPPI iteration.
%
%
\begin{corollary}[Sample complexity, single step]\label{cor:empiricalTransitionMapSampleComplexity}
Let $f$ satisfy Assumption~\ref{ass:objective}, let $\Sigma\in\R^{d\times d}$ be symmetric positive definite, and let $r>0$.
Let $C_s$, $M_1$ and $M_2$ be the constants of Lemma~\ref{lem:normalSuitableForLaplacePrinciple} for $r$ and $\Sigma$, and let $s_f(\epsilon)$ be as in \eqref{eq:laplacePrinciples:epsf}.
For $\epsilon>0$ and $\delta\in(0,1)$, set
\begin{equation}\label{eq:empiricalTransitionMapSampleComplexity:constants}
    \bar C_\epsilon:=C_{s_f(\epsilon/4)},
    \qquad
    \bar S_\epsilon:=S_{\epsilon/4},
    \qquad
    \Lambda_\epsilon:=\max\Big\{\frac{8M_1}{\epsilon},\,\frac{16M_2}{\epsilon^2}\Big\},
\end{equation}
and
\begin{equation}\label{eq:empiricalTransitionMapSampleComplexity:thresholds}
    \alpha_{\min}(\epsilon):=\frac{2}{\bar S_\epsilon}\max\Big\{0,\;\log\frac{\Lambda_\epsilon}{\bar C_\epsilon}\Big\},
    \qquad
    N_{\min}(\alpha,\epsilon,\delta):=\frac{6\,e^{\alpha \bar S_\epsilon/2}}{\delta\,\bar C_\epsilon}.
\end{equation}
Then for all $\alpha>0$ with $\alpha\ge\alpha_{\min}(\epsilon)$, all $N\in\Np$ with $N\ge N_{\min}(\alpha,\epsilon,\delta)$, and all $x\in\ball_r(x_\ast)$,
\begin{equation}\label{eq:empiricalTransitionMapSampleComplexity:claim}
    \Pb\Big(\big\|\widehat{\mathcal{T}}_{\alpha,N}(x)-x_\ast\big\|_2<\epsilon\Big)\;\ge\;1-\delta .
\end{equation}
In particular, if $\alpha_{\min}(\epsilon)>0$, then for $\alpha=\alpha_{\min}(\epsilon)$ the sample size condition reads
\begin{equation}\label{eq:empiricalTransitionMapSampleComplexity:atAlphaMin}
    N\;\ge\;N_{\min}\big(\alpha_{\min}(\epsilon),\epsilon,\delta\big)=\frac{6}{\delta\,\bar C_\epsilon}\max\Big\{1,\frac{\Lambda_\epsilon}{\bar C_\epsilon}\Big\}.
\end{equation}
\end{corollary}
\begin{proof}
Fix $\alpha$, $N$ and $x$ as in the statement. Recall that $\alpha_0(\delta',S,C)=\frac2S\max\{0,\log\frac{1}{\delta'C}\}$ by \eqref{eq:laplacePrinciples:alpha0}, which is nondecreasing in $1/\delta'$.

\emph{Bias.}
Lemma~\ref{lem:consensusMapInvariance}(ii) with $\rho=\epsilon/2$ applies for $\alpha\ge\alpha_0\big(\frac{\epsilon}{8M_1},\bar S_\epsilon,\bar C_\epsilon\big)=\frac{2}{\bar S_\epsilon}\max\{0,\log\frac{8M_1}{\epsilon \bar C_\epsilon}\}$, and this holds since $8M_1/\epsilon\le\Lambda_\epsilon$ and $\alpha\ge\alpha_{\min}(\epsilon)$. Hence
\begin{equation}\label{eq:empiricalTransitionMapSampleComplexity:bias}
    \big\|\mathcal{T}_\alpha(x)-x_\ast\big\|_2\;\le\;\tfrac34\cdot\tfrac\epsilon2=\tfrac{3\epsilon}{8}.
\end{equation}

\emph{Variance of the Gibbs tilt.}
Similarly, $16M_2/\epsilon^2\le\Lambda_\epsilon$ gives $\alpha\ge\alpha_0\big(\frac{\epsilon^2}{16M_2},\bar S_\epsilon,\bar C_\epsilon\big)$, so \eqref{eq:laplacePrinciples:alpha0bound} yields
\begin{equation*}
    \frac{e^{-\alpha \bar S_\epsilon/2}}{\bar C_\epsilon}\,M_2\;\le\;\frac{\epsilon^2}{16}.
\end{equation*}

\emph{Sampling error.}
Lemma~\ref{lem:empiricalTransitionMapConcentration}(ii), applied with $\epsilon/4$ in place of both $\epsilon$ and $\epsilon'$ and with $t=\epsilon/2$, together with the previous estimate gives
\begin{equation*}
    \Pb\Big(\big\|\widehat{\mathcal{T}}_{\alpha,N}(x)-\mathcal{T}_\alpha(x)\big\|_2>\tfrac\epsilon2\Big)
    \;\le\;\frac4N\,\frac{e^{\alpha \bar S_\epsilon/2}}{\bar C_\epsilon}\left(1+\frac{4}{\epsilon^2}\Big(\frac{\epsilon^2}{16}+\frac{\epsilon^2}{16}\Big)\right)
    =\frac6N\,\frac{e^{\alpha \bar S_\epsilon/2}}{\bar C_\epsilon}
    \;\le\;\delta,
\end{equation*}
where the last step uses $N\ge N_{\min}(\alpha,\epsilon,\delta)$.

\emph{Conclusion.}
On the complement of the event in the last display, \eqref{eq:empiricalTransitionMapSampleComplexity:bias} and the triangle inequality give $\|\widehat{\mathcal{T}}_{\alpha,N}(x)-x_\ast\|_2\le\frac\epsilon2+\frac{3\epsilon}8<\epsilon$, which proves \eqref{eq:empiricalTransitionMapSampleComplexity:claim}.
Finally, if $\alpha=\alpha_{\min}(\epsilon)>0$, then $e^{\alpha_{\min}(\epsilon)\bar S_\epsilon/2}=\Lambda_\epsilon/\bar C_\epsilon=\max\{1,\Lambda_\epsilon/\bar C_\epsilon\}$, which gives \eqref{eq:empiricalTransitionMapSampleComplexity:atAlphaMin}.
\end{proof}

Let us now prove Theorem \ref{thm:finiteSample}, which considers a finite trajectory of the MPPI iteration.
\begin{proof}[Proof of Theorem~\ref{thm:finiteSample}]
Throughout, for $R>0$ we write $C_s(R)$ and $M_2(R)$ for the constants of Lemma~\ref{lem:normalSuitableForLaplacePrinciple} with radius $R$.
Let $r$, $\rho$, $\beta$, $\gamma$ and $\rho'$ be as in the statement, and set
\begin{equation*}
    \theta:=\frac{1-\rho}{2},
    \qquad
    R_0:=\sigma_{\min}(\Sigma),
\end{equation*}
so that $\rho'=\rho+\theta$ and $\theta+\rho'=1$. By (F2), there is $s_\gamma>0$, depending only on $f$ and $\gamma$, such that $f\le f_\ast+\gamma$ on $\ball_{s_\gamma}(x_\ast)$.
With $\alpha_0(\cdot,\cdot,\cdot)$ from \eqref{eq:laplacePrinciples:alpha0}, $\alpha_0(f,\Sigma,r,\rho)$ from \eqref{eq:meanFieldConvergence:alpha0} and $\alpha_{\mathrm b}$ from \eqref{eq:biasAtMinimiser:alphab}, define
\begin{equation}\label{eq:finiteSample:thresholds}
\begin{aligned}
    \alpha_0'(f,\Sigma,r,\rho)&:=\max\Big\{\alpha_0(f,\Sigma,r,\rho),\ \alpha_0\Big(\frac{\theta^2R_0^2}{4M_2(r)},\,S_{\theta R_0/2},\,C_{s_f(\theta R_0/2)}(r)\Big)\Big\},\\
    \alpha_{\mathrm b}'(f,\Sigma,\beta)&:=\max\Big\{\alpha_{\mathrm b}(f,\Sigma,\beta),\ \alpha_0\Big(\frac{\beta^2}{4M_2(R_0)},\,S_{\beta/2},\,C_{s_f(\beta/2)}(R_0)\Big)\Big\},\\
    C(f,\Sigma,r,\gamma)&:=\frac{6}{C_{s_\gamma}(r)} .
\end{aligned}
\end{equation}
Then $\alpha_0'$ depends only on $f$, $\Sigma$, $r$ and $\rho$, $\alpha_{\mathrm b}'$ only on $f$, $\Sigma$ and $\beta$, and $C$ only on $f$, $\Sigma$, $r$ and $\gamma$. Fix $\alpha$, $K$, $\delta$, $N$ and $\hat x_0$ as in the statement.

\emph{Step 0: Realisation of the iteration.}
Let $(G^{(k)}_i)_{k\in\Nnn,\,i=1,\dots,N}$ be independent $\Norm(0,\Ident_d)$-distributed random vectors, let $\Sigma^{1/2}$ be the symmetric positive definite square root of $\Sigma$, and write $G^{(k)}:=(G^{(k)}_1,\dots,G^{(k)}_N)$. For $x\in\R^d$ and $g\in(\R^d)^N$, let $\widehat T(x,g)$ be the right-hand side of \eqref{eq:empiricalTransitionMap} with $X_i=x+\Sigma^{1/2}g_i$. Since $f$ is measurable and the denominator is strictly positive, $\widehat T$ is Borel measurable, and the MPPI iterates are given by $\hat x_{k+1}=\widehat T(\hat x_k,G^{(k)})$. By induction, $\hat x_k$ is a Borel measurable function of $G^{(0)},\dots,G^{(k-1)}$, and hence independent of $G^{(k)}$. We write $e_k:=\|\hat x_k-x_\ast\|_2$.

\emph{Step 1: Normaliser and variance.}
Let $x\in\ball_r(x_\ast)$, $e:=\|x-x_\ast\|_2$ and $\nu:=\Norm(x,\Sigma)$. By (F1), $\underline f=f_\ast$. Restricting the integral to $\ball_{s_\gamma}(x_\ast)$ and using Lemma~\ref{lem:normalSuitableForLaplacePrinciple} with radius $r$ gives $Z_\alpha(x)\ge e^{-\alpha(f_\ast+\gamma)}C_{s_\gamma}(r)$, i.e.\
\begin{equation}\label{eq:finiteSample:zeta}
    \frac{1}{\zeta_\alpha(x)}\;\le\;\frac{e^{\alpha\gamma}}{C_{s_\gamma}(r)} .
\end{equation}
Next, we show
\begin{equation}\label{eq:finiteSample:variance}
    V_\alpha(x)\;\le\;\tfrac12\,(\theta e+\beta)^2 .
\end{equation}
By Proposition~\ref{prop:snisErrorBound} and its proof, $V_\alpha(x)$ is the total variance of $\eta^\alpha_\nu$, and since the mean minimises the expected squared distance, $V_\alpha(x)\le\E_{Y\distr\eta^\alpha_\nu}\big[\|Y-x_\ast\|_2^2\big]$. We bound the right-hand side with Lemma~\ref{lem:laplacePrinciples}(ii) for $p=2$ in two regimes.
If $e\le R_0$, then $x\in\ball_{R_0}(x_\ast)$, so $\nu$ satisfies \eqref{eq:laplacePrinciples:mass} and \eqref{eq:laplacePrinciples:moment} for $p=2$ with $C_s(R_0)$ and $M_2(R_0)$. With $\epsilon=\beta/2$, \eqref{eq:laplacePrinciples:momentbound}, $\alpha\ge\alpha_{\mathrm b}'(f,\Sigma,\beta)$ and \eqref{eq:laplacePrinciples:alpha0bound} give $\E_{Y\distr\eta^\alpha_\nu}[\|Y-x_\ast\|_2^2]\le\frac{\beta^2}4+\frac{\beta^2}{4}\le\frac12(\theta e+\beta)^2$.
If $R_0<e\le r$, we use the constants $C_s(r)$ and $M_2(r)$. With $\epsilon=\theta R_0/2$ and $\alpha\ge\alpha_0'(f,\Sigma,r,\rho)$, we obtain $\E_{Y\distr\eta^\alpha_\nu}[\|Y-x_\ast\|_2^2]\le\frac{\theta^2R_0^2}4+\frac{\theta^2R_0^2}4<\frac12\theta^2e^2\le\frac12(\theta e+\beta)^2$.

\emph{Step 2: One sampled step.}
For $x\in\R^d$, let
\begin{equation*}
    \psi(x):=\Pb\Big(\big\|\widehat T(x,G^{(0)})-\mathcal{T}_\alpha(x)\big\|_2>\theta\|x-x_\ast\|_2+\beta\Big).
\end{equation*}
Let $x\in\ball_r(x_\ast)$. We apply Proposition~\ref{prop:snisErrorBound} with $t=\theta\|x-x_\ast\|_2+\beta$ and use \eqref{eq:finiteSample:zeta}, \eqref{eq:finiteSample:variance} and \eqref{eq:finiteSample:N}. This yields
\begin{equation}\label{eq:finiteSample:psi}
    \psi(x)
    \;\le\;\frac4N\,\frac{1}{\zeta_\alpha(x)}\Big(1+\frac{V_\alpha(x)}{(\theta\|x-x_\ast\|_2+\beta)^2}\Big)
    \;\le\;\frac4N\,\frac{e^{\alpha\gamma}}{C_{s_\gamma}(r)}\cdot\frac32
    \;=\;\frac{C(f,\Sigma,r,\gamma)\,e^{\gamma\alpha}}{N}
    \;\le\;\frac\delta K .
\end{equation}

\emph{Step 3: Deterministic progress on good events.}
For $k\in\Nnn$, let $A_k:=\big\{\|\hat x_{k+1}-\mathcal{T}_\alpha(\hat x_k)\|_2\le\theta e_k+\beta\big\}$. Suppose that $\hat x_k\in\ball_r(x_\ast)$ and that $A_k$ occurs. Since $\alpha\ge\alpha_0'\ge\alpha_0(f,\Sigma,r,\rho)$ and $\alpha\ge\alpha_{\mathrm b}'\ge\alpha_{\mathrm b}(f,\Sigma,\beta)$, Theorem~\ref{thm:meanFieldOneStep}(i) gives $\|\mathcal{T}_\alpha(\hat x_k)-x_\ast\|_2\le\beta+\rho e_k$, and hence
\begin{equation}\label{eq:finiteSample:recursion}
    e_{k+1}\;\le\;\|\mathcal{T}_\alpha(\hat x_k)-x_\ast\|_2+\theta e_k+\beta\;\le\;2\beta+\rho' e_k .
\end{equation}
In particular, $e_{k+1}\le2\beta+\rho'r\le\theta r+\rho'r=r$, where we used $2\beta\le\theta r$, i.e.\ $4\beta\le(1-\rho)r$. By induction, on $\bigcap_{j=0}^{k-1}A_j$ we have $\hat x_j\in\ball_r(x_\ast)$ for $j=0,\dots,k$, and \eqref{eq:finiteSample:recursion} yields
\begin{equation*}
    e_k\;\le\;\rho'^ke_0+2\beta\sum_{j=0}^{k-1}\rho'^j\;\le\;\rho'^ke_0+\frac{2\beta}{1-\rho'}\;=\;\rho'^ke_0+\frac{4\beta}{1-\rho} .
\end{equation*}
Hence the event in \eqref{eq:finiteSample:claim} contains $\bigcap_{k=0}^{K-1}A_k$.

\emph{Step 4: Probability of the good events.}
Since $\hat x_k$ is independent of $G^{(k)}$, $\widehat T$ is Borel measurable and $\mathcal{T}_\alpha$ is continuous by Lemma~\ref{lem:consensusMapDifferentiable}, Fubini's theorem and \eqref{eq:finiteSample:psi} give $\Pb\big(A_k^c\cap\{\hat x_k\in\ball_r(x_\ast)\}\big)=\E\big[\mathbf 1_{\ball_r(x_\ast)}(\hat x_k)\,\psi(\hat x_k)\big]\le\delta/K$. Decomposing according to the first failure time and using Step~3,
\begin{equation*}
    \Pb\Big(\bigcup_{k=0}^{K-1}A_k^c\Big)
    =\sum_{k=0}^{K-1}\Pb\Big(A_k^c\cap\bigcap_{j=0}^{k-1}A_j\Big)
    \;\le\;\sum_{k=0}^{K-1}\Pb\Big(A_k^c\cap\big\{\hat x_k\in\ball_r(x_\ast)\big\}\Big)
    \;\le\;\delta ,
\end{equation*}
which proves \eqref{eq:finiteSample:claim}.

\emph{Step 5: Accuracy $\epsilon$.}
Let $\epsilon\in(0,2r]$, and apply \eqref{eq:finiteSample:claim} with $\rho=\frac12$ and $\beta=\frac\epsilon{16}$. Then $4\beta=\frac\epsilon4\le\frac r2=(1-\rho)r$, $\rho'=\frac34$ and $\frac{4\beta}{1-\rho}=\frac\epsilon2$. For $k\ge K_\epsilon$ we have $(\frac34)^k\le\frac{\epsilon}{2r}$, and hence $\rho'^k\|\hat x_0-x_\ast\|_2+\frac\epsilon2\le\epsilon$ on the event in \eqref{eq:finiteSample:claim}. This proves \eqref{eq:finiteSample:accuracy}.
\end{proof}
\subsection{Analysis under a Simplified Assumption}\label{app:simplified}
In this section, we prove the results stated under Assumption~\ref{ass:simplified}. Under this assumption, the covariance of the Gibbs tilt inside the well can be controlled by a log-concavity argument, and all thresholds become explicit. Throughout, $\Sigma\in\R^{d\times d}$ is symmetric positive definite, and $\kappa_\Sigma=\sigma_{\max}(\Sigma)^2/\sigma_{\min}(\Sigma)^2$ as in the notation paragraph.

\subsubsection{Assumption \ref{ass:simplified} implies Assumption \ref{ass:objective}}
For $s>0$, we use the abbreviation
\begin{equation}\label{eq:simplified:Ds}
    D_s:=\Gamma\!\left(\frac d2+1\right)\Big(\frac{2\sigma_{\max}(\Sigma)^2}{s^2}\Big)^{d/2},
\end{equation}
so that the constants of Lemma~\ref{lem:normalSuitableForLaplacePrinciple} and \eqref{eq:biasAtMinimiser:constants} read
\begin{equation}\label{eq:simplified:Cs}
    \frac{1}{C_s}=D_s\exp\Big(\frac{(r+s)^2}{2\sigma_{\min}(\Sigma)^2}\Big),
    \qquad
    \frac{1}{C^\circ_s}=D_s\exp\Big(\frac{s^2}{2\sigma_{\min}(\Sigma)^2}\Big).
\end{equation}
\begin{lemma}[Consequences of the simplified assumption]\label{lem:simplifiedImpliesObjective}
Let $f$ satisfy (SF1) and (SF2) of Assumption~\ref{ass:simplified}. Then the following hold.
\begin{enumerate}
    \item[(i)] \emph{Quadratic growth in the well.} $f(y)-f_\ast\ge\frac m2\|y-x_\ast\|_2^2$ for all $y\in\ball_R(x_\ast)$.
    \item[(ii)] \emph{Radial monotonicity.} The function $g(y):=f(y)-f_\ast-\frac m2\|y-x_\ast\|_2^2$ is convex and nonnegative on $\ball_R(x_\ast)$ with $g(x_\ast)=0$, and for every unit vector $u$, the map $t\mapsto g(x_\ast+tu)$ is nondecreasing on $[0,R]$.
    \item[(iii)] \emph{Assumption~\ref{ass:objective}.} $f$ satisfies Assumption~\ref{ass:objective}, $x_\ast$ is its unique global minimiser, and \eqref{eq:objective:separation} holds with
    \begin{equation}\label{eq:simplified:S}
        S_\epsilon:=\min\Big\{\frac m2\epsilon^2,\,\Delta\Big\}\qquad\text{for all }\epsilon>0 .
    \end{equation}
    \item[(iv)] \emph{Small ball below half the gap.} There is $s\in(0,R]$ with $\sup_{y\in\ball_s(x_\ast)}\big(f(y)-f_\ast\big)\le\frac\Delta2$. If $f$ also satisfies (SF3), then $s:=\min\{R,\sqrt{\Delta/M}\}$ is admissible.
\end{enumerate}
\end{lemma}
\begin{proof}
Without loss of generality, $x_\ast=0$ and $f_\ast=0$.

\emph{Ad (i).} Let $y\in\ball_R(0)$ and $u\in(0,1)$. The point $uy$ lies in $\ball_R(0)$, and convexity of $g$ from (SF1) on the segment from $0$ to $y$ gives $g(uy)\le(1-u)g(0)+u\,g(y)$, i.e.
\begin{equation*}
    f(uy)-\frac m2u^2\|y\|_2^2\;\le\;u\Big(f(y)-\frac m2\|y\|_2^2\Big).
\end{equation*}
Since $f(uy)\ge0$, this yields $0\le u f(y)-\frac m2u(1-u)\|y\|_2^2$. Dividing by $u$ and letting $u\downarrow0$ gives $f(y)\ge\frac m2\|y\|_2^2$.

\emph{Ad (ii).} Convexity is (SF1), nonnegativity is (i), and $g(0)=f(0)=0$. Let $u$ be a unit vector and $0\le t_1<t_2\le R$. Writing $t_1=(1-\lambda)\cdot0+\lambda t_2$ with $\lambda=t_1/t_2\in[0,1)$, convexity gives $g(t_1u)\le(1-\lambda)g(0)+\lambda g(t_2u)=\lambda g(t_2u)\le g(t_2u)$, where we used $g(0)=0$ and $g(t_2u)\ge0$.

\emph{Ad (iii).} $f$ is continuous, hence Borel measurable and upper semicontinuous at $0$, and $f\ge0=f(0)$ by assumption, so (F1) and (F2) hold. Let $\epsilon>0$ and $\|y\|_2\ge\epsilon$. If $y\in\ball_R(0)$, then (i) gives $f(y)\ge\frac m2\|y\|_2^2\ge\frac m2\epsilon^2$; otherwise, (SF2) gives $f(y)\ge\Delta$. Hence $f(y)\ge S_\epsilon>0$, which is (F3) with \eqref{eq:simplified:S}, and uniqueness of the global minimiser follows.

\emph{Ad (iv).} Since $f$ is continuous with $f(0)=0$, there is $s\in(0,R]$ with $f\le\frac\Delta2$ on $\ball_s(0)$. Under (SF3), $f(y)\le\frac M2\|y\|_2^2\le\frac M2s^2\le\frac\Delta2$ for $\|y\|_2\le s=\min\{R,\sqrt{\Delta/M}\}$.
\end{proof}
In the remainder of this subsection, $s\in(0,R]$ is fixed as in Lemma~\ref{lem:simplifiedImpliesObjective}(iv). For $\alpha>0$, we abbreviate
\begin{equation}\label{eq:simplified:tau}
    \tau_\alpha:=\frac{1}{1+\alpha m\sigma_{\max}(\Sigma)^2},
\end{equation}
and for $x\in\R^d$, we write $P_x:=\eta^\alpha_{\Norm(x,\Sigma)}$ for the Gibbs tilt, so that $\mathcal{T}_\alpha(x)$ is the mean of $P_x$.

\subsubsection{Mean-field Analysis under the Simplified Assumption}
The simplified Assumption \ref{ass:simplified} allows us to make the mean-field analysis more concrete and quantitative.
We start with investigating how the various components of the Gibbs tilt behave.
\begin{lemma}[The Gibbs tilt inside and outside the well]\label{lem:simplifiedTilt}
Let $f$ satisfy (SF1) and (SF2) of Assumption~\ref{ass:simplified}, let $\alpha>0$ and $x\in\R^d$, and let $e:=\|x-x_\ast\|_2$.
\begin{enumerate}
    \item[(i)] \emph{Outside the well.} On $\R^d\setminus\ball_R(x_\ast)$, the density of $P_x$ with respect to $\Norm(x,\Sigma)$ is bounded by
    \begin{equation}\label{eq:simplifiedTilt:ratio}
        q_\alpha(e):=D_s\exp\Big(\frac{(e+s)^2}{2\sigma_{\min}(\Sigma)^2}-\frac{\alpha\Delta}{2}\Big).
    \end{equation}
    \item[(ii)] \emph{Inside the well.} $P_x(\ball_R(x_\ast))>0$, and the conditional law $P^{\mathrm{in}}_x:=P_x(\,\cdot\mid\ball_R(x_\ast))$ satisfies
    \begin{equation}\label{eq:simplifiedTilt:cov}
        \mathrm{Cov}\big(P^{\mathrm{in}}_x\big)\;\preceq\;\big(\Sigma^{-1}+\alpha m\Ident_d\big)^{-1}\;\preceq\;\sigma_{\max}(\Sigma)^2\,\tau_\alpha\,\Ident_d .
    \end{equation}
    \item[(iii)] \emph{Inside the well, centred at the minimiser.} For $x=x_\ast$,
    \begin{equation}\label{eq:simplifiedTilt:secondmoment}
        \E_{Y\distr P^{\mathrm{in}}_{x_\ast}}\big[\|Y-x_\ast\|_2^2\big]\;\le\;\frac{d}{\lambda_{\min}\big(\Sigma^{-1}+\alpha m\Ident_d\big)}\;=\;d\,\sigma_{\max}(\Sigma)^2\,\tau_\alpha .
    \end{equation}
\end{enumerate}
\end{lemma}
\begin{proof}
Without loss of generality, $x_\ast=0$ and $f_\ast=0$. Write $\nu:=\Norm(x,\Sigma)$, $P:=P_x$, $B:=\ball_R(0)$, and $A_\alpha:=\Sigma^{-1}+\alpha m\Ident_d$.

\emph{Ad (i).} By Lemma~\ref{lem:simplifiedImpliesObjective}(iv), $e^{-\alpha f}\ge e^{-\alpha\Delta/2}$ on $\ball_s(0)$, and Lemma~\ref{lem:normalSuitableForLaplacePrinciple} with radius $e$ (for $e=0$ by the limiting argument in the proof of Lemma~\ref{lem:biasAtMinimiser}) together with \eqref{eq:simplified:Cs} gives
\begin{equation*}
    Z_\alpha(\nu)\;\ge\;e^{-\alpha\Delta/2}\,\nu\big(\ball_s(0)\big)\;\ge\;e^{-\alpha\Delta/2}\,D_s^{-1}\exp\Big(-\frac{(e+s)^2}{2\sigma_{\min}(\Sigma)^2}\Big).
\end{equation*}
On $\R^d\setminus B$, (SF2) gives $e^{-\alpha f}\le e^{-\alpha\Delta}$, so the density $e^{-\alpha f}/Z_\alpha(\nu)$ of $P$ with respect to $\nu$ is at most $q_\alpha(e)$ there.

\emph{Ad (ii).} Since $e^{-\alpha f}>0$ and the Gaussian density is positive, $P(B)>0$. Let $g$ be as in Lemma~\ref{lem:simplifiedImpliesObjective}(ii), and let $w(y):=\mathbf 1_B(y)\,e^{-\alpha g(y)}$, which is log-concave on $\R^d$, as $g$ is convex on the convex set $B$. With $Z_{\mathrm{in}}(x):=\int_Be^{-\alpha f(y)}\Norm(y\mid x,\Sigma)\,\mathrm dy>0$, completing the square gives
\begin{equation*}
    \alpha\frac m2\|y\|_2^2+\frac12(y-x)^\top\Sigma^{-1}(y-x)
    =\frac12\big(y-A_\alpha^{-1}\Sigma^{-1}x\big)^\top A_\alpha\big(y-A_\alpha^{-1}\Sigma^{-1}x\big)+\frac12x^\top Ex,
\end{equation*}
with the symmetric positive semidefinite matrix $E:=\Sigma^{-1}-\Sigma^{-1}A_\alpha^{-1}\Sigma^{-1}$. Hence
\begin{equation}\label{eq:simplifiedTilt:Zin}
    \log Z_{\mathrm{in}}(x)+\frac12x^\top Ex
    =c+\log\int_{\R^d}w(y)\exp\Big(-\frac12\big(y-A_\alpha^{-1}\Sigma^{-1}x\big)^\top A_\alpha\big(y-A_\alpha^{-1}\Sigma^{-1}x\big)\Big)\,\mathrm dy
\end{equation}
with a constant $c$ that does not depend on $x$. The integrand is a log-concave function of $(x,y)\in\R^d\times\R^d$, as the logarithm of the Gaussian factor is a concave quadratic function of $(x,y)$. By Pr\'ekopa's theorem \cite{prekopa1973logarithmic}, the right-hand side of \eqref{eq:simplifiedTilt:Zin} is a concave function of $x$.

On the other hand, the proof of Lemma~\ref{lem:consensusMapDifferentiable} applies verbatim with all integrands multiplied by $\mathbf 1_B$, and shows that $x\mapsto Z_{\mathrm{in}}(x)$ and the mean $\mathcal{T}^{\mathrm{in}}(x)$ of $P^{\mathrm{in}}_x$ are continuously differentiable, with $\nabla\log Z_{\mathrm{in}}(x)=\Sigma^{-1}\big(\mathcal{T}^{\mathrm{in}}(x)-x\big)$ and $\nabla\mathcal{T}^{\mathrm{in}}(x)=\mathrm{Cov}(P^{\mathrm{in}}_x)\Sigma^{-1}$. Hence $\log Z_{\mathrm{in}}$ is twice continuously differentiable with
\begin{equation*}
    \nabla^2\log Z_{\mathrm{in}}(x)=\Sigma^{-1}\mathrm{Cov}(P^{\mathrm{in}}_x)\Sigma^{-1}-\Sigma^{-1}.
\end{equation*}
Concavity of the left-hand side of \eqref{eq:simplifiedTilt:Zin} gives $\nabla^2\log Z_{\mathrm{in}}(x)+E\preceq0$, i.e.
\begin{equation*}
    \Sigma^{-1}\mathrm{Cov}(P^{\mathrm{in}}_x)\Sigma^{-1}\;\preceq\;\Sigma^{-1}-E\;=\;\Sigma^{-1}A_\alpha^{-1}\Sigma^{-1}.
\end{equation*}
Multiplying from both sides by $\Sigma$ gives the first inequality in \eqref{eq:simplifiedTilt:cov}. The second follows from $A_\alpha\succeq(\sigma_{\max}(\Sigma)^{-2}+\alpha m)\Ident_d$ and $(\sigma_{\max}(\Sigma)^{-2}+\alpha m)^{-1}=\sigma_{\max}(\Sigma)^2\tau_\alpha$.

\emph{Ad (iii).} For $x=0$, completing the square shows that $P^{\mathrm{in}}_0$ has the Lebesgue density proportional to $w(y)\exp(-\frac12y^\top A_\alpha y)$. We use polar coordinates $y=tu$ with $t>0$ and $u\in\mathbb S^{d-1}$, and write $a(u):=u^\top A_\alpha u\ge\lambda_{\min}(A_\alpha)$. For fixed $u$, let $\pi_u$ be the probability measure on $(0,\infty)$ with density proportional to $t^{d-1}e^{-a(u)t^2/2}$, let $c_u:=\int_0^\infty t^{d-1}e^{-a(u)t^2/2}\,\mathrm dt$, and let $W_u(t):=w(tu)$. By Lemma~\ref{lem:simplifiedImpliesObjective}(ii), $W_u$ is nonincreasing on $(0,\infty)$, as it equals $e^{-\alpha g(tu)}$ for $t\le R$ and vanishes for $t>R$. Since $t\mapsto t^2$ is increasing, Chebyshev's integral inequality gives
\begin{equation*}
    \int_0^\infty t^2W_u(t)\,\pi_u(\mathrm dt)\;\le\;\int_0^\infty t^2\,\pi_u(\mathrm dt)\int_0^\infty W_u(t)\,\pi_u(\mathrm dt)
    =\frac{d}{a(u)}\int_0^\infty W_u(t)\,\pi_u(\mathrm dt),
\end{equation*}
where we used $\int_0^\infty t^2\,\pi_u(\mathrm dt)=d/a(u)$. Integrating over $u$ with respect to the surface measure, with weights $c_u$,
\begin{equation*}
    \E_{Y\distr P^{\mathrm{in}}_0}\big[\|Y\|_2^2\big]
    =\frac{\int_{\mathbb S^{d-1}}c_u\int_0^\infty t^2W_u\,\mathrm d\pi_u\,\mathrm du}{\int_{\mathbb S^{d-1}}c_u\int_0^\infty W_u\,\mathrm d\pi_u\,\mathrm du}
    \;\le\;\sup_{u\in\mathbb S^{d-1}}\frac{d}{a(u)}
    =\frac{d}{\lambda_{\min}(A_\alpha)} ,
\end{equation*}
and $\lambda_{\min}(A_\alpha)\ge\sigma_{\max}(\Sigma)^{-2}+\alpha m$ gives \eqref{eq:simplifiedTilt:secondmoment}.
\end{proof}
Recall that to prove contractivity of the transition map on the mean-field, we followed \cite{fornasier2026consensus} and bounded the gradient.
Under the Assumption \ref{ass:simplified}, this simplifies again.
\begin{lemma}[Explicit Jacobian and bias bounds]\label{lem:simplifiedBounds}
Let $f$ satisfy (SF1) and (SF2) of Assumption~\ref{ass:simplified}, and let $\alpha>0$.
\begin{enumerate}
    \item[(i)] For all $r>0$ and $x\in\ball_r(x_\ast)$,
    \begin{equation}\label{eq:simplifiedBounds:Q}
        \big\|\nabla\mathcal{T}_\alpha(x)\big\|\;\le\;Q_\alpha(r):=\kappa_\Sigma\,\tau_\alpha+\kappa_\Sigma\Big(1+\frac{(r+R)^2}{\sigma_{\max}(\Sigma)^2}\Big)D_s\exp\Big(\frac{(r+s)^2}{2\sigma_{\min}(\Sigma)^2}-\frac{\alpha\Delta}2\Big).
    \end{equation}
    In particular, $\mathcal{T}_\alpha$ is Lipschitz continuous on $\ball_r(x_\ast)$ with constant $Q_\alpha(r)$.
    \item[(ii)] At the minimiser,
    \begin{equation}\label{eq:simplifiedBounds:bias}
        \big\|\mathcal{T}_\alpha(x_\ast)-x_\ast\big\|_2^2\;\le\;d\,\sigma_{\max}(\Sigma)^2\Big(\tau_\alpha+D_s\exp\Big(\frac{s^2}{2\sigma_{\min}(\Sigma)^2}-\frac{\alpha\Delta}{2}\Big)\Big).
    \end{equation}
\end{enumerate}
\end{lemma}
\begin{proof}
Without loss of generality, $x_\ast=0$ and $f_\ast=0$, and write $B:=\ball_R(0)$.

\emph{Ad (i).} Let $x\in\ball_r(0)$, $P:=P_x$, $\nu:=\Norm(x,\Sigma)$, and let $h$ be the mean of $P^{\mathrm{in}}_x$. Since $P^{\mathrm{in}}_x$ is concentrated on the closed convex set $B$, $h\in B$, and hence $\|x-h\|_2\le r+R$. Let $v$ be a unit vector. Since the mean minimises the expected squared distance,
\begin{equation*}
    v^\top\mathrm{Cov}(P)v
    \;\le\;\E_{Y\distr P}\big[(v^\top(Y-h))^2\big]
    =P(B)\,v^\top\mathrm{Cov}(P^{\mathrm{in}}_x)v+\int_{\R^d\setminus B}\big(v^\top(y-h)\big)^2\,P(\mathrm dy).
\end{equation*}
By Lemma~\ref{lem:simplifiedTilt}(ii), the first term is at most $\sigma_{\max}(\Sigma)^2\tau_\alpha$. By Lemma~\ref{lem:simplifiedTilt}(i) with $e\le r$, the second term is at most
\begin{equation*}
    q_\alpha(r)\,\E_{X\distr\nu}\big[(v^\top(X-h))^2\big]
    =q_\alpha(r)\big(v^\top\Sigma v+(v^\top(x-h))^2\big)
    \;\le\;q_\alpha(r)\big(\sigma_{\max}(\Sigma)^2+(r+R)^2\big).
\end{equation*}
Hence $\|\mathrm{Cov}(P)\|\le\sigma_{\max}(\Sigma)^2\tau_\alpha+q_\alpha(r)(\sigma_{\max}(\Sigma)^2+(r+R)^2)$. By Lemma~\ref{lem:consensusMapDifferentiable}, $\|\nabla\mathcal{T}_\alpha(x)\|\le\|\mathrm{Cov}(P)\|\,\|\Sigma^{-1}\|$, and $\|\Sigma^{-1}\|=\sigma_{\min}(\Sigma)^{-2}$ gives \eqref{eq:simplifiedBounds:Q}. The Lipschitz bound follows as in the proof of Corollary~\ref{cor:consensusMapContraction}, since $\ball_r(0)$ is convex.

\emph{Ad (ii).} Let $P:=P_0$. By Jensen's inequality,
\begin{equation*}
    \|\mathcal{T}_\alpha(0)\|_2^2\;\le\;\E_{Y\distr P}\big[\|Y\|_2^2\big]
    =P(B)\,\E_{Y\distr P^{\mathrm{in}}_0}\big[\|Y\|_2^2\big]+\int_{\R^d\setminus B}\|y\|_2^2\,P(\mathrm dy).
\end{equation*}
By Lemma~\ref{lem:simplifiedTilt}(iii), the first term is at most $d\sigma_{\max}(\Sigma)^2\tau_\alpha$. By Lemma~\ref{lem:simplifiedTilt}(i) with $e=0$, the second term is at most $q_\alpha(0)\,\E_{X\distr\Norm(0,\Sigma)}[\|X\|_2^2]=q_\alpha(0)\trace(\Sigma)\le q_\alpha(0)\,d\,\sigma_{\max}(\Sigma)^2$.
\end{proof}
We are now ready to investigate the convergence behaviour in the mean-field.
First, let us define explicit thresholds. 
For $r>0$, $\rho\in(0,1)$ and $\beta>0$, let
\begin{equation}\label{eq:simplified:thresholds}
\begin{aligned}
    \alpha^{\mathrm S}_{\mathrm c}(r,\rho)&:=\max\Big\{\frac{2}{\rho\,m\,\sigma_{\min}(\Sigma)^2},\ \frac2\Delta\Big[\frac{(r+s)^2}{2\sigma_{\min}(\Sigma)^2}+\log\Big(\frac{2\kappa_\Sigma D_s}{\rho}\Big(1+\frac{(r+R)^2}{\sigma_{\max}(\Sigma)^2}\Big)\Big)\Big]\Big\},\\
    \alpha^{\mathrm S}_{\mathrm b}(\beta)&:=\max\Big\{\frac{2d}{m\,\beta^2},\ \frac2\Delta\Big[\frac{s^2}{2\sigma_{\min}(\Sigma)^2}+\log\frac{2d\,\sigma_{\max}(\Sigma)^2D_s}{\beta^2}\Big]\Big\},\\
    \alpha^{\mathrm S}_0(r,\rho)&:=\max\big\{\alpha^{\mathrm S}_{\mathrm c}(r,\rho),\ \alpha^{\mathrm S}_{\mathrm b}\big((1-\rho)r\big)\big\}.
\end{aligned}
\end{equation}
Note that $\alpha^{\mathrm S}_{\mathrm c}(r,\rho)$ and $\alpha^{\mathrm S}_0(r,\rho)$ depend on $r$ and $\rho$ but not on $\beta$, and that $\alpha^{\mathrm S}_{\mathrm b}(\beta)$ does not depend on $r$.

The next result corresponds to Theorems \ref{thm:meanFieldConvergence} and \ref{thm:meanFieldOneStep} under Assumption \ref{ass:simplified}.
\begin{proposition}[Explicit mean-field convergence]\label{prop:simplifiedMeanField}
Let $f$ satisfy (SF1) and (SF2) of Assumption~\ref{ass:simplified}, and let $r>0$ and $\rho\in(0,1)$.
\begin{enumerate}
    \item[(i)] \emph{Contraction.} For all $\alpha\ge\alpha^{\mathrm S}_{\mathrm c}(r,\rho)$, $\mathcal{T}_\alpha$ is Lipschitz continuous on $\ball_r(x_\ast)$ with constant $\rho$.
    \item[(ii)] \emph{Bias at the minimiser.} For all $\beta>0$ and $\alpha\ge\alpha^{\mathrm S}_{\mathrm b}(\beta)$, $\|\mathcal{T}_\alpha(x_\ast)-x_\ast\|_2\le\beta$.
    \item[(iii)] \emph{Invariance and convergence.} For all $\alpha\ge\alpha^{\mathrm S}_0(r,\rho)$, we have $\mathcal{T}_\alpha(\ball_r(x_\ast))\subseteq\ball_r(x_\ast)$, $\mathcal{T}_\alpha$ has exactly one fixed point $m_\alpha$ in $\ball_r(x_\ast)$, and for every $x_0\in\ball_r(x_\ast)$, the mean-field iterates $x_{k+1}=\mathcal{T}_\alpha(x_k)$ satisfy $x_k\in\ball_r(x_\ast)$ and $\|x_k-m_\alpha\|_2\le\rho^k\|x_0-m_\alpha\|_2$ for all $k\in\Nnn$.
    \item[(iv)] \emph{One-step improvement.} For all $\beta>0$ and $\alpha\ge\max\{\alpha^{\mathrm S}_0(r,\rho),\alpha^{\mathrm S}_{\mathrm b}(\beta)\}$, and all $x\in\ball_r(x_\ast)$,
    \begin{equation}\label{eq:simplifiedMeanField:affine}
        \big\|\mathcal{T}_\alpha(x)-x_\ast\big\|_2\;\le\;\beta+\rho\,\|x-x_\ast\|_2 .
    \end{equation}
    Consequently, the iterates in (iii) satisfy $\|x_k-x_\ast\|_2\le\rho^k\|x_0-x_\ast\|_2+\frac{\beta}{1-\rho}$ for all $k\in\Nnn$, and $\|m_\alpha-x_\ast\|_2\le\frac{\beta}{1-\rho}$.
\end{enumerate}
\end{proposition}
\begin{proof}
\emph{Ad (i).} Let $\alpha\ge\alpha^{\mathrm S}_{\mathrm c}(r,\rho)$. Since $\tau_\alpha\le(\alpha m\sigma_{\max}(\Sigma)^2)^{-1}$ and $\kappa_\Sigma/\sigma_{\max}(\Sigma)^2=\sigma_{\min}(\Sigma)^{-2}$, the first threshold gives $\kappa_\Sigma\tau_\alpha\le(\alpha m\sigma_{\min}(\Sigma)^2)^{-1}\le\rho/2$. Taking logarithms, the second threshold is equivalent to
\begin{equation*}
    \kappa_\Sigma\Big(1+\frac{(r+R)^2}{\sigma_{\max}(\Sigma)^2}\Big)D_s\exp\Big(\frac{(r+s)^2}{2\sigma_{\min}(\Sigma)^2}-\frac{\alpha\Delta}2\Big)\;\le\;\frac\rho2 .
\end{equation*}
Hence $Q_\alpha(r)\le\rho$, and Lemma~\ref{lem:simplifiedBounds}(i) gives the claim.

\emph{Ad (ii).} Let $\alpha\ge\alpha^{\mathrm S}_{\mathrm b}(\beta)$. As in (i), the first threshold gives $\tau_\alpha\le(\alpha m\sigma_{\max}(\Sigma)^2)^{-1}\le\beta^2/(2d\sigma_{\max}(\Sigma)^2)$, and the second gives $D_s\exp(s^2/(2\sigma_{\min}(\Sigma)^2)-\alpha\Delta/2)\le\beta^2/(2d\sigma_{\max}(\Sigma)^2)$. Inserting both into \eqref{eq:simplifiedBounds:bias} gives $\|\mathcal{T}_\alpha(x_\ast)-x_\ast\|_2^2\le\beta^2$.

\emph{Ad (iii).} Let $\alpha\ge\alpha^{\mathrm S}_0(r,\rho)$. By (i) and (ii) with $\beta=(1-\rho)r$, for $x\in\ball_r(x_\ast)$,
\begin{equation*}
    \big\|\mathcal{T}_\alpha(x)-x_\ast\big\|_2
    \;\le\;\big\|\mathcal{T}_\alpha(x)-\mathcal{T}_\alpha(x_\ast)\big\|_2+\big\|\mathcal{T}_\alpha(x_\ast)-x_\ast\big\|_2
    \;\le\;\rho r+(1-\rho)r=r .
\end{equation*}
The remaining claims follow from Banach's fixed-point theorem on the compact set $\ball_r(x_\ast)$ and induction, exactly as in the proof of Theorem~\ref{thm:meanFieldConvergence}.

\emph{Ad (iv).} Let $\beta_r:=\min\{\beta,(1-\rho)r\}$. Since $\alpha^{\mathrm S}_{\mathrm b}(\beta_r)$ equals $\alpha^{\mathrm S}_{\mathrm b}(\beta)$ or $\alpha^{\mathrm S}_{\mathrm b}((1-\rho)r)\le\alpha^{\mathrm S}_0(r,\rho)$, (ii) gives $\|\mathcal{T}_\alpha(x_\ast)-x_\ast\|_2\le\beta_r\le\beta$. Together with (i), the triangle inequality as in (iii) yields \eqref{eq:simplifiedMeanField:affine}. With $e_k:=\|x_k-x_\ast\|_2$, we get $e_{k+1}\le\beta+\rho e_k$, since $x_k\in\ball_r(x_\ast)$ by (iii), and induction gives $e_k\le\rho^ke_0+\beta\sum_{j<k}\rho^j\le\rho^ke_0+\frac{\beta}{1-\rho}$. Applying this with $x_0=m_\alpha$, for which $x_k=m_\alpha$ for all $k$, and letting $k\to\infty$ gives the bound on $m_\alpha$.
\end{proof}
The next result summarises the semi-global convergence behaviour under Assumption \ref{ass:simplified} and is a more detailed version of Corollary \ref{cor:meanFieldSimplified} from the main text.
\begin{corollary}[Explicit accuracy of mean-field MPPI]\label{cor:simplifiedAccuracy}
Let $f$ satisfy (SF1) and (SF2) of Assumption~\ref{ass:simplified}, let $r>0$ and $\epsilon>0$, and set $K_0(r,\epsilon):=\max\{0,\lceil\log_2(2r/\epsilon)\rceil\}$. Then for all
\begin{equation}\label{eq:simplifiedAccuracy:alpha}
    \alpha\;\ge\;\alpha^{\mathrm S}_{\mathrm{mf}}(r,\epsilon):=\max\Big\{\alpha^{\mathrm S}_{\mathrm c}\big(r,\tfrac12\big),\ \alpha^{\mathrm S}_{\mathrm b}\big(\tfrac r2\big),\ \alpha^{\mathrm S}_{\mathrm b}\big(\tfrac\epsilon4\big)\Big\}
\end{equation}
and all $x_0\in\ball_r(x_\ast)$, the mean-field iterates satisfy $\|x_k-x_\ast\|_2\le\epsilon$ for all $k\ge K_0(r,\epsilon)$.
\end{corollary}
\begin{proof}
The right-hand side of \eqref{eq:simplifiedAccuracy:alpha} equals $\max\{\alpha^{\mathrm S}_0(r,\frac12),\alpha^{\mathrm S}_{\mathrm b}(\frac\epsilon4)\}$. Proposition~\ref{prop:simplifiedMeanField}(iv) with $\rho=\frac12$ and $\beta=\frac\epsilon4$ gives $\|x_k-x_\ast\|_2\le2^{-k}r+\frac\epsilon2\le\epsilon$ for $k\ge K_0(r,\epsilon)$, since $2^k\ge2r/\epsilon$ if $K_0(r,\epsilon)>0$, and $r\le\epsilon/2$ otherwise.
\end{proof}

\subsubsection{Finite-sample Results}
We now investigate the behaviour of MPPI in a finite-sample setting under Assumption \ref{ass:simplified}.
Recall that the starting point of our finite-sample analysis is Lemma \ref{lem:empiricalTransitionMapConcentration}, and that the main quantities are the total variance of the Gibbs tilt and the expected normalised weight.
The next result bounds these quantities in our simplified setting.
\begin{lemma}[Explicit variance and effective sample size]\label{lem:simplifiedSNIS}
Let $f$ satisfy (SF1) and (SF2) of Assumption~\ref{ass:simplified}, let $\alpha>0$, $r>0$ and $x\in\ball_r(x_\ast)$.
\begin{enumerate}
    \item[(i)] The total variance $V_\alpha(x)$ from \eqref{eq:empiricalTransitionMapConcentration:variance} satisfies
    \begin{equation}\label{eq:simplifiedSNIS:variance}
        V_\alpha(x)\;\le\;d\,\sigma_{\max}(\Sigma)^2\tau_\alpha+\big(d\,\sigma_{\max}(\Sigma)^2+(r+R)^2\big)\,D_s\exp\Big(\frac{(r+s)^2}{2\sigma_{\min}(\Sigma)^2}-\frac{\alpha\Delta}2\Big).
    \end{equation}
    \item[(ii)] If $f$ also satisfies (SF3) and
    \begin{equation}\label{eq:simplifiedSNIS:alphaZeta}
        \alpha\;\ge\;\alpha^{\mathrm S}_\zeta(r):=\frac1M\max\Big\{\frac{(4d+8)\ln2}{R^2},\ \frac{2r}{R\,\sigma_{\min}(\Sigma)^2}\Big\},
    \end{equation}
    then the expected normalised weight $\zeta_\alpha(x)$ from \eqref{eq:snis:quantities} satisfies
    \begin{equation}\label{eq:simplifiedSNIS:zeta}
    \begin{aligned}
        \frac{1}{\zeta_\alpha(x)}
        &\;\le\;2\det\big(\Ident_d+\alpha M\Sigma\big)^{1/2}\exp\Big(\frac12(x-x_\ast)^\top\big(\Sigma+(\alpha M)^{-1}\Ident_d\big)^{-1}(x-x_\ast)\Big)\\
        &\;\le\;2\big(1+\alpha M\sigma_{\max}(\Sigma)^2\big)^{d/2}\exp\Big(\frac{r^2}{2\sigma_{\min}(\Sigma)^2}\Big).
    \end{aligned}
    \end{equation}
\end{enumerate}
\end{lemma}
\begin{proof}
Without loss of generality, $x_\ast=0$ and $f_\ast=0$, and write $B:=\ball_R(0)$, $P:=P_x$ and $\nu:=\Norm(x,\Sigma)$.

\emph{Ad (i).} Let $h$ be the mean of $P^{\mathrm{in}}_x$, so that $\|x-h\|_2\le r+R$ as in the proof of Lemma~\ref{lem:simplifiedBounds}. Since the mean minimises the expected squared distance,
\begin{equation*}
    V_\alpha(x)\;\le\;\E_{Y\distr P}\big[\|Y-h\|_2^2\big]
    =P(B)\,\trace\,\mathrm{Cov}(P^{\mathrm{in}}_x)+\int_{\R^d\setminus B}\|y-h\|_2^2\,P(\mathrm dy).
\end{equation*}
By Lemma~\ref{lem:simplifiedTilt}(ii), the first term is at most $d\sigma_{\max}(\Sigma)^2\tau_\alpha$. By Lemma~\ref{lem:simplifiedTilt}(i), the second term is at most $q_\alpha(r)\,\E_{X\distr\nu}[\|X-h\|_2^2]=q_\alpha(r)\big(\trace(\Sigma)+\|x-h\|_2^2\big)\le q_\alpha(r)\big(d\sigma_{\max}(\Sigma)^2+(r+R)^2\big)$.

\emph{Ad (ii).} Write $a:=\alpha M$, and for $b>0$ let $G_b:=\E_{X\distr\nu}\big[e^{-\frac b2\|X\|_2^2}\big]$. Completing the square gives the Gaussian integral
\begin{equation}\label{eq:simplifiedSNIS:gaussInt}
    G_b=\det\big(\Ident_d+b\Sigma\big)^{-1/2}\exp\Big(-\frac12x^\top\big(\Sigma+b^{-1}\Ident_d\big)^{-1}x\Big).
\end{equation}
By (SF3), $e^{-\alpha f(y)}\ge e^{-\frac a2\|y\|_2^2}$ for $y\in B$, and $e^{-\alpha f}\ge0$ on $\R^d\setminus B$. Since $f_\ast=0$, this gives
\begin{equation*}
    \zeta_\alpha(x)\;\ge\;\E_{X\distr\nu}\big[e^{-\frac a2\|X\|_2^2}\mathbf 1_B(X)\big]
    \;=\;G_a-\E_{X\distr\nu}\big[e^{-\frac a2\|X\|_2^2}\mathbf 1_{\R^d\setminus B}(X)\big].
\end{equation*}
If $\|y\|_2>R$, then $e^{-\frac a2\|y\|_2^2}\le e^{-aR^2/4}e^{-\frac a4\|y\|_2^2}$, so the last expectation is at most $e^{-aR^2/4}G_{a/2}$.
We compare $G_{a/2}$ with $G_a$ by means of \eqref{eq:simplifiedSNIS:gaussInt}. Let $\sigma_1^2,\dots,\sigma_d^2$ be the eigenvalues of $\Sigma$. Then $\det(\Ident_d+a\Sigma)/\det(\Ident_d+\frac a2\Sigma)=\prod_i(1+a\sigma_i^2)/(1+\frac a2\sigma_i^2)\le2^d$. Moreover, $(\Sigma+a^{-1}\Ident_d)^{-1}-(\Sigma+2a^{-1}\Ident_d)^{-1}$ has the eigenvalues
\begin{equation*}
    \frac{1}{\sigma_i^2+a^{-1}}-\frac{1}{\sigma_i^2+2a^{-1}}=\frac{a^{-1}}{(\sigma_i^2+a^{-1})(\sigma_i^2+2a^{-1})}\;\le\;\frac{1}{a\,\sigma_{\min}(\Sigma)^4},
\end{equation*}
and $\|x\|_2\le r$. Hence $G_{a/2}\le2^{d/2}\exp\big(r^2/(2a\sigma_{\min}(\Sigma)^4)\big)\,G_a$, and
\begin{equation*}
    \zeta_\alpha(x)\;\ge\;\Big(1-2^{d/2}\exp\Big(-\frac{aR^2}4+\frac{r^2}{2a\sigma_{\min}(\Sigma)^4}\Big)\Big)\,G_a .
\end{equation*}
By \eqref{eq:simplifiedSNIS:alphaZeta}, $aR^2/8\ge(\frac d2+1)\ln2$ and $a^2\ge4r^2/(R^2\sigma_{\min}(\Sigma)^4)$, i.e.\ $aR^2/8\ge r^2/(2a\sigma_{\min}(\Sigma)^4)$. Writing $-\frac{aR^2}4=-\frac{aR^2}8-\frac{aR^2}8$, the subtracted term is therefore at most $2^{d/2}\cdot2^{-d/2-1}=\frac12$, and $\zeta_\alpha(x)\ge\frac12G_a$, which together with \eqref{eq:simplifiedSNIS:gaussInt} is the first inequality in \eqref{eq:simplifiedSNIS:zeta}. The second one follows from $\det(\Ident_d+a\Sigma)\le(1+a\sigma_{\max}(\Sigma)^2)^d$ and $(\Sigma+a^{-1}\Ident_d)^{-1}\preceq(\sigma_{\min}(\Sigma)^2+a^{-1})^{-1}\Ident_d\preceq\sigma_{\min}(\Sigma)^{-2}\Ident_d$.
\end{proof}
For $r>0$ and $\beta>0$, define
\begin{equation}\label{eq:simplified:alphav}
    \alpha^{\mathrm S}_{\mathrm v}(r,\beta):=\max\Big\{\frac{4d}{m\,\beta^2},\ \frac2\Delta\Big[\frac{(r+s)^2}{2\sigma_{\min}(\Sigma)^2}+\log\frac{4\big(d\,\sigma_{\max}(\Sigma)^2+(r+R)^2\big)D_s}{\beta^2}\Big]\Big\},
\end{equation}
which ensures $V_\alpha(x)\le\beta^2/2$ for all $x\in\ball_r(x_\ast)$, cf.\ the proof below.
We can now formulate a result on the finite-sample behaviour for finite horizons under our simplified assumption.
\begin{proposition}[Explicit finite-sample convergence]\label{prop:simplifiedFiniteSample}
Let $f$ satisfy Assumption~\ref{ass:simplified}, including (SF3), let $r>0$, $\rho\in(0,1)$ and $\beta>0$ with $2\beta\le(1-\rho)r$, and let
\begin{equation}\label{eq:simplifiedFiniteSample:alpha}
    \alpha\;\ge\;\max\Big\{\alpha^{\mathrm S}_0(r,\rho),\ \alpha^{\mathrm S}_{\mathrm b}(\beta),\ \alpha^{\mathrm S}_{\mathrm v}(r,\beta),\ \alpha^{\mathrm S}_\zeta(r)\Big\}.
\end{equation}
Then for all $K\in\Np$, $\delta\in(0,1)$, $N\in\Np$ with
\begin{equation}\label{eq:simplifiedFiniteSample:N}
    N\;\ge\;\frac{12\,K}{\delta}\,\big(1+\alpha M\sigma_{\max}(\Sigma)^2\big)^{d/2}\exp\Big(\frac{r^2}{2\sigma_{\min}(\Sigma)^2}\Big),
\end{equation}
and all $\hat x_0\in\ball_r(x_\ast)$, the MPPI iterates \eqref{eq:mppiIteration} satisfy
\begin{equation}\label{eq:simplifiedFiniteSample:claim}
    \Pb\Big(\hat x_k\in\ball_r(x_\ast)\ \text{and}\ \|\hat x_k-x_\ast\|_2\le\rho^k\,\|\hat x_0-x_\ast\|_2+\frac{2\beta}{1-\rho}\ \text{ for all }k=0,\dots,K\Big)\;\ge\;1-\delta .
\end{equation}
\end{proposition}
\begin{proof}
Fix $\alpha$, $K$, $\delta$, $N$ and $\hat x_0$ as in the statement, and write $e_k:=\|\hat x_k-x_\ast\|_2$.

\emph{Step 1: Variance and effective sample size.} Let $x\in\ball_r(x_\ast)$. The first threshold in \eqref{eq:simplified:alphav} gives $d\sigma_{\max}(\Sigma)^2\tau_\alpha\le d/(\alpha m)\le\beta^2/4$, and the second one bounds the second term in \eqref{eq:simplifiedSNIS:variance} by $\beta^2/4$. Hence Lemma~\ref{lem:simplifiedSNIS}(i) gives $V_\alpha(x)\le\beta^2/2$. Moreover, Lemma~\ref{lem:simplifiedSNIS}(ii) and \eqref{eq:simplifiedFiniteSample:N} give $\frac{6}{N\zeta_\alpha(x)}\le\frac\delta K$.

\emph{Step 2: One sampled step.} For $x\in\R^d$, let $\psi(x)$ be the probability that the empirical transition map at $x$, computed from $N$ independent samples of $\Norm(x,\Sigma)$, deviates from $\mathcal{T}_\alpha(x)$ by more than $\beta$. For $x\in\ball_r(x_\ast)$, Lemma~\ref{lem:empiricalTransitionMapConcentration}(i) with $t=\beta$ and Step~1 give
\begin{equation*}
    \psi(x)\;\le\;\frac4N\,\frac{1}{\zeta_\alpha(x)}\Big(1+\frac{V_\alpha(x)}{\beta^2}\Big)\;\le\;\frac{6}{N\zeta_\alpha(x)}\;\le\;\frac\delta K .
\end{equation*}

\emph{Step 3: Deterministic progress on good events.} For $k\in\Nnn$, let $A_k:=\{\|\hat x_{k+1}-\mathcal{T}_\alpha(\hat x_k)\|_2\le\beta\}$. Suppose that $\hat x_k\in\ball_r(x_\ast)$ and that $A_k$ occurs. Since $\alpha\ge\max\{\alpha^{\mathrm S}_0(r,\rho),\alpha^{\mathrm S}_{\mathrm b}(\beta)\}$, Proposition~\ref{prop:simplifiedMeanField}(iv) gives
\begin{equation*}
    e_{k+1}\;\le\;\big\|\mathcal{T}_\alpha(\hat x_k)-x_\ast\big\|_2+\beta\;\le\;2\beta+\rho\,e_k\;\le\;2\beta+\rho r\;\le\;r,
\end{equation*}
where the last step uses $2\beta\le(1-\rho)r$. By induction, on $\bigcap_{j=0}^{k-1}A_j$ we have $\hat x_j\in\ball_r(x_\ast)$ for $j=0,\dots,k$, and $e_k\le\rho^ke_0+2\beta\sum_{j<k}\rho^j\le\rho^ke_0+\frac{2\beta}{1-\rho}$. Hence the event in \eqref{eq:simplifiedFiniteSample:claim} contains $\bigcap_{k=0}^{K-1}A_k$.

\emph{Step 4: Probability of the good events.} As in Steps~0 and~4 of the proof of Theorem~\ref{thm:finiteSample}, $\hat x_k$ is independent of the samples used in step $k$, so Fubini's theorem and Step~2 give $\Pb(A_k^c\cap\{\hat x_k\in\ball_r(x_\ast)\})=\E[\mathbf 1_{\ball_r(x_\ast)}(\hat x_k)\psi(\hat x_k)]\le\delta/K$. Decomposing according to the first failure time and using Step~3,
\begin{equation*}
    \Pb\Big(\bigcup_{k=0}^{K-1}A_k^c\Big)
    \;\le\;\sum_{k=0}^{K-1}\Pb\Big(A_k^c\cap\big\{\hat x_k\in\ball_r(x_\ast)\big\}\Big)\;\le\;\delta ,
\end{equation*}
which proves \eqref{eq:simplifiedFiniteSample:claim}.
\end{proof}
Finally, we can use the preceding result to get an explicit bound on the sample complexity to achieve a given accuracy with high probability.
\begin{corollary}[Explicit sample complexity]\label{cor:simplifiedSampleComplexity}
Let $f$ satisfy Assumption~\ref{ass:simplified}, including (SF3), let $r>0$ and $\epsilon\in(0,2r]$, and set $K_0(r,\epsilon):=\max\{0,\lceil\log_2(2r/\epsilon)\rceil\}$. Let
\begin{equation}\label{eq:simplifiedSampleComplexity:alpha}
    \alpha\;\ge\;\alpha^{\mathrm S}_{\mathrm{fs}}(r,\epsilon):=\max\Big\{\alpha^{\mathrm S}_{\mathrm c}\big(r,\tfrac12\big),\ \alpha^{\mathrm S}_{\mathrm b}\big(\tfrac r2\big),\ \alpha^{\mathrm S}_{\mathrm b}\big(\tfrac\epsilon8\big),\ \alpha^{\mathrm S}_{\mathrm v}\big(r,\tfrac\epsilon8\big),\ \alpha^{\mathrm S}_\zeta(r)\Big\},
\end{equation}
and let $K\in\Np$, $\delta\in(0,1)$ and $N$ as in \eqref{eq:simplifiedFiniteSample:N}. Then for all $\hat x_0\in\ball_r(x_\ast)$,
\begin{equation*}
    \Pb\Big(\|\hat x_k-x_\ast\|_2\le\epsilon\ \text{ for all }k\text{ with }K_0(r,\epsilon)\le k\le K\Big)\;\ge\;1-\delta .
\end{equation*}
\end{corollary}
\begin{proof}
Apply Proposition~\ref{prop:simplifiedFiniteSample} with $\rho=\frac12$ and $\beta=\frac\epsilon8$. Then $2\beta=\frac\epsilon4\le\frac r2=(1-\rho)r$, the threshold on $\alpha$ is \eqref{eq:simplifiedFiniteSample:alpha}, and $\frac{2\beta}{1-\rho}=\frac\epsilon2$. For $k\ge K_0(r,\epsilon)$, we have $2^{-k}r\le\frac\epsilon2$, and hence $\|\hat x_k-x_\ast\|_2\le\epsilon$ on the event in \eqref{eq:simplifiedFiniteSample:claim}.
\end{proof}

\subsubsection{Explicit Sample Complexity Bounds}
Finally, a significant advantage of working with the simplified Assumption \ref{ass:simplified} is the possibility to derive explicit sample complexity bounds.
The following result describes the necessary estimates for this.
\begin{lemma}[Order of the explicit thresholds]\label{lem:simplifiedThresholdOrders}
Let $f$ satisfy (SF1) and (SF2) of Assumption~\ref{ass:simplified}, and let $r>0$ and $\epsilon\in(0,2r]$. With $\log_+:=\max\{0,\log\}$,
\begin{equation}\label{eq:simplifiedThresholdOrders:mf}
\begin{aligned}
    \alpha^{\mathrm S}_{\mathrm{mf}}(r,\epsilon)\;\le{}&\frac{32d}{m\,\epsilon^2}+\frac{4}{m\,\sigma_{\min}(\Sigma)^2}+\frac{4(r+R)^2}{\Delta\,\sigma_{\min}(\Sigma)^2}\\
    &+\frac2\Delta\Big[d\log_+\frac{(d+2)\sigma_{\max}(\Sigma)^2}{s^2}+\log(4\kappa_\Sigma)+\log_+\frac{32d\,\sigma_{\max}(\Sigma)^2}{\epsilon^2}\Big].
\end{aligned}
\end{equation}
If $f$ also satisfies (SF3), then
\begin{equation}\label{eq:simplifiedThresholdOrders:fs}
\begin{aligned}
    \alpha^{\mathrm S}_{\mathrm{fs}}(r,\epsilon)\;\le{}&\frac{384d}{m\,\epsilon^2}+\frac{4}{m\,\sigma_{\min}(\Sigma)^2}+\frac{8(r+R)^2}{\Delta\,\sigma_{\min}(\Sigma)^2}\\
    &+\frac1\Delta\Big[3d\log_+\frac{(d+2)\sigma_{\max}(\Sigma)^2}{s^2}+2\log(4\kappa_\Sigma)+4\log_+\frac{256d\,\sigma_{\max}(\Sigma)^2}{\epsilon^2}+5d\Big].
\end{aligned}
\end{equation}
In particular, both thresholds are of the order \eqref{eq:meanFieldSimplified:order}.
\end{lemma}
\begin{proof}
\emph{Preliminaries.} For $t\ge0$, we have $\Gamma(t+1)\le(t+1)^t$: for $t\in[0,1]$, $\Gamma(t+1)\le1\le(t+1)^t$, and for $t>1$ with $n:=\lceil t\rceil-1\ge1$, $\Gamma(t+1)=t(t-1)\cdots(t-n+1)\Gamma(t-n+1)\le t^n\le(t+1)^t$, since $t-n+1\in(1,2]$. With $t=\frac d2$, this gives $D_s\le((d+2)\sigma_{\max}(\Sigma)^2/s^2)^{d/2}$, i.e.\ $\log_+D_s\le\frac d2\log_+\frac{(d+2)\sigma_{\max}(\Sigma)^2}{s^2}$. Moreover, both entries of $\alpha^{\mathrm S}_{\mathrm b}(\beta)$ are nonincreasing in $\beta$, and $\log(1+y)\le y$ for $y\ge0$. Finally, $\max$ of finitely many numbers is at most the sum of upper bounds for their positive parts.

\emph{Mean-field threshold.} Since $\epsilon\le2r$, we have $\frac r2\ge\frac\epsilon4$, so $\alpha^{\mathrm S}_{\mathrm b}(\frac r2)\le\alpha^{\mathrm S}_{\mathrm b}(\frac\epsilon4)$ and $\alpha^{\mathrm S}_{\mathrm{mf}}(r,\epsilon)=\max\{\alpha^{\mathrm S}_{\mathrm c}(r,\frac12),\alpha^{\mathrm S}_{\mathrm b}(\frac\epsilon4)\}$. The entries of $\alpha^{\mathrm S}_{\mathrm c}(r,\frac12)$ are $\frac{4}{m\sigma_{\min}(\Sigma)^2}$ and
\begin{align*}
    &\frac{(r+s)^2}{\Delta\sigma_{\min}(\Sigma)^2}+\frac2\Delta\Big[\log(4\kappa_\Sigma)+\log D_s+\log\Big(1+\frac{(r+R)^2}{\sigma_{\max}(\Sigma)^2}\Big)\Big]\\
    &\;\le\;\frac{3(r+R)^2}{\Delta\sigma_{\min}(\Sigma)^2}+\frac2\Delta\big[\log(4\kappa_\Sigma)+\log_+D_s\big],
\end{align*}
where we used $s\le R$ and $\sigma_{\min}(\Sigma)\le\sigma_{\max}(\Sigma)$. The entries of $\alpha^{\mathrm S}_{\mathrm b}(\frac\epsilon4)$ are $\frac{32d}{m\epsilon^2}$ and
\begin{align*}
    &\frac{s^2}{\Delta\sigma_{\min}(\Sigma)^2}+\frac2\Delta\Big[\log D_s+\log\frac{32d\,\sigma_{\max}(\Sigma)^2}{\epsilon^2}\Big]\\
    &\;\le\;\frac{(r+R)^2}{\Delta\sigma_{\min}(\Sigma)^2}+\frac2\Delta\Big[\log_+D_s+\log_+\frac{32d\,\sigma_{\max}(\Sigma)^2}{\epsilon^2}\Big].
\end{align*}
Adding these bounds and using the bound on $\log_+D_s$ gives \eqref{eq:simplifiedThresholdOrders:mf}.

\emph{Finite-sample threshold.} As above, $\alpha^{\mathrm S}_{\mathrm b}(\frac r2)\le\alpha^{\mathrm S}_{\mathrm b}(\frac\epsilon8)$. The entries of $\alpha^{\mathrm S}_{\mathrm c}(r,\frac12)$ are bounded as above. The entries of $\alpha^{\mathrm S}_{\mathrm b}(\frac\epsilon8)$ are $\frac{128d}{m\epsilon^2}$ and at most $\frac{(r+R)^2}{\Delta\sigma_{\min}(\Sigma)^2}+\frac2\Delta\big[\log_+D_s+\log_+\frac{128d\sigma_{\max}(\Sigma)^2}{\epsilon^2}\big]$. The entries of $\alpha^{\mathrm S}_{\mathrm v}(r,\frac\epsilon8)$ are $\frac{256d}{m\epsilon^2}$ and, writing $d\sigma_{\max}(\Sigma)^2+(r+R)^2=d\sigma_{\max}(\Sigma)^2\big(1+\frac{(r+R)^2}{d\sigma_{\max}(\Sigma)^2}\big)$,
\begin{align*}
    &\frac{(r+s)^2}{\Delta\sigma_{\min}(\Sigma)^2}+\frac2\Delta\Big[\log D_s+\log\frac{256\big(d\sigma_{\max}(\Sigma)^2+(r+R)^2\big)}{\epsilon^2}\Big]\\
    &\;\le\;\frac{3(r+R)^2}{\Delta\sigma_{\min}(\Sigma)^2}+\frac2\Delta\Big[\log_+D_s+\log_+\frac{256d\,\sigma_{\max}(\Sigma)^2}{\epsilon^2}\Big].
\end{align*}
For $\alpha^{\mathrm S}_\zeta(r)$, note that $\Delta\le\frac M2R^2$: for $\|y-x_\ast\|_2=R$, the points $x_\ast+(1+\frac1n)(y-x_\ast)$ lie outside $\ball_R(x_\ast)$, so (SF2) and continuity give $f(y)-f_\ast\ge\Delta$, while (SF3) gives $f(y)-f_\ast\le\frac M2R^2$. Hence $\frac{(4d+8)\ln2}{MR^2}\le\frac{(2d+4)\ln2}{\Delta}\le\frac{5d}{\Delta}$ and $\frac{2r}{MR\sigma_{\min}(\Sigma)^2}\le\frac{rR}{\Delta\sigma_{\min}(\Sigma)^2}\le\frac{(r+R)^2}{4\Delta\sigma_{\min}(\Sigma)^2}$. Adding all bounds, using $\log_+\frac{128d\sigma_{\max}(\Sigma)^2}{\epsilon^2}\le\log_+\frac{256d\sigma_{\max}(\Sigma)^2}{\epsilon^2}$ and the bound on $\log_+D_s$, gives \eqref{eq:simplifiedThresholdOrders:fs}.

\emph{Order.} In \eqref{eq:simplifiedThresholdOrders:mf} and \eqref{eq:simplifiedThresholdOrders:fs}, the terms in square brackets are at most $d$ times a factor that is polylogarithmic in $d$, $\kappa_\Sigma$, $\sigma_{\max}(\Sigma)/\epsilon$ and $\sigma_{\max}(\Sigma)/s$, since $d\ge1$. This gives \eqref{eq:meanFieldSimplified:order}.
\end{proof}
We can finally provide the missing proofs of the simplified, concrete results from the main text.
\begin{proof}[Proof of Corollary~\ref{cor:meanFieldSimplified}]
With $\alpha^{\mathrm S}_{\mathrm{mf}}(r,\epsilon)$ from \eqref{eq:simplifiedAccuracy:alpha}, the claim on the iterates is Corollary~\ref{cor:simplifiedAccuracy}, and the order \eqref{eq:meanFieldSimplified:order} is Lemma~\ref{lem:simplifiedThresholdOrders}. The radius $s$ is the one of Lemma~\ref{lem:simplifiedImpliesObjective}(iv), which can be any radius with $f\le f_\ast+\frac\Delta2$ on $\ball_s(x_\ast)$.
\end{proof}
\begin{proof}[Proof of Corollary~\ref{cor:finiteSampleSimplified}]
With $\alpha^{\mathrm S}_{\mathrm{fs}}(r,\epsilon)$ from \eqref{eq:simplifiedSampleComplexity:alpha}, the claim is Corollary~\ref{cor:simplifiedSampleComplexity}, since \eqref{eq:finiteSampleSimplified:N} is \eqref{eq:simplifiedFiniteSample:N}, and the order is Lemma~\ref{lem:simplifiedThresholdOrders}.
\end{proof}

\subsubsection{Quadratic objectives}
Finally, we consider a particularly simple class of functions fulfilling Assumption \ref{ass:simplified}, the quadratic functions.
\begin{lemma}[Quadratic objectives]\label{lem:quadraticSimplified}
Let $f(x)=f_\ast+\frac12(x-x_\ast)^\top H(x-x_\ast)$ with $f_\ast\in\R$, $x_\ast\in\R^d$ and $H\in\R^{d\times d}$ symmetric positive definite, and let $\lambda_-:=\lambda_{\min}(H)$ and $\lambda_+:=\lambda_{\max}(H)$.
\begin{enumerate}
    \item[(i)] For every $R>0$, $f$ satisfies Assumption~\ref{ass:simplified} with $m=\lambda_-$, $M=\lambda_+$ and $\Delta=\frac{\lambda_-}2R^2$, and $s=R\sqrt{\lambda_-/(2\lambda_+)}$ is admissible in Lemma~\ref{lem:simplifiedImpliesObjective}(iv).
    \item[(ii)] Let $r>0$ and $\epsilon\in(0,2r]$. For all $\alpha>\frac{32d}{\lambda_-\epsilon^2}+\frac{12}{\lambda_-\sigma_{\min}(\Sigma)^2}$ and $x_0\in\ball_r(x_\ast)$, the mean-field iterates satisfy $\|x_k-x_\ast\|_2\le\epsilon$ for all $k\ge K_0(r,\epsilon)$.
    \item[(iii)] Let $r>0$ and $\epsilon\in(0,2r]$. For all $\alpha>\frac{384d}{\lambda_-\epsilon^2}+\frac{20}{\lambda_-\sigma_{\min}(\Sigma)^2}$, $K\in\Np$, $\delta\in(0,1)$, $N\in\Np$ with
    \begin{equation*}
        N\;\ge\;\frac{12K}{\delta}\big(1+\alpha\lambda_+\sigma_{\max}(\Sigma)^2\big)^{d/2}\exp\Big(\frac{r^2}{2\sigma_{\min}(\Sigma)^2}\Big),
    \end{equation*}
    and $\hat x_0\in\ball_r(x_\ast)$, the MPPI iterates satisfy
    \begin{equation*}
        \Pb\Big(\|\hat x_k-x_\ast\|_2\le\epsilon\ \text{ for all }k\text{ with }K_0(r,\epsilon)\le k\le K\Big)\;\ge\;1-\delta .
    \end{equation*}
\end{enumerate}
\end{lemma}
\begin{proof}
\emph{Ad (i).} $f$ is continuous and attains its global minimum $f_\ast$ at $x_\ast$. Since $H-\lambda_-\Ident_d\succeq0$, the function $y\mapsto f(y)-\frac{\lambda_-}2\|y-x_\ast\|_2^2=f_\ast+\frac12(y-x_\ast)^\top(H-\lambda_-\Ident_d)(y-x_\ast)$ is convex on $\R^d$, which gives (SF1) with $m=\lambda_-$. For $\|y-x_\ast\|_2>R$, we have $f(y)-f_\ast\ge\frac{\lambda_-}2\|y-x_\ast\|_2^2\ge\frac{\lambda_-}2R^2$, which is (SF2) with $\Delta=\frac{\lambda_-}2R^2$. Moreover, $f(y)-f_\ast\le\frac{\lambda_+}2\|y-x_\ast\|_2^2$ for all $y$, which is (SF3) with $M=\lambda_+\ge m$. Finally, $\sqrt{\Delta/M}=R\sqrt{\lambda_-/(2\lambda_+)}\le R$, so Lemma~\ref{lem:simplifiedImpliesObjective}(iv) shows that $s=R\sqrt{\lambda_-/(2\lambda_+)}$ is admissible.

\emph{Ad (ii).} For $R>0$, insert $m=\lambda_-$, $\Delta=\frac{\lambda_-}2R^2$ and $s=R\sqrt{\lambda_-/(2\lambda_+)}$ from (i) into \eqref{eq:simplifiedThresholdOrders:mf}, and denote the resulting right-hand side by $B_{\mathrm{mf}}(R)$. Then $\frac{4(r+R)^2}{\Delta\sigma_{\min}(\Sigma)^2}=\frac{8(r+R)^2}{\lambda_-R^2\sigma_{\min}(\Sigma)^2}$, and the term in square brackets is multiplied by $\frac2\Delta=\frac{4}{\lambda_-R^2}$ and is nonincreasing in $R$, since $s$ is increasing in $R$. Hence $B_{\mathrm{mf}}(R)$ is nonincreasing in $R$, and
\begin{equation*}
    \lim_{R\to\infty}B_{\mathrm{mf}}(R)=\frac{32d}{\lambda_-\epsilon^2}+\frac{4}{\lambda_-\sigma_{\min}(\Sigma)^2}+\frac{8}{\lambda_-\sigma_{\min}(\Sigma)^2}
    =\frac{32d}{\lambda_-\epsilon^2}+\frac{12}{\lambda_-\sigma_{\min}(\Sigma)^2}.
\end{equation*}
Let $\alpha$ be strictly larger than this limit. Then there is $R>0$ with $\alpha\ge B_{\mathrm{mf}}(R)$, and Lemma~\ref{lem:simplifiedThresholdOrders}, applied with this $R$, gives $\alpha\ge\alpha^{\mathrm S}_{\mathrm{mf}}(r,\epsilon)$. Corollary~\ref{cor:simplifiedAccuracy} yields the claim.

\emph{Ad (iii).} Analogously, insert (i) into \eqref{eq:simplifiedThresholdOrders:fs} and denote the right-hand side by $B_{\mathrm{fs}}(R)$. Then $\frac{8(r+R)^2}{\Delta\sigma_{\min}(\Sigma)^2}=\frac{16(r+R)^2}{\lambda_-R^2\sigma_{\min}(\Sigma)^2}$, the term in square brackets is multiplied by $\frac1\Delta=\frac{2}{\lambda_-R^2}$ and is nonincreasing in $R$, so $B_{\mathrm{fs}}(R)$ is nonincreasing in $R$ and converges to $\frac{384d}{\lambda_-\epsilon^2}+\frac{4+16}{\lambda_-\sigma_{\min}(\Sigma)^2}$ as $R\to\infty$. For $\alpha$ strictly larger than this limit, there is $R>0$ with $\alpha\ge B_{\mathrm{fs}}(R)\ge\alpha^{\mathrm S}_{\mathrm{fs}}(r,\epsilon)$ by Lemma~\ref{lem:simplifiedThresholdOrders}, and Corollary~\ref{cor:simplifiedSampleComplexity} yields the claim, since the sample size condition \eqref{eq:simplifiedFiniteSample:N} with $M=\lambda_+$ does not depend on $R$.
\end{proof}

\subsection{Details for Section~\ref{sec:oc}}\label{app:proofsOC}
We recall how the cost \eqref{eq:control:cost} is rewritten as the quadratic function \eqref{eq:control:quadratic}.
We identify a control sequence $a=(a_0,\dots,a_{T-1})\in(\R^p)^T$ with the stacked vector $(a_0^\top,\dots,a_{T-1}^\top)^\top\in\R^{pT}$, and the state sequence $(s_0,\dots,s_{T-1})$ with the stacked vector $\mathbf s:=(s_0^\top,\dots,s_{T-1}^\top)^\top\in\R^{nT}$. Note that the final state $s_T$ does not enter \eqref{eq:control:cost}.

\emph{States as an affine function of the controls.}
Iterating \eqref{eq:control:lti} gives the variation-of-constants formula
\begin{equation*}
    s_t=A^ts_0+\sum_{k=0}^{t-1}A^{t-1-k}Ba_k,\qquad t=0,\dots,T-1 .
\end{equation*}
In stacked form, this reads $\mathbf s=\Phi_Ts_0+\Gamma_Ta$ with
\begin{equation}\label{eq:control:PhiGamma}
    \Phi_T:=\begin{pmatrix}\Ident_n\\ A\\ A^2\\ \vdots\\ A^{T-1}\end{pmatrix}\in\R^{nT\times n},
    \qquad
    \Gamma_T:=\begin{pmatrix}
        0 & 0 & \cdots & 0 & 0\\
        B & 0 & \cdots & 0 & 0\\
        AB & B & \ddots & \vdots & \vdots\\
        \vdots & \ddots & \ddots & 0 & 0\\
        A^{T-2}B & \cdots & AB & B & 0
    \end{pmatrix}\in\R^{nT\times pT},
\end{equation}
i.e.\ the $(t,k)$-th block of $\Gamma_T$, $t,k=0,\dots,T-1$, is $A^{t-1-k}B$ if $k<t$ and $0$ otherwise.

\emph{The cost as a quadratic form.}
Let $\bar Q:=\Ident_T\otimes Q\in\R^{nT\times nT}$, $\bar R:=\Ident_T\otimes R\in\R^{pT\times pT}$ and $\bar S:=\Ident_T\otimes S\in\R^{nT\times pT}$ be the block diagonal matrices with diagonal blocks $Q$, $R$ and $S$, respectively. Then \eqref{eq:control:cost} becomes
\begin{equation*}
    J_T(s_0,a)=\mathbf s^\top\bar Q\mathbf s+a^\top\bar Ra+2\,\mathbf s^\top\bar Sa .
\end{equation*}
Inserting $\mathbf s=\Phi_Ts_0+\Gamma_Ta$ and using the symmetry of $\bar Q$, we get
\begin{align*}
    \mathbf s^\top\bar Q\mathbf s&=a^\top\Gamma_T^\top\bar Q\Gamma_Ta+2a^\top\Gamma_T^\top\bar Q\Phi_Ts_0+s_0^\top\Phi_T^\top\bar Q\Phi_Ts_0,\\
    2\,\mathbf s^\top\bar Sa&=2a^\top\bar S^\top\Gamma_Ta+2a^\top\bar S^\top\Phi_Ts_0
    =a^\top\big(\Gamma_T^\top\bar S+\bar S^\top\Gamma_T\big)a+2a^\top\bar S^\top\Phi_Ts_0 .
\end{align*}
Collecting the terms that are quadratic in $a$, bilinear in $(a,s_0)$, and quadratic in $s_0$ gives \eqref{eq:control:quadratic} with
\begin{equation}\label{eq:control:HFY}
\begin{aligned}
    H_T&:=\bar R+\Gamma_T^\top\bar Q\Gamma_T+\Gamma_T^\top\bar S+\bar S^\top\Gamma_T\in\R^{pT\times pT},\\
    F_T&:=\big(\Gamma_T^\top\bar Q+\bar S^\top\big)\Phi_T\in\R^{pT\times n},
    \qquad
    Y_T:=\Phi_T^\top\bar Q\Phi_T\in\R^{n\times n}.
\end{aligned}
\end{equation}
Since $Q$ and $R$ are symmetric, so are $H_T$ and $Y_T$.

\emph{Minimiser.}
If $H_T$ is positive definite, completing the square in \eqref{eq:control:quadratic}, i.e.\ using
\begin{equation*}
    a^\top H_Ta+2a^\top F_Ts_0
    =\big(a+H_T^{-1}F_Ts_0\big)^\top H_T\big(a+H_T^{-1}F_Ts_0\big)-s_0^\top F_T^\top H_T^{-1}F_Ts_0,
\end{equation*}
gives
\begin{equation}\label{eq:control:expansion}
    J_T(s_0,a)=J_T^\ast(s_0)+\big(a-a_\ast(s_0)\big)^\top H_T\big(a-a_\ast(s_0)\big)\qquad\text{for all }a\in\R^{pT},
\end{equation}
with $a_\ast(s_0)=-H_T^{-1}F_Ts_0$ and $J_T^\ast(s_0)=s_0^\top(Y_T-F_T^\top H_T^{-1}F_T)s_0$. In particular, $a_\ast(s_0)$ is the unique solution of \eqref{eq:control:lqocp} and $J_T^\ast(s_0)$ is its optimal value. If $H_T$ is not positive definite, then there is $v\neq0$ with $v^\top H_Tv\le0$, and $\lambda\mapsto J_T(s_0,a+\lambda v)$ is a polynomial of degree at most two with nonpositive leading coefficient, so a minimiser of \eqref{eq:control:lqocp} cannot be unique.

Recall from the main text that our MPPI results depend on spectral properties of $H_T$, which are in turn determined by those of $Q,R,S$, and $A,B$.
For convenience, we provide some elementary bounds that are well-known in the control theory literature.
\begin{lemma}\label{lem:control:HTbounds}
Let $W:=\begin{pmatrix}Q&S\\S^\top&R\end{pmatrix}\in\R^{(n+p)\times(n+p)}$ be positive definite. Then $H_T$ is positive definite, and
\begin{gather*}
    \lambda_{\min}(H_T)\;\ge\;\lambda_{\min}(W),\\
    \lambda_{\max}(H_T)\;\le\;\lambda_{\max}(W)\big(1+\|\Gamma_T\|^2\big)\;\le\;\lambda_{\max}(W)\Big(1+\Big(\sum_{j=0}^{T-2}\|A^jB\|\Big)^2\Big).
\end{gather*}
In particular, if $A$ is Schur stable, then $\lambda_{\max}(H_T)/\lambda_{\min}(H_T)$ is bounded uniformly in $T$.
\end{lemma}
\begin{proof}
Let $a\in\R^{pT}$, let $s_t=s(t;0,a)$, so that $\mathbf s=\Gamma_Ta$, and let $z_t:=(s_t^\top,a_t^\top)^\top$. By \eqref{eq:control:quadratic} with $s_0=0$ and \eqref{eq:control:cost},
\begin{equation*}
    a^\top H_Ta=J_T(0,a)=\sum_{t=0}^{T-1}z_t^\top Wz_t .
\end{equation*}
Hence $a^\top H_Ta\ge\lambda_{\min}(W)\sum_{t=0}^{T-1}\big(\|s_t\|_2^2+\|a_t\|_2^2\big)\ge\lambda_{\min}(W)\|a\|_2^2$, which gives positive definiteness and the lower bound, and $a^\top H_Ta\le\lambda_{\max}(W)\big(\|a\|_2^2+\|\Gamma_Ta\|_2^2\big)$, which gives the first upper bound.
Let $L_T\in\R^{T\times T}$ be the lower shift matrix, i.e.\ $(L_T)_{t,k}=1$ if $t=k+1$ and $0$ otherwise. Then $\Gamma_T=\sum_{j=0}^{T-2}L_T^{j+1}\otimes A^jB$ by \eqref{eq:control:PhiGamma}, and since $\|L_T^{j+1}\|\le1$ and $\|X\otimes Y\|=\|X\|\,\|Y\|$, we get $\|\Gamma_T\|\le\sum_{j=0}^{T-2}\|A^jB\|$.
Finally, if $A$ is Schur stable, then there are $C_A\ge1$ and $\rho_A\in(0,1)$ with $\|A^j\|\le C_A\rho_A^j$ for all $j$, so $\sum_{j=0}^{T-2}\|A^jB\|\le C_A\|B\|/(1-\rho_A)$ for all $T$.
\end{proof}

Finally, we provide a proof of Corollary~\ref{cor:control:explicit}.
\begin{proof}[Proof of Corollary~\ref{cor:control:explicit}]
By \eqref{eq:control:expansion}, $f(a)=J_T^\ast(s_0)+\frac12(a-a_\ast)^\top(2H_T)(a-a_\ast)$, so $f$ is the quadratic objective of Lemma~\ref{lem:quadraticSimplified} with $d=pT$, $x_\ast=a_\ast$, $H=2H_T$, $\lambda_-=2\lambda_{\min}(H_T)$ and $\lambda_+=2\lambda_{\max}(H_T)$. Moreover, $\Sigma=\Ident_T\otimes\Sigma_c$ has the same eigenvalues as $\Sigma_c$, so $\sigma_{\min}(\Sigma)=\sigma_{\min}(\Sigma_c)$ and $\sigma_{\max}(\Sigma)=\sigma_{\max}(\Sigma_c)$. Since
\begin{align*}
    \frac{32d}{\lambda_-\epsilon^2}+\frac{12}{\lambda_-\sigma_{\min}(\Sigma)^2}&=\frac{16pT}{\lambda_{\min}(H_T)\epsilon^2}+\frac{6}{\lambda_{\min}(H_T)\sigma_{\min}(\Sigma_c)^2},\\
    \frac{384d}{\lambda_-\epsilon^2}+\frac{20}{\lambda_-\sigma_{\min}(\Sigma)^2}&=\frac{192pT}{\lambda_{\min}(H_T)\epsilon^2}+\frac{10}{\lambda_{\min}(H_T)\sigma_{\min}(\Sigma_c)^2},
\end{align*}
and $\alpha\lambda_+=2\alpha\lambda_{\max}(H_T)$, Lemma~\ref{lem:quadraticSimplified}(ii) gives (i) and Lemma~\ref{lem:quadraticSimplified}(iii) gives (ii).
\end{proof}

\end{document}